\documentclass[12pt,twoside]{article} 
\usepackage{amsmath,amssymb,latexsym,theorem,bbm,shapepar,natbib,epsfig}
\newfont{\msa}{msam10 scaled\magstep1}
\newfont{\ssmsa}{msam9}

\def\Sign{\mathop{\hbox{\rm sign}}}
    
\newcommand{\Var}{\mathrm{Var}}
\newcommand{\Cov}{\mathrm{Cov}}

\newcommand{\bone}{\mathbbm{1}}
\newcommand{\stoch}{\stackrel{\scriptstyle\mathsf{P}}{\longrightarrow}}
\newcommand{\distr}{\stackrel{\scriptstyle{\cal D}}{\longrightarrow}}
\newcommand{\qmean}{\stackrel{\scriptstyle\mathsf{L}_2}{\longrightarrow}}

\newcommand{\proofend}{\hfill$\square$}

\numberwithin{equation}{section}

\theorembodyfont{\em}
\newtheorem{Lem}{Lemma}[section]
\newtheorem{Thm}[Lem]{Theorem}
\newtheorem{Pro}[Lem]{Proposition}

\newtheorem{Cor}[Lem]{Corollary}
\theorembodyfont{\rm}

\title{Parameter estimation in a spatial unit root
  autoregressive
  model}

\author{{\sc S\'andor Baran}\\ 
         Faculty of Informatics, University of Debrecen\\
         Kassai \'ut 26, H--4028 Debrecen, Hungary \\[2mm]
        {\sc Gyula Pap}\\
        Bolyai Institute, University of
             Szeged \\ Aradi v\'ertan\'uk tere 1, H-6720 Szeged, 
             Hungary. 
}

\date{}

\begin{document}
\pagestyle{myheadings}

\maketitle

\begin{abstract}
Spatial unilateral autoregressive model \ $X_{k,\ell}
  =\alpha X_{k-1,\ell}+\beta X_{k,\ell-1}+\gamma X_{k-1,\ell-1}
   +\varepsilon_{k,\ell}$ \ is investigated in the unit root case,
   that is when the parameters are on the boundary of the domain of
   stability that forms a tetrahedron  with
 vertices \ $(1,1,-1)$, \ $(1,-1,1)$,\ $(-1,1,1)$ \ and \ $(-1,-1,-1)$. \
 It is shown 
   that the limiting distribution of the least squares estimator of
   the parameters is normal and the rate of convergence is \ $n$ \
   when the parameters are in the faces or on the edges of the
   tetrahedron, while on the vertices the rate is \ $n^{3/2}$.

\noindent {\bf Key words:\/} Spatial unilateral autoregressive processes, 
unit root models.
\end{abstract}

\section{Introduction}
   \label{sec:sec1}
\markboth{\rm Baran, Pap}{\rm Parameter estimation in a unit root
  AR model} 

The analysis of spatial autoregressive models is of interest 
in many different fields of science such as 
 geography, geology, biology and agriculture. A detailed discussion of these
applications is given by \citet{br2} where the  
authors considered a special case of the so called 
unilateral autoregressive model having the form 
 \begin{equation}
    \label{brmod}
  X_{k_,\ell}
  =\sum_{i=0}^{p_1}\sum_{j=0}^{p_2}\alpha_{i,j}X_{k-i,\ell-j}
   +\varepsilon_{k,\ell}, 
  \qquad\alpha_{0,0}=0.
\end{equation}
A particular case of the  above model is the
 so-called doubly geometric spatial autoregressive process
 $$
  X_{k,\ell}
  =\alpha X_{k-1,\ell}+\beta X_{k,\ell-1}-\alpha\beta X_{k-1,\ell-1}
   +\varepsilon_{k,\ell},
 $$
 introduced by \citet{martin1}. This was the first spatial 
autoregressive model for which unstability has been studied.
It is, in fact, the simplest spatial model, since the product structure
 \ $\varphi(x,y)=xy-\alpha x-\beta y+\alpha\beta=(x-\alpha)(y-\beta)$ \ of
 its characteristic polynomial ensures that it can be considered as some kind
 of combination of two autoregressive processes on the line, and several
 properties can be derived by the analogy of one--dimensional autoregressive
 processes.
This model has been used by \citet{jain} in the study of image processing,
 by \citet{martin2}, \citet{cg}, \citet{br1} in agricultural trials and 
by \citet{tj1} in digital filtering. 

In the stable case when \ $|\alpha|<1$ \ and
 \ $|\beta|<1$, \ asymptotic normality of several estimators
 \ $(\widehat\alpha_{m,n},\widehat\beta_{m,n})$ \ of \
 $(\alpha,\beta)$ \ based on the observations
 \ $\{X_{k,\ell}:\text{$1\leq k\leq m$ \ and \ $1\leq\ell\leq n$}\}$ \ has been
 shown (e.g. \citet{tj2,tj3} or \citet{br3,br2}), namely,
 $$
  \sqrt{mn}
  \begin{bmatrix}
   \widehat\alpha_{m,n}-\alpha\\ 
   \widehat\beta_{m,n}-\beta
  \end{bmatrix}
  \distr{\mathcal N}(0,\Sigma_{\alpha,\beta})
 $$
as  \ $m,n\to\infty$ \ with \ $m/n\to\,\textup{constant}>0$ \ with some 
covariance matrix \ $\Sigma_{\alpha,\beta}$.
 
In the unstable  case when \ $\alpha=\beta=1$, \ in contrast
 to the classical first order autoregressive time series model, where 
the appropriately normed least squares estimator (LSE) of the autoregressive 
parameter converges to a fraction of functionals of the
standard Brownian motion (see e.g. \citet{phil} or \citet{cw}),
the sequence of Gauss--Newton estimators
 \ $(\widehat\alpha_{n,n},\widehat\beta_{n,n})$ \ of \
 $(\alpha,\beta)$ \ has been shown to 
 be asymptotically normal (see \citet{bhat1} and  \citet{bhat2}). 
In the unstable case \ $\alpha=1$, \ $|\beta|<1$ \ the LSE turns out
to be asymptotically normal again \citep{bhat1}.  

\citet{bpz1} discussed  a special case of the
model \eqref{brmod}, namely, when
 \ $p_1=p_2=1$, \ $\alpha_{0,1}=\alpha_{1,0}=:\alpha$ \ and \
 $\alpha_{1,1}=0$, which is the simplest spatial model, that can not
 be reduced somehow to autoregressive models on the line. This model
 is  stable  in case \ $|\alpha|<1/2$
 \ (see e.g. \citet{whittle}, \citet{besag} or \citet{br2}), and unstable  if
 \ $|\alpha|=1/2$. \ In \citet{bpz1} the asymptotic normality of the
LSE of the unknown parameter \ $\alpha$ \ 
is proved both in stable and unstable cases. The case  $p_1=p_2=1$, \
$\alpha_{1,0}=:\alpha, \ \alpha_{0,1}=:\beta$ \ and \
 $\alpha_{1,1}=0$ \ was studied by \citet{paul} and \citet{bpz2}. This model
 is  stable  in case \ $|\alpha|+|\beta|<1$
 \  and unstable  if \ $|\alpha|+|\beta|=1$ \ \citep{br2}.
 \citet{paul} determined the exact asymptotic behaviour of the
 variances of the process, while \citet{bpz2} proved the asymptotic
 normality of the LSE of the parameters both in
 stable and unstable cases.

\begin{figure}[tb]
 \begin{center}
\leavevmode
\epsfig{file=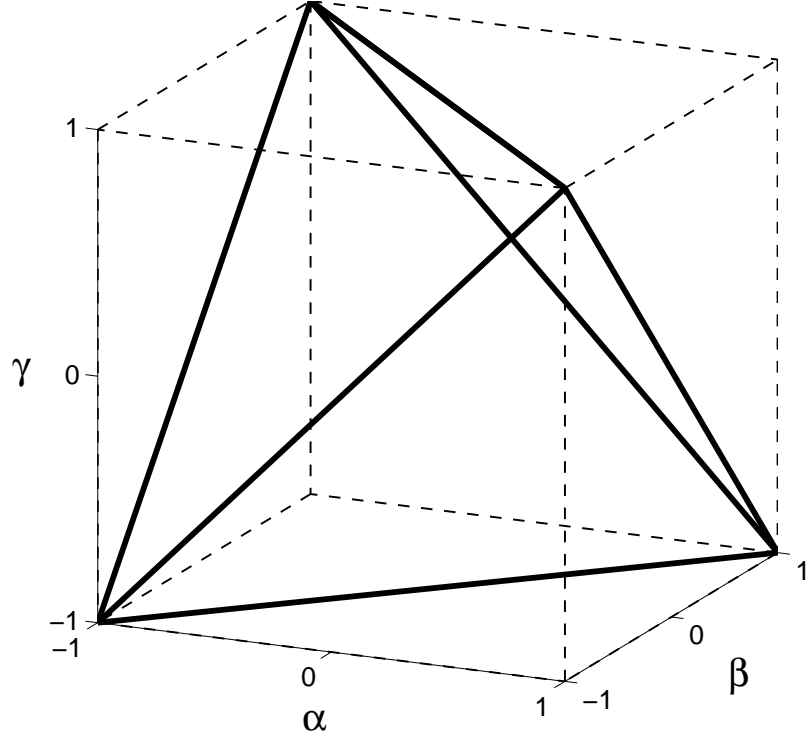,height=11cm}
\vskip1mm
\caption{The domain of stability of model \eqref{model}.}  
\label{fig:fig1}
  \end{center}  
\end{figure}
In the present paper we study the asymptotic properties of a more 
complicated  special case of the unilateral model \eqref{brmod} with 
 \ $p_1=p_2=1$, \   $\alpha_{1,0}=:\alpha$, \ $\alpha_{0,1}=:\beta$ \ and \
$\alpha_{1,1}=:\gamma$. \ In a recent paper \citet{gk} proved the
asymptotic normality of minimum distance estimators in the case when
model equation is valid on \ ${\mathbb Z}^2$. \ Here we deal with a 
model with boundary conditions, namely, we consider the spatial
autoregressive process 
\ $\{X_{k,\ell}:k,\ell\in{\mathbb Z},\,k, \ell\geq0\}$ \ defined as
 \begin{equation}\label{model}
  \begin{cases} X_{k,\ell}=
                     \alpha X_{k-1,\ell}+\beta X_{k,\ell-1}+\gamma
                     X_{k-1,\ell-1}+\varepsilon_{k,\ell}, 
                    & \text{for \ $k,\ell \geq 1$,}\\
                    X_{k,0}=X_{0,\ell}=0,  & \text{for \ $k,\ell\geq0$.}
                   \end{cases}
 \end{equation}
This process has already been examined in \citet{baran}
where the asymptotic behaviour of the variances is clarified. 
The model is stable if \ $(\alpha, \beta, \gamma) \in {\mathcal S}$, \
where \ ${\mathcal S}$ \ is the open tetrahedron with vertices
\begin{equation*}
   {\mathcal V} := \{ (1,1,-1), \, (1,-1,1), \, (-1,1,1), \, (-1,-1,-1) \}
\end{equation*}
(see Figure \ref{fig:fig1}). This result was proved by  \citet{br2}
where the tetrahedron \ ${\mathcal S}$ \ was described by conditions
\ $\ |\alpha |<1, \ |\beta |<1$ \ and \ $|\gamma|<1$, \ 
$|1+\alpha ^2-\beta ^2  -\gamma ^2|>2|\alpha +\beta\gamma|$ \ and \
$1-\beta ^2>|\alpha +\beta\gamma|$. \  Short calculation
shows that condition of stability means that \ $\
|\alpha |<1, \ |\beta |<1$ 
\ and \ $|\gamma|<1$, \ and inequalities 
\begin{equation*}
\alpha-\beta -\gamma <1, \qquad  -\alpha+\beta -\gamma <1, \qquad
-\alpha-\beta +\gamma<1, \qquad  \alpha+\beta +\gamma <1 
\end{equation*}
hold. 
Obviously, in case \ $\alpha\beta\gamma\geq 0$ \ the above set of conditions
reduces to \ $|\alpha|+|\beta|+|\gamma|<1$. 

The model is unstable if \ $(\alpha, \beta, \gamma)$ \ lies on the boundary of
  \ ${\mathcal S}$, \ when one can distinguish three cases:

\noindent
{\sc Case A.}  
The parameters are in the interior of the faces of the boundary of \
${\mathcal S}$, 
 \ i.e. $(\alpha, \beta, \gamma) \in {\mathcal F}$, \ where \ ${\mathcal F} :=
   {\mathcal F}_+ \cup {\mathcal F}_-$ 
 \ with
 \begin{align*}
  {\mathcal F}_+ &:= \{ (\alpha, \beta, \gamma) \in (-1, 1)^3
               : \alpha \beta \gamma \geq 0, \,
                 |\alpha| + |\beta| + |\gamma| = 1 \} \\ 
        &\phantom{\quad\:\:}
            \cup \{ (\alpha, \beta, \gamma) \in (-1, 1)^3
                    : \alpha \beta \gamma < 0, \,
                      |\alpha| + |\beta| - |\gamma| = 1 \} , \\ 
  {\mathcal F}_- &:= \{ (\alpha, \beta, \gamma) \in (-1, 1)^3
               : \alpha \beta \gamma < 0, \,
                 |\alpha| - |\beta| + |\gamma| = 1 \} \\
        &\phantom{\quad\:\:}
            \cup \{ (\alpha, \beta, \gamma) \in (-1, 1)^3
               : \alpha \beta \gamma < 0, \,
                 - |\alpha| + |\beta| + |\gamma| = 1 \} . 
 \end{align*}

\noindent
{\sc Case B.} \ The parameters are in the interior of the edges of the
boundary of \ ${\mathcal S}$, 
 \ i.e.\  $(\alpha, \beta, \gamma) \in {\mathcal E}$, \ where
 \ ${\mathcal E} := {\mathcal E}_1 \cup {\mathcal E}_2 \cup {\mathcal
   E}_3$ \ with 
 \begin{align*}
  {\mathcal E}_1 &:= \{ (1, \beta, \gamma) : \beta \in (-1, 1), \,
  \gamma = - \beta \} 
            \cup \{ (-1, \beta, \gamma)
                    : \beta \in (-1, 1), \, \gamma = \beta \} , \\
  {\mathcal E}_2 &:= \{ (\alpha, 1, \gamma)
               : \alpha \in (-1, 1), \, \gamma = - \alpha \}
            \cup \{ (\alpha, -1, \gamma)
                    : \alpha \in (-1, 1), \, \gamma = \alpha \} , \\
  {\mathcal E}_3 &:= \{ (\alpha, \beta, 1)
               : \alpha \in (-1, 1), \, \beta = - \alpha \}
            \cup \{ (\alpha, \beta, -1)
                    : \alpha \in (-1, 1), \, \beta = \alpha \} .
 \end{align*}
Observe that in each of the above three cases exactly two of the defining
equations of set \ ${\mathcal F}$ \ are satisfied. In this way Case B can be
considered as an extension of Case A to the situation when \
$\alpha\beta\gamma \leq 0$ \ and one of the
parameters equals \ $\pm 1$, \ while the other two parameters have
absolute values less than one. 

\noindent
Further, observe that in the first two cases \
$\gamma=-\alpha \beta$, \ so we 
obtain special cases of the doubly geometric model. If \ $(\alpha,
\beta, \gamma) \in {\mathcal E}_1$  \ then for \ $k\in{\mathbb N}$ \ the
difference \ 
$\Delta_{1,\alpha}X_{k,\ell}:=X_{k,\ell}-\alpha X_{k-1,\ell}$ \ is a
classical AR(1) process, i.e. \ $\Delta_{1,\alpha}X_{k,\ell}=\beta
\Delta_{1,\alpha}X_{k,\ell-1} +\varepsilon_{k,\ell}$. \ Similarly, if
\  $(\alpha,\beta, \gamma) \in {\mathcal E}_2$  \ then \
$\Delta_{2,\beta}X_{k,\ell}=\alpha 
\Delta_{2,\beta}X_{k-1,\ell} +\varepsilon_{k,\ell}$, \ where \
$\Delta_{2,\beta}X_{k,\ell}:=X_{k,\ell}-\beta X_{k,\ell-1}, \ \ell\in{\mathbb N}$. \

\noindent
{\sc Case C.} \ The parameters are in the vertices of the boundary of
the domain of stability,  i.e. \ $(\alpha, \beta, \gamma) \in {\mathcal V}$.

For a set \ $H\subset\{(k,\ell)\in {\mathbb Z}^2:k,\ell\geq1\}$, \ the
least squares estimator \
$(\widehat\alpha_H,\widehat\beta_H,\widehat\gamma_H)$ \ of \
$(\alpha,\beta,\gamma)$ \ based on the 
observations \ $\{X_{k,\ell}:(k,\ell)\in H\}$ \ can be obtained by 
minimizing the sum of squares
 \begin{equation*}
\sum_{(k,\ell)\in H}
    \big(X_{k,\ell}-\alpha X_{k-1,\ell}-\beta
    X_{k,\ell-1}-\gamma X_{k-1,\ell-1}\big)^2 
\end{equation*}
with respect to \ $\alpha, \beta $ \ and \ $\gamma$,  \ and it has the form
 \begin{equation*}
 \begin{bmatrix}\widehat\alpha_H\\ \widehat\beta_H \\
   \widehat\gamma_H \end{bmatrix} 
   =\left(\sum_{(k,\ell)\in H}
       \begin{bmatrix}
        X_{k-1,\ell} \\
        X_{k,\ell-1} \\
        X_{k-1,\ell-1}
       \end{bmatrix}
        \begin{bmatrix}
        X_{k-1,\ell} \\
        X_{k,\ell-1} \\
        X_{k-1,\ell-1}
       \end{bmatrix}^{\top}
     \right)^{-1}
      \sum_{(k,\ell)\in H} X_{k,\ell}
       \begin{bmatrix}
        X_{k-1,\ell} \\
        X_{k,\ell-1} \\
        X_{k-1,\ell-1}
       \end{bmatrix}.
\end{equation*}
For \ $n,m\in{\mathbb N}$  \ consider the rectangle
 \begin{equation*}
R_{n,m}:=\{(i,j)\in{\mathbb N}^2:
   \text{$1\leq i\leq n$ \ and \ $1\leq j\leq m$}\}.
\end{equation*}
For simplicity, we shall write \ $R_n:=R_{n,n}$ \ for \ $n\in{\mathbb N}$.

\begin{Thm}\label{main} \ Let \
  $\{\varepsilon_{k,\ell}:k,\ell\in{\mathbb N}\}$ \ be independent
  random 
variables with \ ${\mathsf E}\,\varepsilon_{k,\ell}=0$, \ 
$\Var\,\varepsilon_{k,\ell}=1$ \ and \ $\sup\{{\mathsf
  E}\,\varepsilon_{k,\ell}^8:k,\ell\in{\mathbb N}\}<\infty$. 
\ Assume that model \eqref{model} is satisfied.

\noindent
If \ $(\alpha, \beta, \gamma) \in {\mathcal F}_+$ \ then  
\begin{equation*}
(nm)^{1/2} \big(\widehat\alpha_{R_{n,m}}-\alpha,
                   \widehat\beta_{R_{n,m}}-\beta, 
                    \widehat\gamma_{R_{n,m}}-\gamma\big)^{\top}\distr {\mathcal
  N}\Big(0,{\mathcal H}^{\top}{\mathcal K}_{\alpha,\beta}^{-1}{\mathcal H}\Big ) 
\end{equation*}\
as \ $m,n\to\infty$ \ with $m/n\to {\mathrm constant}>0$, \ where 
\begin{equation*}
{\mathcal H}:=\begin{bmatrix}1&0&-1\\0&1&-1 \end{bmatrix}, \qquad \qquad
{\mathcal K}_{\alpha,\beta}  := 
\begin{bmatrix}
\kappa_{\beta,\alpha}^{(1)}& \kappa_{\alpha,\beta}^{(2)}\\
\kappa_{\alpha,\beta}^{(2)}& \kappa_{\alpha,\beta}^{(1)} 
\end{bmatrix},
\end{equation*}
with
\begin{equation*}
\kappa^{(1)}_{\alpha,\beta}:=(1-\beta^2)^{-1}+(1-\alpha)^2\varrho
^{(1)}_{\alpha,\beta}, \qquad \kappa^{(2)}_{\beta,\alpha}:=(1-\alpha)(1-\beta)\varrho
^{(2)}_{\alpha,\beta},
\end{equation*}
and 
\begin{align*}
\varrho^{(1)}_{\alpha,\beta}:=\sum _{k=0}^{\infty}\sum
_{\ell=0}^{\infty} &\Big ({\mathsf
  P}\big(S_{k,\ell}^{(|\alpha|,1-|\beta|)}=k+1\big)-{\mathsf
  P}\big(S_{k,\ell}^{(|\alpha|,1-|\beta|)} =k\big) \Big)^2,\\
\varrho^{(2)}_{\alpha,\beta}:=\sum _{k=0}^{\infty}\sum
_{\ell=0}^{\infty} &\Big ({\mathsf
  P}\big(S_{k,\ell}^{(|\alpha|,1-|\beta|)} =k+1\big)-{\mathsf
  P}\big(S_{k,\ell}^{(|\alpha|,1-|\beta|)}=k\big) \Big) \\ 
  &\times\Big ({\mathsf
  P}\big(S_{\ell,k}^{(|\beta|,1-|\alpha|)}=\ell+1\big)-{\mathsf
  P}\big((S_{\ell,k}^{(|\beta|,1-|\alpha|)}=\ell\big) \Big),
\end{align*}
where \ $S_{k,\ell}^{(\nu,\mu)}:=\xi_k^{(\nu)}+\eta_{\ell }^{(\mu)}$ \
and \ $\xi_k^{(\nu)}$ \ and  \ 
$\eta_{\ell }^{(\mu)}$ \ are independent binomial random variables
with parameters \ $(k,\nu)$ \ and \ $(\ell,\mu)$, \
respectively.

\noindent
If \ $(\alpha, \beta, \gamma) \in {\mathcal E}_1\cup {\mathcal E}_2$ \ 
then
\begin{equation*}
(nm)^{1/2} \big(\widehat\alpha_{R_{n,m}}-\alpha,
                   \widehat\beta_{R_{n,m}}-\beta, 
                    \widehat\gamma_{R_{n,m}}-\gamma\big)^{\top}\distr {\mathcal
  N}\Big(0, {\overline\Sigma} _{\alpha,\beta}\Big )
\end{equation*}\
as \ $m,n\to\infty$ \ with $m/n\to {\mathrm constant}>0$, \ where 
\begin{equation*}
\Sigma_{\alpha,\beta}
:=\begin{cases}
      \begin{bmatrix}
        1 & |\beta|\Sign(\alpha) & |\beta|\\
        |\beta|\Sign(\alpha) & 1 & \Sign(\alpha) \\
        |\beta| & \Sign(\alpha) & 1
      \end{bmatrix}, \qquad \text{if \ $|\alpha|=1$,} \\[10mm]
      \begin{bmatrix}
        1 & |\alpha|\Sign(\beta) & \Sign(\beta) \\
        |\alpha\Sign(\beta)| & 1 & |\alpha| \\
        \Sign(\beta) & |\alpha|  & 1
      \end{bmatrix}, \qquad \text{if \ $|\beta|=1$,} 
   \end{cases}
\end{equation*}
and \ $\overline A$ \ denotes the adjoint of a matrix \ $A$.

If \ $(\alpha, \beta, \gamma) \in {\mathcal V}$ \ then
\begin{equation*}
(nm)^{3/4} \big(\widehat\alpha_{R_{n,m}}-\alpha,
                   \widehat\beta_{R_{n,m}}-\beta, 
                    \widehat\gamma_{R_{n,m}}-\gamma\big)^{\top}\distr {\mathcal
  N}\Big(0,\Theta_{\alpha,\beta}\Big )
\end{equation*}\
as \ $m,n\to\infty$ \ with $m/n\to {\mathrm constant}>0$, \ where 
\begin{equation*}
\Theta_{\alpha,\beta}:=2\begin{bmatrix}
        1 & 0 & -\beta\\
        0& 1 & -\alpha \\
        -\beta & -\alpha & 2
      \end{bmatrix}.
\end{equation*}
\end{Thm}

We remark that results given Theorem \ref{main} do not cover the cases \
$(\alpha, \beta, \gamma) \in {\mathcal F}_-\cup{\mathcal E}_3$. \
The main problem is that in these cases we could not
handle the asymptotic behaviour of the covariance structure. A
more detailed explanation and some results on the missing cases 
can be found in \citet{baran}.  Another problem is that we were not able
to find closed forms of \ 
$\varrho^{(i)}_{\alpha,\beta}, \ i=1,2,$ \ and in this way we don't know
how they depend on the parameters.

For the sake of simplicity, we carry out the proof of Theorem
\ref{main} only for \ $m=n$.
\ The general case can be handled with slight modifications.
We can write
\begin{equation*}
\begin{bmatrix}\widehat\alpha_{R_n}-\alpha\\
                   \widehat\beta_{R_n}-\beta \\
                    \widehat\gamma_{R_n}-\gamma\end{bmatrix} 
    =B_n^{-1}A_n,
\end{equation*}
with
\begin{equation*}
 A_n:=\sum_{(k,\ell)\in R_n}\varepsilon_{k,\ell}
       \begin{bmatrix}
        X_{k-1,\ell} \\
        X_{k,\ell-1} \\
        X_{k-1,\ell-1}
       \end{bmatrix}, \ \ \
  B_n:=\sum_{(k,\ell)\in R_n}
                \begin{bmatrix}
        X_{k-1,\ell} \\
        X_{k,\ell-1} \\
        X_{k-1,\ell-1}
       \end{bmatrix}
        \begin{bmatrix}
        X_{k-1,\ell} \\
        X_{k,\ell-1} \\
        X_{k-1,\ell-1}
       \end{bmatrix}^{\top}.
\end{equation*}
Concerning the asymptotic behaviour of \ $A_n$ \ and \ $B_n$ \ we
can prove the following propositions.
\begin{Pro}\label{Bn} \
If   \ $(\alpha, \beta, \gamma) \in {\mathcal F}_+$ \ then
\begin{equation*}
n^{-5/2}B_n\qmean \sigma_{\alpha,\beta}^2\Psi _{\alpha,\beta} \qquad
\text{as \ $n\to\infty$,}
\end{equation*}
where
\begin{align*}
\sigma_{\alpha,\beta}^2&:=\frac23\bigg( \frac{1-|\alpha|\lor|
  \beta|}{\pi(|\alpha|+|\beta|)}\bigg)^{1/2}\Bigg(\frac
1{(1-|\alpha|)(1-|\beta|)}-\frac 1{5(1-|\alpha|\land|\beta|)^2}\Bigg), \\
\Psi _{\alpha,\beta}&:=\begin{bmatrix}
        1 & \Sign (\alpha\beta) & \Sign (\beta)\\
        \Sign(\alpha\beta) & 1 & \Sign (\alpha) \\
        \Sign(\beta) & \Sign (\alpha) & 1
      \end{bmatrix}.
\end{align*}
If \ $(\alpha, \beta, \gamma) \in {\mathcal E}_1\cup{\mathcal E}_2$ \ then
\begin{equation*}
n^{-3}B_n\qmean \big(2(1-\gamma^2)\big)^{-1}\Sigma_{\alpha, \beta}
\qquad \text{as \ $n\to\infty$.} 
\end{equation*}
If \ $(\alpha, \beta, \gamma) \in {\mathcal V}$ \ then
\begin{align*}
n^{-4}B_n&\distr \int\limits_0^1\!\! \int\limits_0^1{\mathcal
  W}^2(s,t)\,{\mathrm d}s\,{\mathrm d} t 
 \begin{bmatrix}\alpha\\
                 \beta \\
                  \alpha\beta\end{bmatrix}\begin{bmatrix}\alpha\\
                 \beta \\
                  \alpha\beta\end{bmatrix}^{\top}
               \qquad \text{as \ $n\to\infty$,} 
\end{align*}
where \ ${\mathcal W}(s,t)$ \ is a standard Wiener sheet.

\end{Pro}

\begin{Pro}\label{An} \
If  \ $(\alpha, \beta, \gamma) \in {\mathcal F}_+$ \ then
\begin{equation*}
n^{-5/4}A_n\distr{\mathcal N}\big (0,\sigma_{\alpha,\beta}^2
\Psi _{\alpha,\beta}\big ) \qquad \text{as \ $n\to\infty$.}
\end{equation*} 
If \ $(\alpha, \beta, \gamma) \in {\mathcal E}_1\cup{\mathcal E}_2$ \ then
\begin{equation*}
n^{-3/2}A_n\distr{\mathcal N}\Big
(0,\big(2(1-\gamma^2)\big)^{-1}\Sigma_{\alpha, \beta}
\Big )  \qquad \text{as \ $n\to\infty$.}
\end{equation*}
If \ $(\alpha, \beta, \gamma) \in {\mathcal V}$ \ then
\begin{equation*}
n^{-2}A_n\distr \int\limits_0^1\!\!\int\limits_0^1{\mathcal
  W}(s,t){\mathcal W}({\mathrm d}s,{\mathrm d} t) 
 \begin{bmatrix}\alpha\\
                 \beta \\
                  \alpha\beta\end{bmatrix} \qquad \text{as \ $n\to\infty$.}
\end{equation*}
\end{Pro}

As the limits in Proposition \ref{Bn} are singular, the statements of
Theorem \ref{main} can not be obtained directly from Propositions
\ref{An} and \ref{Bn}. Hence, one has to use the same idea as in
\citet{bpz2} and consider \ $B_n^{-1}={\overline B}_n/{\det (B_n)}$.

\begin{Pro}\label{DetB}
\ If  \ $(\alpha, \beta, \gamma) \in {\mathcal F}_+$ \  then 
\begin{equation*}
n^{-13/2} \det (B_n)\stoch \sigma^2_{\alpha,\beta} \det \big({\mathcal
  K}_{\alpha,\beta}\big)>0 \qquad \text{as \ $n\to\infty$.} 
\end{equation*}
If \ $(\alpha, \beta, \gamma) \in {\mathcal E}_1\cup{\mathcal E}_2$ \ then
\begin{equation*}
n^{-8} \det (B_n)\stoch \big(2(1-\gamma^2)\big)^{-2}  \qquad \text{as \
  $n\to\infty$.} 
\end{equation*}
If \ $(\alpha, \beta, \gamma) \in {\mathcal V}$ \ then
\begin{equation*}
n^{-10} \det (B_n)\distr \frac 14  \int\limits_0^1\!\!\int\limits_0^1{\mathcal
  W}^2(s,t)\,{\mathrm d}s\,{\mathrm d} t  \qquad \text{as \
  $n\to\infty$.} 
\end{equation*}
\end{Pro}

We remark that using higher moment conditions on the innovations \
$\varepsilon _{k,\ell}$, \ after tedious but straightforward
calculations, instead of stochastic convergence one can also prove
${\mathsf L}_2$ convergence in the first two statements of Proposition
\ref{DetB}.   

Further, if we take appropriate linear transformations of \ $A_n$ \ we
have asymptotic normality in all of the unstable cases considered. Let \
$C_n:={\mathcal H}A_n=\big(C_n^{(1)},\  C_n^{(2)}\big)^{\top}$, where  
\begin{align}
   \label{Cn1}
C_n^{(1)}&:=\big (1, \ 0, \ -1 \big)A_n=\sum_{(k,\ell )\in
  {R_n}} \big(X_{k-1,\ell }-X_{k-1,\ell-1}\big)\varepsilon_{k,\ell},\\
C_n^{(2)}&:=\big (0, \ 1, \ -1 \big)A_n=\sum_{(k,\ell )\in
  {R_n}} \big(X_{k,\ell -1}-X_{k-1,\ell-1}\big)\varepsilon_{k,\ell}. \label{Cn2}
\end{align}

\begin{Pro}\label{Cnlim} \
If  \ $(\alpha, \beta, \gamma) \in {\mathcal F}_+$ \  then
\begin{equation*}
n^{-1}C_n\distr{\mathcal N}\big (0, {\mathcal K}_{\alpha,\beta} \big )
\qquad \text{as \ $n\to\infty$.} 
\end{equation*} 
If \ $(\alpha, \beta, \gamma) \in {\mathcal E}_1$ \ then
\begin{equation*}
n^{-3/2}C_n^{(1)}\distr{\mathcal N}\big (0, (1-\gamma)^{-1}\big )
\quad \text{and} \quad  
n^{-1}C_n^{(2)}\distr{\mathcal N}\big (0, (1-\gamma^2)^{-1}\big )
\qquad \text{as \ $n\to\infty$.} 
\end{equation*}
If \ $(\alpha, \beta, \gamma) \in {\mathcal E}_2$ \ then
\begin{equation*}
n^{-1}C_n^{(1)}\distr{\mathcal N}\big (0, (1-\gamma^2)^{-1}\big )
\quad \text{and} \quad  
n^{-3/2}C_n^{(2)}\distr{\mathcal N}\big (0, (1-\gamma)^{-1}\big )
\qquad \text{as \ $n\to\infty$.} 
\end{equation*}
If \ $(\alpha, \beta, \gamma) \in {\mathcal V}$ \ then
\begin{equation*}
n^{-3/2}C_n\distr{\mathcal N}\big (0, {\mathcal I}_2/2\big ) \qquad
\text{as \ $n\to\infty$,} 
\end{equation*}
where \ ${\mathcal I}_2$ \ denotes the two-by-two unit matrix.
\end{Pro}

\begin{Pro}\label{LimAn} \
If  \ $(\alpha, \beta, \gamma) \in {\mathcal F}_+$ \  then
\begin{equation*}
n^{-11/2} {\overline B}_n A_n\distr {\mathcal
  N}\Big(0,\sigma_{\alpha,\beta}^4\det \big({\mathcal K}_{\alpha,\beta}\big)
{\mathcal H}^{\top}{\overline{\mathcal K}}_{\alpha,\beta} {\mathcal H}\Big )
\qquad \text{as \ 
  $n\to\infty$.}
\end{equation*} 
If \ $(\alpha, \beta, \gamma) \in {\mathcal E}_1\cup{\mathcal E}_2$
\ then
\begin{equation*}
n^{-7} {\overline B}_n A_n\distr {\mathcal
  N}\Big(0,\big(2(1-\gamma^2\big))^{-4}{\overline\Sigma} _{\alpha,\beta}\Big )
\qquad \text{as \ 
  $n\to\infty$.}
\end{equation*} 
\end{Pro}

The aim of the following discussion is to show that it suffices to
prove Theorem \ref{main} for \ $\alpha \geq 0$, \ $\beta \geq 0$ \ and
\ $\gamma \geq 0$ \ if \ $\alpha\beta\gamma\geq 0$ \ and for  \
$\alpha >0$, \ $\beta >0$ \ and \ $\gamma < 0$ \ if \
$\alpha\beta\gamma <0$. \ In this way instead of \ ${\mathcal F}_+, \
{\mathcal V}, \  
{\mathcal E}_1$ \ and \ ${\mathcal E}_2$ \ it suffices
to use their subsets
\begin{alignat*}{2}
{\mathcal F}_{++}&=\!\{(\alpha,\beta,\gamma) : 0\leq\alpha,\beta
<1,|\gamma |<1, \alpha \!+\!\beta \!+\!\gamma =1\},\quad &&{\mathcal
  V}_{+}\!=\!\{(1,1,-1)\}, \\
{\mathcal E}_{1+}&=\!\{(1,\beta,\gamma) : 0\leq\beta=-\gamma <1\}, 
&&{\mathcal E}_{2+}\!=\!\{(\alpha,1,\gamma) : 0\leq\alpha=-\gamma <1\}, 
\end{alignat*}
respectively.

First we note that direct calculations imply 
 \begin{align}\label{marep}
  X_{k,\ell}&=\sum_{i=1}^k\sum_{j=1}^{\ell}G(k-i,\ell-j;\alpha,\beta,\gamma)
            \varepsilon_{i,j}\\ 
           &=\sum_{i=1}^k\sum_{j=1}^{\ell}\binom{k+\ell -i-j}{\ell
             -j}\alpha^{k-i}\beta^{\ell-j}F\bigg(i-k,j-\ell\,;i+j-k-\ell\,;
           -\frac{\gamma}{\alpha\beta}\bigg) \varepsilon_{i,j}, \label{marepf}
 \end{align}
$k,\ell \geq 1$, \ where \eqref{marepf} holds only for \
$\alpha\beta\ne 0$,
\begin{equation*}
G(m,n;\alpha,\beta,\gamma):=\sum_{r=0}^{m\land n} \frac{(m+n-r)!}{(m-r)!(n-r)!r!}
                     \alpha^{m-r}\beta^{n-r}\gamma ^r, \qquad m,n\in
                     {\mathbb N}\cup \{0\},
\end{equation*}
and \ $F(-n,b;c;z)$ \ is the Gauss hypergeometric
function defined by
\begin{equation*}
F(-n,b;c;z):=\sum_{r=0}^{n }
\frac {(-n)_r(b)_r}{(c)_rr!}z^r, \qquad n\in{\mathbb N},  \quad
b,c,z\in{\mathbb C},
\end{equation*}
and \ $(a)_r:=a(a+1)\dots (a+r-1)$ \ (for the definition in more
general cases see e.g. \citet{be}).

Observe that as for \ $m,n\in{\mathbb N}$ \ we have \
$F(-n,-m;-n-m;1)=\binom{m+n}n^{-1}$ \ and \ $F(-n,-m;-n-m;0)=1$, \ 
moving average representations of the doubly geometric model of
\citet{martin1} and of the spatial models studied by \citet{paul} and
\citet{bpz1,bpz2}, respectively, are special forms of \eqref{marepf}.

Now, put \
$\widetilde\varepsilon_{k,\ell}:=(-1)^{k+\ell}\varepsilon_{k,\ell}$ \
for \ $k,\ell \in{\mathbb N}$. \ Then \
$\{\widetilde\varepsilon_{k,\ell}:k,\ell\in{\mathbb N}\}$ \ are
independent random variables with \ ${\mathsf
  E}\,\widetilde\varepsilon_{k,\ell}=0$, \ 
$\Var\,\widetilde\varepsilon_{k,\ell}=1$  \ and \ $\sup\{{\mathsf
  E}\,\widetilde\varepsilon_{k,\ell}^8:k,\ell\in{\mathbb
  N}\}<\infty$. 
 \ Consider the process \ $\{\widetilde
X_{k,\ell}:k,\ell\in{\mathbb Z},\,k, \ell\geq0\}$ \  defined as 
 \begin{equation*}
  \begin{cases} \widetilde X_{k,\ell}=
                     -\alpha \widetilde X_{k-1,\ell}-\beta \widetilde
                     X_{k,\ell-1}+\gamma \widetilde
                     X_{k-1,\ell-1}+\widetilde\varepsilon_{k,\ell},  
                    & \text{for \ $k,\ell \geq 1$,}\\
                    \widetilde X_{k,0}=\widetilde X_{0,\ell}=0,  &
                    \text{for \ $k,\ell\geq0$.} 
                   \end{cases}
 \end{equation*}
Then by representation \eqref{marep} for \ $k,\ell\in{\mathbb N}$ \
we have
\begin{align*}
\widetilde X_{k,\ell}&=\sum_{i=1}^k\sum_{j=1}^{\ell}
  G(k-i,\ell-j;-\alpha,-\beta,\gamma)
  \widetilde\varepsilon_{i,j} \\
&=\sum_{i=1}^k\sum_{j=1}^{\ell}
  (-1)^{k+\ell-i-j}G(k-i,\ell-j;\alpha,\beta,\gamma)
  \widetilde\varepsilon_{i,j}=(-1)^{k+\ell}X_{k,\ell}.
\end{align*}
Hence
\begin{align*}
 \widetilde A_n&:=\sum_{(k,\ell)\in R_n} \widetilde\varepsilon_{k,\ell}
       \begin{bmatrix}
        \widetilde X_{k-1,\ell} \\
        \widetilde X_{k,\ell-1}   \\
        \widetilde X_{k-1,\ell-1}
       \end{bmatrix}=
       \begin{bmatrix}
        -1 & 0 & 0\\0 &-1 &0\\0 &0 &1
       \end{bmatrix}A_n, \\
  \widetilde B_n&:= \sum_{(k,\ell)\in R_n}
               \begin{bmatrix}
                 \widetilde X_{k-1,\ell} \\
                 \widetilde X_{k,\ell-1}   \\
                 \widetilde X_{k-1,\ell-1}
               \end{bmatrix}
               \begin{bmatrix}
                 \widetilde X_{k-1,\ell} \\
                 \widetilde X_{k,\ell-1}   \\
                 \widetilde X_{k-1,\ell-1}
               \end{bmatrix}^{\top}
               =\begin{bmatrix}
        -1 & 0 & 0\\0 &-1 &0\\0 &0 &1
       \end{bmatrix}B_n
      \begin{bmatrix}
        -1 & 0 & 0\\0 &-1 &0\\0 &0 &1
       \end{bmatrix}.
\end{align*}
Consequently, in order to prove Propositions \ref{Bn} -- \ref{LimAn} for
\ $\alpha\leq 0$ \ and \ $\beta \leq 0$ \ it 
suffices to prove them for \ $\alpha\geq0$ \ and \ $\beta\geq0.$

Next, put \
$\widehat\varepsilon_{k,\ell}:=(-1)^k\varepsilon_{k,\ell}$ \
for \ $k,\ell \in{\mathbb N}$. \ Then \
$\{\widehat\varepsilon_{k,\ell}:k,\ell\in{\mathbb N}\}$ \ are again
independent random variables with \ ${\mathsf
  E}\,\widehat\varepsilon_{k,\ell}=0$, \ 
$\Var\,\widehat\varepsilon_{k,\ell}=1$, \ and \ $\sup\{{\mathsf
  E}\,\widehat\varepsilon_{k,\ell}^8:k,\ell\in{\mathbb
  N}\}<\infty$. 
 \ Consider the process \ $\{\widehat
X_{k,\ell}:k,\ell\in{\mathbb Z},\,k, \ell\geq0\}$ \  defined as 
 \begin{equation*}
  \begin{cases} \widehat X_{k,\ell}=
                     -\alpha \widehat X_{k-1,\ell}+\beta \widehat
                     X_{k,\ell-1}-\gamma \widehat
                     X_{k-1,\ell-1}+\widehat\varepsilon_{k,\ell},  
                    & \text{for \ $k,\ell \geq 1$,}\\
                    \widehat X_{k,0}=\widehat X_{0,\ell}=0,  &
                    \text{for \ $k,\ell\geq0$.} 
                   \end{cases}
 \end{equation*}
Then by representation \eqref{marep} for \ $k,\ell\in{\mathbb N}$ \
we have
\begin{align*}
\widehat X_{k,\ell}&=\sum_{i=1}^k\sum_{j=1}^{\ell}
 G(k-i,\ell-j;-\alpha,\beta,-\gamma)
\widehat\varepsilon_{i,j} \\
&=\sum_{i=1}^k\sum_{j=1}^{\ell} (-1)^{k-i}
 G(k-i,\ell-j;\alpha,\beta,\gamma)
\widehat\varepsilon_{i,j}=(-1)^kX_{k,\ell}, 
\end{align*}
Hence
\begin{align*}
 \widehat A_n&:=\sum_{(k,\ell)\in R_n}\widehat \varepsilon_{k,\ell}
       \begin{bmatrix}
        \widehat X_{k-1,\ell}  \\
        \widehat X_{k,\ell-1}  \\
        \widehat X_{k-1,\ell-1}
       \end{bmatrix}=
       \begin{bmatrix}
        -1 & 0 & 0\\0 &1 &0\\0 &0 &-1
       \end{bmatrix}A_n, \\
  \widehat B_n&:= \sum_{(k,\ell)\in R_n}
       \begin{bmatrix}
        \widehat X_{k-1,\ell}  \\
        \widehat X_{k,\ell-1}  \\
        \widehat X_{k-1,\ell-1}
       \end{bmatrix}
       \begin{bmatrix}
        \widehat X_{k-1,\ell}  \\
        \widehat X_{k,\ell-1}  \\
        \widehat X_{k-1,\ell-1}
       \end{bmatrix}^{\top}
               =\begin{bmatrix}
        -1 & 0 & 0\\0 &1 &0\\0 &0 &-1
       \end{bmatrix}B_n
      \begin{bmatrix}
        -1 & 0 & 0\\0 &1 &0\\0 &0 &-1
       \end{bmatrix}.
\end{align*}
Thus, case
\ $\alpha \leq 0$ \ and \ $\beta \geq 0$ \ can also be obtained from case \
$\alpha\geq0$ \ and \ $\beta\geq0$. 

In a similar way we have \ $\bar X_{k,\ell}=(-1)^{\ell}X_{k,\ell}$, \
where \ $\{\bar X_{k,\ell}:k,\ell\in{\mathbb Z},\,k, \ell\geq0\}$ \
is defined as 
 \begin{equation*}
  \begin{cases} \bar X_{k,\ell}=
                     \alpha \bar X_{k-1,\ell}-\beta \bar
                     X_{k,\ell-1}-\gamma \bar
                     X_{k-1,\ell-1}+\bar\varepsilon_{k,\ell},  
                    & \text{for \ $k,\ell \geq 1$,}\\
                    \bar X_{k,0}=\bar X_{0,\ell}=1.  &
                    \text{for \ $k,\ell\geq0$,} 
                   \end{cases}
 \end{equation*}
with \ $\bar\varepsilon_{k,\ell}:=(-1)^{\ell}\varepsilon_{k,\ell}$. \ Hence
\begin{align*}
 \bar A_n&:=\sum_{(k,\ell)\in R_n} \bar \varepsilon_{k,\ell}
       \begin{bmatrix}
        \bar X_{k-1,\ell} \\
        \bar X_{k,\ell-1}  \\
        \bar X_{k-1,\ell-1}
       \end{bmatrix}=
       \begin{bmatrix}
        1 & 0 & 0\\0 &-1 &0\\0 &0 &-1
       \end{bmatrix}A_n, \\
  \bar B_n&:= \sum_{(k,\ell)\in R_n}
       \begin{bmatrix}
        \bar X_{k-1,\ell} \\
        \bar X_{k,\ell-1}  \\
        \bar X_{k-1,\ell-1}
       \end{bmatrix}
        \begin{bmatrix}
        \bar X_{k-1,\ell} \\
        \bar X_{k,\ell-1}  \\
        \bar X_{k-1,\ell-1}
       \end{bmatrix}^{\top}
       =\begin{bmatrix}
        1 & 0 & 0\\0 &-1 &0\\0 &0 &-1
        \end{bmatrix}B_n
        \begin{bmatrix}
        1 & 0 & 0\\0 &-1 &0\\0 &0 &-1
         \end{bmatrix}.
\end{align*}
Thus, case
\ $\alpha \geq 0$ \ and \ $\beta \leq 0$ \ can be obtained from case \
$\alpha\geq0$ \ and \ $\beta\geq0$, \ too.

\section{Covariance structure}
   \label{sec:sec2}

By representations \eqref{marep} and \eqref{marepf} we obtain that for 
\ $k_1,\ell_1,k_2,\ell_2\in {\mathbb N}$ \  and  \
$\alpha,\beta,\gamma \in{\mathbb R}$ \ we have
 \begin{align}\label{eq:eq2.1}
  \Cov \big(X_{k_1,\ell_1}&,X_{k_2,\ell_2})\!=\!\!\!\sum_{i=1}^{k_1\land k_2}
    \sum_{j=1}^{\ell_1\land
      \ell_2}\!G(k_1-i,\ell_1-j;\alpha,\beta,\gamma) 
     G(k_2-i,\ell_2-j;\alpha,\beta,\gamma)\\
           &=\!\!\!\sum_{i=1}^{k_1\land k_2}\sum_{j=1}^{\ell_1\land
             \ell_2} \binom{k_1+\ell_1 -i-j}{\ell_1 -j}
            \binom{k_2+\ell_2 -i-j}{\ell_2
              -j}\alpha^{k_1+k_2-2i}\beta^{\ell_1+\ell_2-2j} \label{eq:eq2.2}\\
          &\phantom{=\!\!\!\sum_{i=1}^{k_1\land
              k_2}}
           \times \!F\bigg(i\!-\!k_1,j\!-\!\ell_1;i\!+\!j\!-\!k_1\!-\!\ell_1;
           -\frac{\gamma}{\alpha\beta}\bigg)
           F\bigg(i\!-\!k_2,j\!-\!\ell_2;i\!+\!j\!-\!k_2\!-\!\ell_2; 
           -\frac{\gamma}{\alpha\beta}\bigg), \nonumber
 \end{align}  
where \ $x\land y:=\min\{x,y\}, \ x,y \in {\mathbb R}$, \ an empty
sum is defined to be equal to \ $0$, \ and
\eqref{eq:eq2.2} holds only for \ $\alpha\beta\ne 0$.

The following lemma \citep[Corollary 2.2]{baran} helps us to find a
more convenient form of the covariances.
\begin{Lem}
   \label{binorep} \
If \ $0\leq\alpha,\beta<1$ \ and  
\ $\alpha+\beta+\gamma=1$, \ then
\begin{equation*}
G(m,n;\alpha,\beta,\gamma)
={\mathsf P}\Big (S_{m,n}^{(\alpha,1-\beta)}=m \Big)=
{\mathsf P}\Big (S_{n,m}^{(\beta,1-\alpha)}=n \Big).
\end{equation*}
\end{Lem}

With the help of \eqref{eq:eq2.1} and Lemma \ref{binorep} one can find
upper bounds for the covariances \citep[Theorem 2.4]{baran}.
\begin{Lem} 
\label{covbound} \
If \ $(\alpha, \beta, \gamma) \in {\mathcal F}_+$ \ then
\begin{equation*}
\big |\Cov (X_{k_1,\ell_1}, X_{k_2,\ell_2})\big |\leq \
C_{\alpha,\beta} \sqrt{k_1+\ell_1+k_2+\ell_2}
 \end{equation*}
with some constant \ $C_{\alpha,\beta}>0$.

\noindent
If \ $(\alpha, \beta, \gamma) \in {\mathcal E}_1$ \
or \ $(\alpha, \beta, \gamma) \in {\mathcal E}_2$
\ then 
\begin{equation*}
\big |\Cov (X_{k_1,\ell_1}, X_{k_2,\ell_2})\big |\leq \ (k_1\land
k_2) \frac {|\gamma|^{|\ell_1 -\ell_2|}}{1-\gamma^2} \quad \text{or}
\quad
\big |\Cov (X_{k_1,\ell_1}, X_{k_2,\ell_2})\big |\leq \ (\ell_1\land
\ell_2) \frac {|\gamma|^{|k_1 -k_2|}}{1-\gamma^2},
\end{equation*}
respectively.

\noindent
If \ $(\alpha, \beta, \gamma) \in {\mathcal V}$ \ then 
\begin{equation*}
\Cov (X_{k_1,\ell_1}, X_{k_2,\ell_2})=(k_1\land k_2)
(\ell_1\land \ell_2)\,\alpha ^{|k_1-k_2|} \beta ^{|\ell_1-\ell_2|}.
\end{equation*}
\end{Lem}

For \ $n\in{\mathbb N}$ \ let us introduce the piecewise constant random
fields
\begin{equation*}
Y^{(n)}_{1,0}(s,t):=X_{[ns]+1,[nt]}, \quad
Y^{(n)}_{0,1}(s,t):=X_{[ns],[nt]+1}, \quad 
Y^{(n)}_{0,0}(s,t):=X_{[ns],[nt]}, \\
\end{equation*}
for \ $0\leq s,t\in {\mathbb R}$.

The following result is a natural, but non-trivial generalization
of Proposition 2.5 of \citet{bpz2}. 

\begin{Pro}
   \label{covdiff} \
If \ $(\alpha, \beta, \gamma) \in {\mathcal F}_{++}$ \ then
there exists a constant \ $K_{\alpha,\beta} >0$ \ such that
\begin{equation*}
\Big | \Cov
\big(Y^{(n)}_{i,j}(s_1,t_1),Y^{(n)}_{0,0}(s_2,t_2)\big)-\Cov
\big(Y^{(n)}_{0,0}(s_1,t_1),Y^{(n)}_{0,0}(s_2,t_2) \big) \Big | \leq
K_{\alpha,\beta} 
\end{equation*}
for all \ $n\in {\mathbb N}$, \ $0<s_1,t_1,s_2,t_2\in {\mathbb R}$, \
with \ $(i,j)\in \big\{(0,1),(1,0)\big\}$.
\end{Pro}

In the proof of Proposition \ref{covdiff} we make use of the following
lemmas. Lemma \ref{lclt} is an obvious generalization of Theorem 2.4 of
\citet{bpz2}, while Lemma \ref{lcltdiff} can be easily obtained from a
generalization of Theorem 2.6 of \citet{bpz2} using Taylor series expansion.

\begin{Lem}
   \label{lclt}
\ Let \ $k,\ell\in {\mathbb N}$, \ let \ $0<\mu, \ \nu<1$ \ be  real numbers
and let \ $\xi_k^{(\nu)}$ \ and  \
$\eta_{\ell }^{(\mu)}$ \ be independent binomial random variables
with parameters \ $(k,\nu)$ \ and \ $(\ell,\mu)$, \
respectively. Further, let \ $S_{k,\ell}^{(\nu,\mu)}:=\xi_k^{(\nu)}+\eta_{\ell
}^{(\mu)}$ \ and let
\begin{equation*}
m_{k,\ell}:={\mathsf E} S_{k,\ell}^{(\nu,\mu)} ,\qquad
b_{k,\ell}:=\Var \big (S_{k,\ell}^{(\nu,\mu)} \big),\qquad  
    x_{j,k,\ell}:=(j-m_{k,\ell})/\sqrt{b_{k,\ell}}.
\end{equation*}
Then for all \ $k,\ell \in{\mathbb N}$ \ and \
$j\in\{0,1,\ldots,k+\ell\}$, \ we have
\begin{equation*}
\left|{\mathsf P}\big (S_{k,\ell}^{(\nu,\mu)}=j\big )
 -\frac1{\sqrt {2\pi b_{k,\ell}}}\exp\left (-x_{j,k,\ell}^2/2\right )\right|
    \leq\frac{C_{\mu,\nu}}{b_{k,\ell}},
\end{equation*}
where \ $C_{\mu,\nu}>0$ \ is a constant depending only on \
$\mu$ \ and \ $\nu$ \ (and not 
depending on  \ $k,\ell,j$).
\end{Lem}

\begin{Lem}
  \label{lcltdiff}
Using notations of Lemma \ref{lclt}, for \ $j\in\{0,1,\ldots,k+\ell-1\}$ \ let
\begin{equation*}
\Delta_{j,k,\ell}:=\Big ({\mathsf P}\big
(S_{k,\ell}^{(\nu,\mu)}=j+1\big)-{\mathsf
P}\big (S_{k,\ell}^{(\nu,\mu)}=j\big)\Big)+\frac{x_{j,k,\ell}}{\sqrt {2\pi}
  b_{k,\ell}}\exp\big(-x_{j,k,\ell}^2/2\big ). 
\end{equation*}
Then there exists a constant \ $C_{\mu,\nu}>0$ \ depending only on \
$\mu$ \ and \ $\nu$ \ (and not 
depending on  \ $k,\ell,j$) such that
\begin{equation*}
\big |\Delta_{j,k,\ell} \big | \leq \frac {C_{\mu,\nu}}{b_{k,\ell}^{3/2}}.
\end{equation*}
\end{Lem}

\begin{Cor} 
   \label{probdiff}
\ Let \ $0<\mu, \ \nu<1$ \ be  real numbers. There exists a constant \
$C_{\mu,\nu}>0$ \ such that for all \ $k,\ell\in {\mathbb N}$ \ and \
$j\in\{0,1,\ldots,k+\ell-1\}$ \ we have
\begin{equation*}
\Big |{\mathsf P}\big(S_{k,\ell}^{(\nu,\mu)}=j+1\big)-{\mathsf
P}(S_{k,\ell}^{(\nu,\mu)}=j\big)\Big| \leq \frac {C_{\mu,\nu}}{b_{k,\ell}}.
\end{equation*}
\end{Cor}

\noindent
{\bf Proof of Proposition \ref{covdiff}} \ Without loss of generality
we may assume \ $(i,j)=(1,0)$. \ Let
\begin{equation*}
\omega^{(n)}_{\alpha,\beta}(s_1,t_1,s_2,t_2):=\Cov
\big(Y^{(n)}_{1,0}(s_1,t_1),Y^{(n)}_{0,0}(s_2,t_2)\big)-\Cov
\big(Y^{(n)}_{0,0}(s_1,t_1),Y^{(n)}_{0,0}(s_2,t_2) \big).
\end{equation*}
Consider first the case \ $[ns_1]\geq [ns_2]$ \ and \ $[nt_1]\geq [nt_2]$. \
Obviously, one may assume \ $[ns_2]\geq 1$ \ and \ $[nt_2]\geq 1$. \ From
the definition of random fields \ $Y^{(n)}_{1,0}$ \ and \
$Y^{(n)}_{0,0}$ \ with the help of Lemma \ref{binorep} and
\eqref{eq:eq2.1} we obtain
\begin{align*}
\omega^{(n)}_{\alpha,\beta}(s_1,t_1,s_2,t_2)=\sum
_{k=0}^{[ns_2]-1}&\sum _{\ell=0}^{[nt_2]-1} (1-\alpha)
\Big ({\mathsf P}
\big(S_{[ns_1]-[ns_2]+k,[nt_1]-[nt_2]+\ell}^{(\alpha,1-\beta)}=
[ns_1]-[ns_2]+k+1\big)\\ &-{\mathsf P}
\big(S_{[ns_1]-[ns_2]+k,[nt_1]-[nt_2]+\ell}^{(\alpha,1-\beta)}=
[ns_1]-[ns_2]+k\big)\Big) {\mathsf P}\big(S_{k,\ell}^{(\alpha,1-\beta)}=k\big).
\end{align*}
Hence, one can use the local versions of the CLT given in Lemmas
\ref{lclt} and \ref{lcltdiff} yielding approximation
\begin{align*}
\omega^{(n)}_{\alpha,\beta}(s_1,t_1,s_2,t_2)\approx
\widetilde E^{(n)}_{\alpha,\beta}(s_1,t_1,s_2,t_2):=& -\frac
{1-\alpha}{2\pi} \sum_{k=1}^{[ns_2]-1}
\sum_{\ell=1}^{[nt_2]-1} f(b_{k,\ell},a_{k,\ell}) \\
=&-\frac {1-\alpha}{2\pi} \int\limits_1^{[ns_2]}
\int\limits_1^{[nt_2]} f(b_{[y],[z]},a_{[y],[z]})\,{\mathrm d}z\, {\mathrm d}y,
\end{align*}
where 
\begin{equation*}
f(u,v):=\frac {v+g_{\alpha,\beta}}{u^{1/2}(u+q_{\alpha,\beta})^{3/2}}
\exp\Big(-\frac {v^2}{2u}\Big) \exp\Big(-\frac
  {(v+g_{\alpha,\beta})^2}{2(u+q_{\alpha,\beta})}\Big),
\end{equation*}
and
\begin{alignat*}{2}
b_{k,\ell}&:=\alpha(1\!-\!\alpha)k+\beta(1\!-\!\beta)\ell, \quad 
&&q_{\alpha,\beta}:=\alpha(1\!-\!\alpha)\big([ns_1]-[ns_2]\big)
+\beta(1\!-\!\beta)\big([nt_1]-[nt_2]\big)\!\geq\! 0, \\ 
a_{k,\ell}&:=(1\!-\!\alpha)k-(1\!-\!\beta)\ell, \quad 
&&g_{\alpha,\beta}:=(1\!-\!\alpha)\big([ns_1]-[ns_2]\big)
-(1\!-\!\beta)\big([nt_1]-[nt_2]\big).
\end{alignat*}
Using  Lemmas \ref{lclt} and \ref{lcltdiff}, as for \ $z\geq0$ \ we
have \ $z\exp(-z)\leq 1$, \ direct calculations show
that for the error 
\begin{equation*}
\widetilde \Delta_{\alpha,\beta}^{(n)}(s_1,t_1,s_2,t_2):=
\omega^{(n)}_{\alpha,\beta}(s_1,t_1,s_2,t_2) -\widetilde
E^{(n)}_{\alpha,\beta}(s_1,t_1,s_2,t_2) 
\end{equation*}
we have
\begin{equation*}
\big |\widetilde \Delta_{\alpha,\beta}^{(n)}(s_1,t_1,s_2,t_2) \big
|\leq C_{\alpha,\beta}\Big (
\widetilde \Delta_{\alpha,\beta}^{(n,1)}(s_1,t_1,s_2,t_2) +\widetilde
\Delta_{\alpha,\beta}^{(n,2)}(s_1,t_1,s_2,t_2)+ +\widetilde
\Delta_{\alpha,\beta}^{(n,3)}(s_1,t_1,s_2,t_2)\Big ), 
\end{equation*}
where \ $C_{\alpha,\beta}$ \ is a positive constant and
\begin{align*}
\widetilde\Delta_{\alpha,\beta}^{(n,1)}(s_1,t_1,s_2,t_2)&:=\sum
_{k=1}^{[ns_2]-1}\sum _{\ell=1}^{[nt_2]-1} \frac
1{b_{k,\ell}^{5/2}}+\sum _{k=1}^{[ns_2]-1}\Big (\frac
1{b_{k,1}^2}+\frac {\alpha ^k}{b_{k,0}}\Big) + 
\sum _{\ell=1}^{[nt_2]-1}\Big (\frac 1{b_{1,\ell}^2}+\frac {\beta 
  ^{\ell}}{b_{0,\ell}}\Big) +1,\\
\widetilde\Delta_{\alpha,\beta}^{(n,2)}(s_1,t_1,s_2,t_2)&:=\sum
_{k=2}^{[ns_2]-1}\sum _{\ell=2}^{[nt_2]-1} \frac 1{b_{k,\ell}^2}\exp
  \Big(-\frac {a_{k,\ell}^2}{2b_{k,\ell}}\Big), \\
\widetilde\Delta_{\alpha,\beta}^{(n,3)}(s_1,t_1,s_2,t_2)&:=\sum
_{k=2}^{[ns_2]-1}\sum _{\ell=2}^{[nt_2]-1} \frac
1{b_{k,\ell}(b_{k,\ell}+q_{\alpha,\beta})}  \exp\Big(-\frac
  {(a_{k,\ell}+g_{\alpha,\beta})^2}{4(b_{k,\ell}+q_{\alpha,\beta})}\Big).
\end{align*}
Obviously, as e.g.
\begin{equation*}
\sum_{k=1}^{[ns_2]-1}\sum _{\ell=1}^{[nt_2]-1} \frac
1{b_{k,\ell}^{5/2}}\leq
  2\big(\alpha\beta(1-\alpha)(1-\beta)\big)^{-3/2}+2D_{\alpha,\beta}
  \quad \text{with} \quad D_{\alpha,\beta}:= 
  \big(\alpha\beta(1-\alpha)(1-\beta)b_{1,1}^{1/2}\big)^{-1}, 
\end{equation*}
it is not difficult to show that \
$\widetilde\Delta_{\alpha,\beta}^{(n,1)}(s_1,t_1,s_2,t_2)$ \ is 
bounded 
from above with an upper bound not depending on \ $s_1,t_1,s_2,t_2$ \
and \ $n$. \
Further, 
\begin{align}
   \label{eq:eq2.3}
\widetilde\Delta_{\alpha,\beta}^{(n,2)}(s_1,t_1,s_2,t_2)&
 \leq 4\int\limits_1^{[ns_2]}
\int\limits_1^{[nt_2]} \frac 1{b_{y,z}^2}
\exp \Big (- \frac{a_{y,z}^2}{2b_{y,z}}\Big){\mathrm d}z\,{\mathrm d}y \\
&\leq \frac 4{(\alpha +\beta)(1-\alpha )(1-\beta )}
\int\limits_{b_{1,1}}^{b_{[ns_2],[nt_2]}}  
\int\limits_{a_{1,[nt_2]}}^{a_{[ns_2],1}} \frac 1{u^2}
\exp \Big (- \frac{v^2}{2u}\Big){\mathrm d}v\,{\mathrm d}u. \nonumber
\end{align}
As for some real constants \ $r$ \ and \ $m>0$
\begin{equation}
   \label{eq:eq2.4}
\int\ \exp \Big (-
\frac{(v-r)^2}{mu}\Big){\mathrm d}v=\frac {\sqrt{\pi mu}}2 
\widetilde\Phi \bigg (\frac {v-r}{\sqrt{mu}} \bigg) 
\end{equation}
holds, where  \ $\widetilde\Phi(x):=2\Phi (\sqrt {2}x)-1, \ x\in{\mathbb
  R},$ \ is the Gauss error function defined with the help of the cdf \
$\Phi (x)$ \ of the standard normal distribution,  using
\eqref{eq:eq2.3} and \eqref{eq:eq2.4} with \ $m=2$ \ and \ $r=0$ \ we have 
\begin{equation*}
\widetilde\Delta_{\alpha,\beta}^{(n,2)}(s_1,t_1,s_2,t_2)
\leq 8\sqrt{2\pi}H_{\alpha,\beta}<\infty,
\quad \text{where} \quad H_{\alpha,\beta}:= 
  \big((\alpha+\beta)(1-\alpha)(1-\beta)b_{1,1}^{1/2}\big)^{-1}.
\end{equation*}
Using the same ideas we also obtain
\begin{align*}
\widetilde\Delta_{\alpha,\beta}^{(n,3)}(s_1,t_1,s_2,t_2)&\leq 
\frac 4{(\alpha \!+\!\beta)(1\!-\!\alpha )(1\!-\!\beta )} \!\!\!\!
\int\limits_{b_{1,1}}^{b_{[ns_2],[nt_2]}}  
\int\limits_{a_{1,[nt_2]}}^{a_{[ns_2],1}}\!\!\!
\frac
1{(u+q_{\alpha,\beta})u}
\exp \Big(-\frac
  {(v+g_{\alpha,\beta})^2}{4(u+q_{\alpha,\beta})}\Big) {\mathrm
    d}v\,{\mathrm d}u \\ &\leq
16\sqrt{\pi}H_{\alpha,\beta}<\infty,
\end{align*}
so there exists a constant \ $0<K_{\alpha,\beta}^{(1)}<\infty$ \ not
depending on \ $s_1,t_1,s_2,t_2$ \ and \ $n$ \ such that
\begin{equation}
  \label{eq:eq2.5}
\big |\widetilde \Delta_{\alpha,\beta}^{(n)}(s_1,t_1,s_2,t_2) \big
|\leq K_{\alpha,\beta}^{(1)} .
\end{equation}

Further, let
\begin{equation*}
\Delta_{\alpha,\beta}^{(n)}(s_1,t_1,s_2,t_2):=
\widetilde E^{(n)}_{\alpha,\beta}(s_1,t_1,s_2,t_2)-
E^{(n)}_{\alpha,\beta}(s_1,t_1,s_2,t_2),  
\end{equation*}
where
\begin{align*}
E^{(n)}_{\alpha,\beta}(s_1,t_1,s_2,t_2):=-\frac {1-\alpha}{2\pi}\int\limits_1^{[ns_2]}
\int\limits_1^{[nt_2]} f(b_{y,z},a_{y,z})\,{\mathrm d}z\, {\mathrm d}y.
\end{align*}
Obviously,
\begin{equation*}
\big| \Delta_{\alpha,\beta}^{(n)}(s_1,t_1,s_2,t_2)\big| \leq \frac
{1-\alpha}{2\pi}\Big( 
\Delta_{\alpha,\beta}^{(n,1)}(s_1,t_1,s_2,t_2)
+\Delta_{\alpha,\beta}^{(n,2)}(s_1,t_1,s_2,t_2)\Big), 
\end{equation*}
where
\begin{align*}
\Delta ^{(n,1)}_{\alpha,\beta}(s_1,t_1,s_2,t_2)&:=\int\limits_1^{[ns_2]} 
\int\limits_1^{[nt_2]}
\big|f(b_{[y],[z]},a_{[y],[z]})-f(b_{[y],[z]},a_{y,z})\big|\,{\mathrm 
  d}z\, {\mathrm d}y, \\ 
\Delta ^{(n,2)}_{\alpha,\beta}(s_1,t_1,s_2,t_2)&:= \int\limits_1^{[ns_2]}
\int\limits_1^{[nt_2]} \big|f(b_{[y],[z]},a_{y,z})-f(b_{y,z},a_{y,z})\big|\,{\mathrm
  d}z\, {\mathrm d}y. 
\end{align*}
As \ $\big|z-[z]\big|<1, \ z\in{\mathbb R}$, \ and for \ $z\geq0$ \ we
have \ $z\exp(-z)\leq 1$, \ and \ $|1-\exp(-z)|\leq |z|$, \ while  for
\ $z\geq 1$, \ $z\geq [z]>z/2$ \ holds,  long but straightforward
calculations yield
\begin{equation*}
\Delta_{\alpha,\beta}^{(n,1)}(s_1,t_1,s_2,t_2) \leq
\ 24\int\limits_1^{[ns_2]}
\int\limits_1^{[nt_2]} \frac 1{b_{y,z}^{5/2}} \leq 48D_{\alpha,\beta}.
\end{equation*}
Now, using similar ideas we obtain
\begin{align*}
&\Delta_{\alpha,\beta}^{(n,2)}(s_1,t_1,s_2,t_2) \!\leq\! \int\limits_1^{[ns_2]}
\int\limits_1^{[nt_2]} \frac 2{b_{y,z}^2}\exp \Big (\!-
\frac{a_{[y],[z]}^2}{2b_{y,z}}\Big){\mathrm d}z\,{\mathrm d}y 
+ \!\!\int\limits_1^{[ns_2]}
\int\limits_1^{[nt_2]}\! \frac {|a_{y,z}\!+\!g_{\alpha,\beta}|\land
  |a_{[y],[z]}\!+\!g_{\alpha,\beta}|}{b_{y,z}^{1/2}(b_{y,z}\!+\!q_{\alpha,\beta})^{3/2}} \\
&\phantom{=}\times \bigg|\exp \Big (\!- \frac{a_{y,z}^2}{2b_{y,z}}\Big)
\exp \Big (\!- 
\frac{(a_{y,z}\!+\!g_{\alpha,\beta})^2}{2(b_{y,z}\!+\!q_{\alpha,\beta})}\Big) 
-\exp \Big (\!-
\frac{a_{[y],[z]}^2}{2b_{y,z}}\Big)\exp \Big (\!-\frac{(a_{[y],[z]}
  \!+\!g_{\alpha,\beta})^2}{2(b_{y,z}\!+\!q_{\alpha,\beta})}\Big)\!\bigg|{\mathrm 
  d}z\,{\mathrm d}y. 
\end{align*}
Hence,
\begin{align*}
\Delta_{\alpha,\beta}^{(n,2)}(s_1,t_1,s_2,t_2)&\leq
8\ln(2)b_{1,1}^{1/2}D_{\alpha,\beta}+\int\limits_1^{[ns_2]+1}
\int\limits_1^{[nt_2]+1} \frac 8{b_{y,z}^2}\exp \Big (\!-
\frac{a_{[y],[z]}^2}{2b_{y,z}}\Big){\mathrm d}z\,{\mathrm d}y \\
+\int\limits_1^{[ns_2]}
\int\limits_1^{[nt_2]}& \frac {|a_{y,z}\!+\!g_{\alpha,\beta}|\land
  |a_{[y],[z]}\!+\!g_{\alpha,\beta}|}{b_{y,z}^{1/2}(b_{y,z}\!+\!q_{\alpha,\beta})^{3/2}}
\Bigg(\frac {|a_{y,z}^2-a_{[y],[z]}^2|}{2b_{y,z}}  
\exp \Big (\!- \frac{a_{y,z}^2\land a_{[y],[z]}^2}{2b_{y,z}}\Big) \\
&\times
\exp \Big(\!-\frac{(a_{[y],[z]}
  \!+\!g_{\alpha,\beta})^2}{2(b_{y,z}\!+\!q_{\alpha,\beta})}\Big)
+\frac
{|(a_{y,z}\!+\!g_{\alpha,\beta})^2-
  (a_{[y],[z]}\!+\!g_{\alpha,\beta})^2|}{2(b_{y,z}\!+\!q_{\alpha,\beta})}  
\\ 
&\times 
\exp \Big (\!- \frac{a_{y,z}^2}{2b_{y,z}}\Big)\exp \Big(
\!-\frac{(a_{[y],[z]}
  \!+\!g_{\alpha,\beta})^2\land
  (a_{y,z}\!+\!g_{\alpha,\beta})^2}{2(b_{y,z}\!+\!q_{\alpha,\beta})}\Big)\!\Bigg)  
{\mathrm d}z\,{\mathrm d}y,  \nonumber
\end{align*}
and as \ $\big|a_{y,z}^2-a_{[y],[z]}^2\big|\leq 4+4|a_{y,z}\land a_{[y],[z]}| $, \
short calculation shows
\begin{align*}
\Delta_{\alpha,\beta}^{(n,2)}&(s_1,t_1,s_2,t_2)\leq
\,16\sqrt{2\pi}\big(D_{\alpha,\beta}\!+\!H_{\alpha,\beta}\big)+8
\int\limits_1^{[ns_2]} 
\int\limits_1^{[nt_2]} \bigg(\frac 1{b_{y,z}^{5/2}}+\frac
1{b_{y,z}^2}\exp \Big (- \frac{a_{y,z}^2}{2b_{y,z}}\Big)\bigg){\mathrm
  d}z\,{\mathrm d}y \\
 +&
8\int\limits_1^{[ns_2]} 
\int\limits_1^{[nt_2]}\frac
1{b_{y,z}(b_{y,z}\!+\!q_{\alpha,\beta})}\exp \Big(\!-\frac{(a_{[y],[z]} 
  \!+\!g_{\alpha,\beta})^2}{4(b_{y,z}\!+\!q_{\alpha,\beta})}\Big){\mathrm
  d}z\,{\mathrm d}y \leq
32\sqrt{2\pi}\big(D_{\alpha,\beta}\!+\!H_{\alpha,\beta}\big)\\ 
+&32\!\!\!\int\limits_1^{[ns_2]+1}\! 
\int\limits_1^{[nt_2]+1}\!\!\!\!\frac
1{b_{y,z}^{3/2}(b_{y,z}\!+\!q_{\alpha,\beta})^{1/2}}\exp \Big(\!-\frac{(a_{y,z} 
  \!+\!g_{\alpha,\beta})^2}{4(b_{y,z}\!+\!q_{\alpha,\beta})}\Big){\mathrm
  d}z\,{\mathrm d}y\!\leq\!
80\sqrt{2\pi}\big(D_{\alpha,\beta}\!+\!H_{\alpha,\beta}\big)\!<\!\infty.
\end{align*}
Thus, there exists a constant \ $0<K_{\alpha,\beta}^{(2)}<\infty$ \  such that
\begin{equation}
  \label{eq:eq2.6}
\big |\widetilde \Delta_{\alpha,\beta}^{(n)}(s_1,t_1,s_2,t_2) \big
|\leq K_{\alpha,\beta}^{(2)}.
\end{equation}

Finally, one has to show the existence of a constant \
$0<K_{\alpha,\beta}^{(3)}<\infty$ \ such that   
\begin{equation}
  \label{eq:eq2.7}
\big |E_{\alpha,\beta}^{(n)}(s_1,t_1,s_2,t_2) \big
|\leq K_{\alpha,\beta}^{(3)},
\end{equation}
which together with \eqref{eq:eq2.5} and \eqref{eq:eq2.6} implies the
statement of the proposition.
 
Consider first the case  \ $b_{[ns_2],1}\leq b_{1,[nt_2]}$. In this case 
\begin{equation*}
 E^{(n)}_{\alpha,\beta}(s_1,t_1,s_2,t_2)\!=\! \frac
 {(1\!-\!\alpha)H_{\alpha,\beta}}{2\pi b_{1,1}^{-1/2}}
\Big(\!E^{(n,1)}_{\alpha,\beta}(s_1,t_1,s_2,t_2)\!+\!
E^{(n,2)}_{\alpha,\beta}(s_1,t_1,s_2,t_2)\!+\! 
E^{(n,3)}_{\alpha,\beta}(s_1,t_1,s_2,t_2)\!\Big), 
\end{equation*}
where
\begin{align*}
E^{(n,1)}_{\alpha,\beta}(s_1,t_1,s_2,t_2)&:=
\int\limits_{b_{1,1}}^{b_{[ns_2],1}}
\int\limits_{-u/\beta+(\alpha+\beta)(1-\alpha)/\beta}^{u/\alpha
  -(\alpha+\beta)(1-\beta)/\alpha} 
f(u,v)\,{\mathrm d}v\,{\mathrm
  d}u \\
=\int\limits_{b_{1,1}}^{b_{[ns_2],1}}\!\!\! F\bigg(u,&\,\frac{u-
  (\alpha+\beta)(1-\beta)}{\alpha} \bigg){\mathrm d}u-\!\!\!
\int\limits_{b_{1,1}}^{b_{[ns_2],1}}\!\!\! F\bigg(u,\frac{-u+
  (\alpha+\beta)(1-\alpha)}{\beta} \bigg){\mathrm d}u , \\
E^{(n,2)}_{\alpha,\beta}(s_1,t_1,s_2,t_2)&:=\int\limits_{b_{[ns_2],1}}^{b_{1,[nt_2]}}
\int\limits_{-u/\beta+(\alpha+\beta)(1-\alpha)/\beta}^{-u/\beta+(\alpha+\beta)
  (1-\alpha)[ns_2]/\beta}
f(u,v)\,{\mathrm d}v\,{\mathrm d}u \\
=\int\limits_{b_{[ns_2],1}}^{b_{1,[nt_2]}}\!\!\! F\bigg(u,&\,\frac{-u+
  (\alpha+\beta)(1-\alpha)[ns_2]}{\beta} \bigg){\mathrm d}u-\!\!\!
\int\limits_{b_{[ns_2],1}}^{b_{1,[nt_2]}}\!\!\! F\bigg(u,\frac{-u+
  (\alpha+\beta)(1-\alpha)}{\beta} \bigg){\mathrm d}u ,
\end{align*}
\begin{align*}
E^{(n,3)}_{\alpha,\beta}(s_1,t_1,s_2,t_2)&:=
\int\limits_{b_{1,[nt_2]}}^{b_{[ns_2],[nt_2]}}
\int\limits_{u/\alpha-(\alpha+\beta)(1-\beta)[nt_2]/\alpha}^{-u/\beta+(\alpha+\beta)
  (1-\alpha)[ns_2]/\beta}
f(u,v)\,{\mathrm d}v\,{\mathrm d}u\\
=\int\limits_{b_{1,[nt_2]}}^{b_{[ns_2],[nt_2]}}\!\!\! F\bigg(u,&\,\frac{-u+
  (\alpha+\beta)(1-\alpha)[ns_2]}{\beta} \bigg){\mathrm d}u-\!\!\!
\int\limits_{b_{1,[nt_2]}}^{b_{[nt_2],[nt_2]}}\!\!\! F\bigg(u,\frac{u-
  (\alpha+\beta)(1-\beta)[nt_2]}{\alpha} \bigg){\mathrm d}u,
\end{align*}
with
\begin{align*}
F(u,v):=\int &f(u,v)\,{\mathrm d}v=-\frac
{u^{1/2}}{(u+q_{\alpha,\beta})^{1/2}(2u+q_{\alpha,\beta})}  
\exp \Big (-\frac{v^2}{2u}\Big)
\exp \Big (-\frac{(v+g_{\alpha,\beta})^2}{2(u+q_{\alpha,\beta})}\Big) \\
&+\frac {\sqrt{2\pi }
  g_{\alpha,\beta}}{2 
  (2u+q_{\alpha,\beta})^{3/2}} \exp \Big (-
\frac{g_{\alpha,\beta}^2}{2(2u+q_{\alpha,\beta})}\Big)\widetilde\Phi
\bigg(\frac{(2u+q_{\alpha,\beta})v+g_{\alpha,\beta}u}{(2u
  (2u+q_{\alpha,\beta})(u+q_{\alpha,\beta}))^{1/2}} \bigg).
\end{align*}
Combining similar terms we obtain
\begin{align}
   \label{eq:eq2.8}
E^{(n)}_{\alpha,\beta}&(s_1,t_1,s_2,t_2)=\frac 1{2\pi
  (\alpha\!+\!\beta)(1\!-\!\beta)} \Bigg(\int\limits_{b_{1,1}}^{b_{[ns_2],1}}\!\!\! 
F\bigg(u,\,\frac{u\!-\! 
  (\alpha+\beta)(1-\beta)}{\alpha} \bigg){\mathrm d}u \\
&-\int\limits_{b_{1,1}}^{b_{1,[nt_2]}}\!\!\! F\bigg(u,\,\frac{-u\!+\! 
  (\alpha\!+\!\beta)(1\!-\!\alpha)}{\beta} \bigg){\mathrm d}u 
+\int\limits_{b_{[ns_2],1}}^{b_{[ns_2],[nt_2]}}\!\!\! 
F\bigg(u,\,\frac{-u\!+\! 
  (\alpha\!+\!\beta)(1\!-\!\alpha)[ns_2]}{\beta} \bigg){\mathrm d}u
\nonumber \\  
&-\int\limits_{b_{1,[nt_2]}}^{b_{[ns_2],[nt_2]}}\!\!\! 
F\bigg(u,\,\frac{u\!-\! 
  (\alpha\!+\!\beta)(1\!-\!\beta)[nt_2]}{\alpha}\bigg){\mathrm
  d}u\Bigg). \nonumber 
\end{align}
As for \ $u\geq 0$
\begin{equation}
  \label{eq:eq2.9}
\big|F(u,v)\big|\leq \frac 1{2u}  
\exp \Big (-\frac{v^2}{2u}\Big)+\frac {\sqrt{2\pi }
  |g_{\alpha,\beta}|}{2 
  (2u+q_{\alpha,\beta})^{3/2}} \exp \Big (-
\frac{g_{\alpha,\beta}^2}{2(2u+q_{\alpha,\beta})}\Big),
\end{equation}
using Taylor series expansion we have
\begin{align}
  \label{eq:eq2.10}
\int\limits_{b_{1,1}}^{b_{[ns_2],1}} &\,\Bigg|F\bigg(u,\,\frac{u-
  (\alpha\!+\!\beta)(1\!-\!\beta)}{\alpha} \bigg)\Bigg|{\mathrm d}u
\leq  \!\!\!
\int\limits_{b_{1,1}}^{b_{[ns_2],1}}\!\! \frac 1{2u}\exp \bigg(\!-\frac
12\Big (\frac {u^{1/2}}{\alpha}\!-\!\frac
{(\alpha\!+\!\beta)(1\!-\!\beta)}{\alpha 
  u^{1/2}}\Big)^2\bigg){\mathrm d}u \\
+&\!\!\!\int\limits_{b_{1,1}}^{b_{[ns_2],1}} \!\!\!\frac {\sqrt{2\pi } 
  |g_{\alpha,\beta}|}{2 
  (2u\!+\!q_{\alpha,\beta})^{3/2}} \exp \Big (-\!
\frac{g_{\alpha,\beta}^2}{2(2u\!+\!q_{\alpha,\beta})}\Big){\mathrm d}u 
\leq \!\!\! \int\limits_{b_{1,1}}^{b_{[ns_2],1}} \!\!\bigg (\frac 1{2u}\exp
\Big(\!-\frac u{2\alpha^2}\Big)\!+\!\frac
{(\alpha\!+\!\beta)(1\!-\!\beta)}{2\alpha 
  u^{3/2}} \bigg){\mathrm d}u \nonumber  \\
+&\sqrt{\frac
  {\pi}2}\int\limits_{|g_{\alpha,\beta}|/(2b_{[ns_2],1}
  +q_{\alpha,\beta})^{1/2}}^{|  
  g_{\alpha,\beta}|/(2b_{1,1}+q_{\alpha,\beta})^{1/2}}
\!\!\!\!\!\!\!\!\!\exp \big 
(\!-w^2/2\big) {\mathrm d}w \leq \frac
{(\alpha+\beta)(1-\beta)}{\alpha b_{1,1}^{1/2}}+\frac 12E_1\Big(\frac
{b_{1,1}}{2\alpha ^2}\Big) \nonumber\\
+&\frac {\pi}2\Bigg (\widetilde \Phi \bigg(
\frac{|g_{\alpha,\beta}|}{(2b_{1,1}\!+\!q_{\alpha,\beta})^{1/2}}\bigg) 
\!-\!\widetilde \Phi \bigg(
\frac{|g_{\alpha,\beta}|}{(2b_{[ns_2],1}+q_{\alpha,\beta})^{1/2}}\bigg)\Bigg) 
\leq \frac
{(\alpha\!+\!\beta)(1\!-\!\beta)}{\alpha b_{1,1}^{1/2}}\!+\!\frac
12E_1\Big(\frac 
{b_{1,1}}{2\alpha ^2}\Big)\!+\!\pi, \nonumber  
\end{align} 
where \ $E_1(x)$ \ is the exponential integral \citep{as}
\begin{equation*}
\frac 12 {\mathrm e}^{-x}\log \Big
(1+\frac 2x\Big)<E_1(x):=\int\limits_x^{\infty}\frac {{\mathrm
    e}^{-u}}u{\mathrm d}u<  
{\mathrm e}^{-x}\log \Big (1+\frac 1x\Big), \qquad x>0.
\end{equation*}
Since the right hand side of \eqref{eq:eq2.9} is an upper bound for \
$|F(u,-v)|$ \ as well, same calculations as in  \eqref{eq:eq2.10}
yield
\begin{equation}
  \label{eq:eq2.11}
\int\limits_{b_{1,1}}^{b_{[ns_2],1}} \Bigg|F\bigg(u,\,\frac{-u+
  (\alpha\!+\!\beta)(1\!-\!\alpha)}{\beta} \bigg)\Bigg|{\mathrm d}u
\leq\frac
{(\alpha\!+\!\beta)(1\!-\!\alpha)}{\beta b_{1,1}^{1/2}}\!+\!\frac
12E_1\Big(\frac 
{b_{1,1}}{2\beta ^2}\Big)\!+\!\pi.
\end{equation}
Now, using again \eqref{eq:eq2.9}, as  \ $b_{1,1}\leq b_{[ns_2],1}$, \ we have 
\begin{equation}
  \label{eq:eq2.12}
\int\limits_{b_{[ns_2],1}}^{b_{[ns_2],[nt_2]}} \Bigg|F\bigg(u,\,\frac{-u+
  (\alpha\!+\!\beta)(1\!-\!\alpha)[ns_2]}{\beta} \bigg)\Bigg|{\mathrm d}u 
\leq\frac
{(\alpha\!+\!\beta)(1\!-\!\alpha)}{\beta b_{1,1}^{1/2}}+\pi+{\mathcal
  E}^{(n)}_{\alpha , \beta}(s_2,t_2), 
\end{equation}
where
\begin{equation*}
{\mathcal E}^{(n)}_{\alpha,\beta}(s_2,t_2):=
\int\limits_{b_{[ns_2],1}}^{b_{[ns_2],[nt_2]}}\frac 
1{2u}\exp \bigg(\!- 
\frac {(u\!-\!(\alpha\!+\!\beta)(1\!-\!\alpha)[ns_2])^2}{2\beta^2 
  u}\bigg){\mathrm d}u. 
\end{equation*}
Assume first \ $(1-\beta)[nt_2]\leq (1-\alpha)[ns_2]$ \ implying \
$b_{[ns_2],[nt_2]} \leq (\alpha+\beta)(1-\alpha)[ns_2]$. \ In
this case
\begin{align}
{\mathcal E}^{(n)}_{\alpha,\beta}(s_2,t_2)&\leq \frac {(2\pi
  b_{[ns_2],[nt_2]})^{1/2}\beta}{2b_{[ns_2],1}} \!\!
\int\limits_{b_{[ns_2],1}}^{b_{[ns_2],[nt_2]}}\!\! \frac  
1{\sqrt{2\pi b_{[ns_2],[nt_2]}}\beta}\exp \bigg(\!- 
\frac {(u\!-\!(\alpha\!+\!\beta)(1\!-\!\alpha)[ns_2])^2}{2b_{[ns_2],[nt_2]}\beta^2 
  }\bigg){\mathrm d}u \nonumber \\
&\leq \frac{(2\pi(\alpha\!+\!\beta)(1\!-\!\alpha)[ns_2])^{1/2}
  \beta}{2\alpha(1\!-\!\alpha)[ns_2]} \leq \frac {2\beta}{\alpha}\bigg(\frac {
  \alpha\!+\!\beta}{1\!-\!\alpha}\bigg)^{1/2} \leq \frac
2{\alpha(1\!-\!\alpha)^{1/2}}. \label{eq:eq2.13} 
\end{align}
Further, assume \ $(1-\beta)[nt_2]> (1-\alpha)[ns_2]>1-\beta$ \
implying \ $b_{[ns_2],[nt_2]} > (\alpha+\beta)(1-\alpha)[ns_2]
>b_{[ns_2],1}$. \ Thus,
\begin{align}
{\mathcal E}^{(n)}_{\alpha,\beta}(s_2,t_2)\leq &
\int\limits_{b_{[ns_2],1}}^{(\alpha+\beta)(1-\alpha)[ns_2]}\frac  
1{2b_{[ns_2],1}}\exp \bigg(\!- 
\frac {(u\!-\!(\alpha\!+\!\beta)(1\!-\!\alpha)[ns_2])^2}{2\beta^2 
  (\alpha+\beta)(1-\alpha)[ns_2]}\bigg){\mathrm d}u \nonumber \\
&+\int\limits_{(\alpha+\beta)(1-\alpha)[ns_2]}^{b_{[ns_2],[nt_2]}}\frac   
1{2u}\exp \bigg(\!- \frac 1{2\beta^2}\Big(u^{1/2}\!-\!\frac
{(\alpha\!+\!\beta)(1\!-\!\alpha)[ns_2]}{
  u^{1/2}}\Big)^2\bigg){\mathrm d}u  \label{eq:eq2.14} \\
\leq &\, \frac {2\beta}{\alpha}\bigg(\frac {
  \alpha\!+\!\beta}{1\!-\!\alpha}\bigg)^{1/2}\!\!\!+
\int\limits_{\sqrt{(\alpha+\beta)(1-\alpha)[ns_2]}}^{\sqrt{b_{[ns_2],[nt_2]}}}\frac   
1w\exp\bigg(\!-\frac{\big(w-\sqrt{(\alpha\!+\!\beta)
    (1\!-\!\alpha)[ns_2]}\big)^2}{2\beta^2}\bigg){\mathrm d}w
\nonumber \\
\leq &\, \frac 2{\alpha(1\!-\!\alpha)^{1/2}}+\bigg(\frac
{2\pi\beta^2}{(\alpha\!+\!\beta)(1\!-\!\alpha)[ns_2]}\bigg)^{1/2} \!\!\!
\leq \frac 2{\alpha(1\!-\!\alpha)^{1/2}}+\frac
{4\beta}{\big((\alpha\!+\!\beta)(1\!-\!\alpha)\big)^{1/2}}. \nonumber 
\end{align}
Finally, suppose \ $(1-\alpha)[ns_2]\leq 1-\beta$ \
implying \ $(\alpha+\beta)(1-\alpha)[ns_2] \leq b_{[ns_2],1}$. \ In
this case
\begin{equation*}
{\mathcal E}^{(n)}_{\alpha,\beta}(s_2,t_2)\leq  \frac 1{b_{[ns_2],1}^{1/2}}\!\!
\int\limits_{\sqrt{b_{[ns_2],1}}}^{\sqrt{b_{[ns_2],[nt_2]}}}\!\!\!
\exp \bigg(\!- \frac 1{2\beta^2}\Big(w-\frac
{(\alpha\!+\!\beta)(1\!-
  \!\alpha)[ns_2]}{\sqrt{b_{[ns_2],1}}}\Big)^2\bigg){\mathrm d}w 
\leq \frac {4\beta}{b_{1,1}^{1/2}}\leq \frac 4{\alpha(1\!-\!\alpha)^{1/2}}, 
\end{equation*}
that together with \eqref{eq:eq2.12} -- \eqref{eq:eq2.14} yields
\begin{equation}
  \label{eq:eq2.15}
\int\limits_{b_{[ns_2],1}}^{b_{[ns_2],[nt_2]}} \Bigg|F\bigg(u,\,\frac{-u\!+\!
  (\alpha\!+\!\beta)(1\!-\!\alpha)[ns_2]}{\beta} \bigg)\Bigg|{\mathrm d}u 
\leq\frac 7{\alpha\beta(1\!-\!\alpha)^{1/2}}+\frac
4{\big((\alpha\!+\!\beta)(1\!-\!\alpha))^{1/2}}+\pi.
\end{equation}
Using similar arguments one can easily show
\begin{equation}
  \label{eq:eq2.16}
\int\limits_{b_{1,[nt_2]}}^{b_{[ns_2],[nt_2]}} \Bigg|F\bigg(u,\,\frac{u\!-\!
  (\alpha\!+\!\beta)(1\!-\!\beta)[nt_2]}{\alpha} \bigg)\Bigg|{\mathrm d}u 
\leq\frac 7{\alpha\beta(1\!-\!\beta)^{1/2}}+\frac
4{\big((\alpha\!+\!\beta)(1\!-\!\beta))^{1/2}}+\pi.
\end{equation}
Now, representation \eqref{eq:eq2.8} together with \eqref{eq:eq2.10},
\eqref{eq:eq2.11}, \eqref{eq:eq2.15} and \eqref{eq:eq2.16} directly
implies \eqref{eq:eq2.7} that completes the proof in the case \
$[ns_1]\geq [ns_2]$,  \ $[nt_1] \geq [nt_2]$ \ and \ $b_{[ns_2],1}
\leq b_{1,[nt_2]}$.  

Further, if \ $b_{[ns_2],1} >b_{1,[nt_2]}$ \ representation
\eqref{eq:eq2.8} also holds for \
$E^{(n)}_{\alpha,\beta}(s_1,t_1,s_2,t_2)$ \  directly implying 
\eqref{eq:eq2.7}. 

By symmetry, case \ $[ns_1] < [ns_2]$,  \ $[nt_1] < [nt_2]$ \ can be
handled in the same way as case \ $[ns_1]\geq [ns_2]$,  \ $[nt_1] \geq
[nt_2]$, \ while in case 
\ $[ns_1] < [ns_2]$,  \ $[nt_1] \geq [nt_2]$ \ we have 
\begin{align*}
\omega^{(n)}_{\alpha,\beta}(s_1,t_1,s_2,t_2)=\sum
_{k=0}^{[ns_1]-1}\sum _{\ell=0}^{[nt_2]-1} (1-\alpha)&
\Big ({\mathsf P}
\big(S_{k,[nt_1]-[nt_2]+\ell}^{(\alpha,1-\beta)}=k+1\big)-{\mathsf P}
\big(S_{k,[nt_1]-[nt_2]+\ell}^{(\alpha,1-\beta)}=k\big)\Big) \\ &\times{\mathsf
  P}\big(S_{[ns_2]-[ns_1]+k,\ell}^{(\alpha,1-\beta)}=[ns_2]-[ns_1]+k\big)
+\sum_{\ell=0}^{[nt_2]-1}\beta^{\ell}. 
\end{align*}
Thus, local versions of the CLT given in Lemmas
\ref{lclt} and \ref{lcltdiff} yield approximation 
\begin{equation*}
\omega^{(n)}_{\alpha,\beta}(s_1,t_1,s_2,t_2)\approx
E^{(n)}_{\alpha,\beta}(s_1,t_1,s_2,t_2):=-\frac {1-\alpha}{2\pi}
\int\limits_1^{[ns_1]} 
\int\limits_1^{[nt_2]} f(b_{y,z},a_{y,z})\,{\mathrm d}z\, {\mathrm d}y,
\end{equation*}
where now
\begin{equation}
   \label{eq:eq2.17}
f(u,v):=\frac {v-g_{2,\alpha,\beta}}{(u\!+\!\alpha g_{1,\alpha,\beta})^{1/2}
  (u\!+\!\beta g_{2,\alpha,\beta})^{3/2}} 
\exp\Big(\!-\frac
{(v\!+\!g_{1,\alpha,\beta})^2}{2(u\!+\!\alpha g_{1,\alpha,\beta})}\Big) 
\exp\Big(\!-\frac {(v\!-\!g_{2,\alpha,\beta})^2}{2(u\! +\!\beta
  g_{2,\alpha,\beta})}\Big), 
\end{equation}
and
\begin{equation*}
g_{1,\alpha,\beta}:=(1\!-\!\alpha)\big([ns_2]-[ns_1]\big), \qquad 
g_{2,\alpha,\beta}:=(1\!-\!\beta)\big([nt_1]-[nt_2]\big).
\end{equation*}
Using similar ideas as in case \ $[ns_1]\geq [ns_2]$,  \ $[nt_1] \geq
[nt_2]$ \  one can show that the error of the approximation is bounded
with a bound not depending on \ $s_1,t_1,s_2,t_2$ \ and \ $n$. 

Substituting \ $[ns_1]$ \ for  \ $[ns_2]$ \ in \eqref{eq:eq2.8} we
obtain a representation of \ $E^{(n)}_{\alpha,\beta}(s_1,t_1,s_2,t_2)$
\ with the help of integrals of \ $F(u,v):=\int f(u,v)\,{\mathrm d}v$
\ where \ $f(u,v)$ \ is defined by \eqref{eq:eq2.17}. Further, for \ $u\geq
0$ \ we have
\begin{align*}
\big |F(u,v)\big | \leq &\,
\frac 1{(u\!+\!\alpha g_{1,\alpha,\beta})^{1/2}
  (u\!+\!\beta g_{2,\alpha,\beta})^{1/2}} 
\exp\Big(\!-\frac
{(v\!+\!g_{1,\alpha,\beta})^2}{2(u\!+\!\alpha g_{1,\alpha,\beta})}\Big) 
\exp\Big(\!-\frac {(v\!-\!g_{2,\alpha,\beta})^2}{2(u\! +\!\beta
  g_{2,\alpha,\beta})}\Big)
\\\ &+\frac
{\sqrt{2\pi }|g_{1,\alpha,\beta}+g_{2,\alpha,\beta}|}{2  
  (2u+q_{1,\alpha,\beta}+q_{2,\alpha,\beta})^{3/2}} \exp \Big (-
\frac{(g_{1,\alpha,\beta}+g_{2,\alpha,\beta})^2}{2(2u+q_{1,\alpha,\beta}+
q_{2,\alpha,\beta})}\Big),
\end{align*}
that is the analogue of \eqref{eq:eq2.9}. Hence, after long but
straightforward calculations similar to the proofs of
\eqref{eq:eq2.10}  and \eqref{eq:eq2.15} one can verify \eqref{eq:eq2.7} that
completes the proof in case \ $[ns_1] < [ns_2]$,  \ $[nt_1] \geq
[nt_2]$. \ Finally, case \ $[ns_1] \geq [ns_2]$,  \ $[nt_1] < [nt_2]$
\ follows by symmetry. \proofend

\begin{Pro}
   \label{covlim}
\ Let \ $0< s_1,t_1,s_2,t_2 \in {\mathbb R}$  \ and let \ 
$(q_1,q_2),(r_1,r_2)\in \big\{ (0,1),(1,0),(0,0)\big \}$. 

If \ $(\alpha, \beta, \gamma) \in {\mathcal F}_{++}$ \ then
\begin{equation*}
\frac 1{n^{1/2}}\Cov
\big(Y^{(n)}_{q_1,q_2}(s_1,t_1),Y^{(n)}_{r_1,r_2}(s_2,t_2)\big)\!\to\!
\frac {\left((1-\alpha)s_1\right)^{1/2}\land
  \left((1-\beta)t_1\right)^{1/2}}{\pi^{1/2}(\alpha+\beta)^{1/2}
  (1-\alpha)(1-\beta)}, \quad \text{if  $s_1\!=\!s_2$  and 
  $t_1\!=\!t_2$,}   
\end{equation*}
otherwise, if \  $(1-\alpha)(s_1-s_2)\ne (1-\beta)(t_1-t_2)$ \ it
tends to \ $0$, as \ $n\to\infty$. \
Moreover, convergence to \ $0$ \ has an exponential rate.

If \ $(\alpha, \beta, \gamma) \in {\mathcal E}_{1+}$  \ then  
\begin{equation*}
\lim_{n\to\infty}\frac 1n \Cov
\big(Y^{(n)}_{q_1,q_2}(s_1,t_1),Y^{(n)}_{r_1,r_2}(s_2,t_2)\big)=
\begin{cases}
(s_1\land s_2)\beta ^{|q_2-r_2|}(1-\gamma ^2)^{-1} , & \text{if \
  $t_1=t_2$,} \\
0, &\text{otherwise,}
\end{cases}
\end{equation*}
while for \ $(\alpha, \beta, \gamma) \in {\mathcal E}_{2+}$ \ we have
\begin{equation*}
\lim_{n\to\infty}\frac 1n \Cov
\big(Y^{(n)}_{q_1,q_2}(s_1,t_1),Y^{(n)}_{r_1,r_2}(s_2,t_2)\big)=
\begin{cases}
(t_1\land t_2)\alpha ^{|q_1-r_1|}(1-\gamma ^2)^{-1}, & \text{if \
  $s_1=s_2$,} \\
0, &\text{otherwise.}
\end{cases}
\end{equation*}
Moreover, convergences to \ $0$ \ in both cases have exponential rates.

If \ $(\alpha, \beta, \gamma) \in {\mathcal V}_+$ \ then
\begin{equation*}
\lim_{n\to\infty} \frac 1{n^2}\Cov
\big(Y^{(n)}_{q_1,q_2}(s_1,t_1),Y^{(n)}_{r_1,r_2}(s_2,t_2)\big)=(s_1\land
s_2)(t_1\land t_2). 
\end{equation*}
\end{Pro}

The proof of the above Proposition is strongly based on the following
Lemma that is an obvious generalization of Theorem 2.3 of
\citet{bpz2}. The statement of the Lemma can be obtained from
Hoeffding's inequality \citep{hoef}.

\begin{Lem}
   \label{expbound}
Using notations of Lemma \ref{lclt} let
\begin{equation*}
\theta := \frac{\nu k+\mu \ell}{k+\ell} \quad \text{and} \quad 
I_{\theta}(x):=\begin{cases}
              x\log \frac x{\theta}+(1-x)\log \frac {1-x}{1-\theta}, &
              \ x\in [0,1], \\
              \infty, & \text{otherwise}.
              \end{cases} 
\end{equation*}
Then for \ $x\ne \theta$ \ we have \
$I_{\theta}(x)>I_{\theta}(\theta)=0$, \ and
\begin{align*}
&{\mathrm P}\big(S_{k,\ell}^{(\nu,\mu)}\geq (k+\ell)x \big)\leq \exp
\big(-(k+\ell)I_{\theta}(x) \big), \quad \text{for all \
  $x>\theta$,}\\
&{\mathrm P}\big(S_{k,\ell}^{(\nu,\mu)}\leq (k+\ell)x \big)\leq \exp
\big(-(k+\ell)I_{\theta}(x) \big), \quad \text{for all \ $x<\theta$.}
\end{align*}
\end{Lem}

\noindent
\textbf{Proof of Proposition \ref{covlim}} \
Let \ $(\alpha, \beta, \gamma) \in {\mathcal F}_{++}$ \
and let \ 
$0<s,t\in{\mathbb R}$. \ By Theorem 1.1 of \citet{baran} we have
\begin{equation*}
\lim_{n\to\infty}\frac 1{n^{1/2}}\Var
\big(Y^{(n)}_{0,0}(s,t)\big) = \frac
{\left((1-\alpha)s\right)^{1/2}\land 
  \left((1-\beta)t\right)^{1/2}}{\pi^{1/2}(\alpha+\beta)^{1/2}
  (1-\alpha)(1-\beta)}.
\end{equation*}
Now, as e.g.
\begin{align*}
\Cov\big(Y^{(n)}_{1,0}(s,t),Y^{(n)}_{0,1}(s,t)\big)\!=\!
\Cov\big(Y^{(n)}_{0,1}(s,t),Y^{(n)}_{0,0}(s\!+\!1/n,t)\big)\!-\! 
\Cov\big(Y^{(n)}_{0,0}(s,t),Y^{(n)}_{0,0}(s\!+\!1/n,t)\big)&\\
+\Cov\big(Y^{(n)}_{1,0}(s,t),Y^{(n)}_{0,0}(s,t)\big)-
\Var\big(Y^{(n)}_{0,0}(s,t),Y^{(n)}_{0,0}(s,t)\big)+\Var
\big(Y^{(n)}_{0,0}(s,t)\big)&,
\end{align*}
by Proposition \ref{covdiff} the limit of \ $n^{-1/2}\Cov
\big(Y^{(n)}_{q_1,q_2}(s_1,t_1),Y^{(n)}_{r_1,r_2}(s_2,t_2)\big)$ \ as
\ $n\to\infty$ \
equals the limit of \ $n^{1/2}\Var 
\big(Y^{(n)}_{0,0}(s,t)\big)$  \ for all \ $(q_1,q_2),(r_1,r_2)\in
\big\{ (0,1),(1,0),(0,0)\big \}$. \ 

Assume first \ $s_1>s_2, \ t_1>t_2$ \ that implies \ $[ns_1]+q_1\geq
[ns_2]+r_1$ \ and \ $[nt_1]+q_2\geq [nt_2]+r_2$ \ if \ $n\in{\mathbb
  N}$ \ is large 
enough. In this case Lemma \ref{binorep} and \eqref{eq:eq2.1} imply
\begin{equation*}
\Cov
\big(Y^{(n)}_{q_1,q_2}(s_1,t_1),Y^{(n)}_{r_1,r_2}(s_2,t_2)\big)=\!\! 
\sum_{k=0}^{[ns_2]+r_1-1}
\sum_{\ell =0}^{[nt_2]+r_2-1}\!\! {\mathrm
  P}\big(S_{k,\ell}^{(\alpha,1-\beta)}=k\big){\mathrm
  P}\big(S_{g_{1,n}+k,g_{2,n}+\ell}^{(\alpha,1-\beta)}=g_{1,n}+k\big),
\end{equation*}
where
\begin{equation*}
g_{1,n}:=\big |[ns_1]-[ns_2]+q_1-r_1 \big|, \qquad
g_{2,n}:=\big |[nt_1]-[nt_2]+q_2-r_2 \big|.
\end{equation*}
We are going to apply Lemma \ref{expbound} for the terms of the above
sum. Let
\begin{align*}
\theta_n:=&\,\frac{\alpha(g_{1,n}+k)+(1-\beta)
    (g_{2,n}+\ell)}{g_{1,n}+k+g_{2,n}+\ell}  
\to \frac{\alpha(s_1-s_2)+(1-\beta) (t_1-t_2)}{s_1-s_2+t_1-t_2}=:\theta, \\
\omega_n:=&\,\frac {g_{1,n}+k}{g_{1,n}+k+g_{2,n}+\ell}-\theta_n 
\to \frac{(1-\alpha)(s_1-s_2)-(1-\beta)
  (t_1-t_2)}{s_1-s_2+t_1-t_2}=:\omega,
\end{align*}
as \ $n\to\infty$. \ If \ $(1-\alpha)(s_1-s_2)>(1-\beta)(t_1-t_2)$ \
then \ $\omega>0$. \ 
Hence, for sufficiently large \ $n\in{\mathbb N}$ \ we have \
$\omega_n\geq \omega/2 >0$ \ and in this way 
\begin{equation*}
{\mathrm P}\big(S_{g_{1,n}+k,g_{2,n}+\ell}^{(\alpha,1-\beta)}=
g_{1,n}+k\big)\leq{\mathrm
  P}\big(S_{g_{1,n}+k,g_{2,n}+\ell}^{(\alpha,1-\beta)}\geq
(g_{1,n}+k+g_{2,n}+\ell)(\theta_n+\omega/2)\big)
\end{equation*}
for all \ $k\in\{0,\ldots,[ns_2]+r_1-1\}$ \ and \
$\ell\in\{0,\ldots,[nt_2]+r_2-1\}$.  \  Further,  for sufficiently
large \ $n\in{\mathbb N}$ \ and for all \
$k\in\{0,\ldots,[ns_2]+r_1-1\}$ \ and \ 
$\ell\in\{0,\ldots,[nt_2]+r_2-1\}$  \
\begin{equation*}
g_{1,n}+k+g_{2,n}+\ell =[ns_1]-[ns_2]+[nt_1]-[nt_2]+q_1-r_1+q_2-
  r_2+k+\ell \geq (s_1-s_2+t_1-t_2)n/2
\end{equation*}
holds, so Lemma \ref{expbound} yields
\begin{equation*}
{\mathrm
  P}\big(S_{g_{1,n}+k,g_{2,n}+\ell}^{(\alpha,1-\beta)}\geq
(g_{1,n}+k+g_{2,n}+\ell)(\theta_n+\omega/2)\big)\leq \exp \big(
-n(s_1-s_2+t_1-t_2)I_{\theta_n}(\theta_n+\omega/2)/2 \big).
\end{equation*}
Since \ $\omega>0$ \ implies \ $I_{\theta_n}(\theta_n+\omega/2)>0$, \
with the help of the above inequality we obviously obtain
\begin{equation}
   \label{eq:eq2.18}
n^{-1/2}\Cov \big(Y^{(n)}_{q_1,q_2}(s_1,t_1),Y^{(n)}_{r_1,r_2}(s_2,t_2)\big) \to 0
\end{equation}
in exponential rate as \ $n\to\infty$. \ If \
$(1-\alpha)(s_1-s_2)<(1-\beta)(t_1-t_2)$ \ then \ $\omega <0$. \ 
Hence, for sufficiently large \ $n\in{\mathbb N}$ \ we have \
$\omega_n\leq \omega/2 <0$ \ and in this way 
\begin{equation*}
{\mathrm P}\big(S_{g_{1,n}+k,g_{2,n}+\ell}^{(\alpha,1-\beta)}=
g_{1,n}+k\big)\leq {\mathrm
  P}\big(S_{g_{1,n}+k,g_{2,n}+\ell}^{(\alpha,1-\beta)}\leq
(g_{1,n}+k+g_{2,n}+\ell)(\theta_n+\omega/2)\big)
\end{equation*}
for all \ $k\in\{0,\ldots,[ns_2]+r_1-1\}$ \ and \
$\ell\in\{0,\ldots,[nt_2]+r_2-1\}$. \ Using again Lemma \ref{expbound}
we obtain 
\begin{equation*}
{\mathrm
  P}\big(S_{g_{1,n}+k,g_{2,n}+\ell }^{(\alpha,1-\beta)}\leq
(g_{1,n}+k+g_{2,n}+\ell)(\theta_n+\omega/2)\big)\leq \exp \big(
-n(s_1-s_2+t_1-t_2)I_{\theta_n}(\theta_n+\omega/2)/2 \big),
\end{equation*}
which directly implies \eqref{eq:eq2.18}.

Case \ $s_1<s_2, \ t_1<t_2$ \ follows by symmetry. Let \ $s_1<s_2, \
t_1>t_2$, \ so for sufficiently large \ $n\in {\mathbb N}$ \ we have   
 $[ns_1]+q_1\leq [ns_2]+r_1$ \ and \ $[nt_1]+q_2\geq [nt_2]+r_2$ \ and
\begin{equation*}
\Cov
\big(Y^{(n)}_{q_1,q_2}(s_1,t_1),Y^{(n)}_{r_1,r_2}(s_2,t_2)\big)=\!\! 
\sum_{k=0}^{[ns_1]+q_1-1}
\sum_{\ell =0}^{[nt_2]+r_2-1}\!\! {\mathrm
  P}\big(S_{g_{1,n}+k,\ell}^{(\alpha,1-\beta)}=g_{1,n}+k\big){\mathrm
  P}\big(S_{k,g_{2,n}+\ell}^{(\alpha,1-\beta)}=k\big).
\end{equation*}
Now,  let
\begin{alignat*}{2}
\theta_n^{(1)}:=&\,\frac{\alpha(g_{1,n}+k)+(1-\beta)
    \ell}{g_{1,n}+k+\ell} \to \alpha, \qquad
  \theta_n^{(2)}:=&\,\frac{\alpha k+(1-\beta)(g_{2,n}+ \ell)}{k+g_{2,n}+\ell}  
\to 1-\beta , \\
\omega_n^{(1)}:=&\,\frac {g_{1,n}+k}{g_{1,n}+k+\ell}-\theta_n^{(1)} 
\to 1-\alpha, \qquad
\omega_n^{(2)}:=&\,\frac k{k+g_{2,n}+\ell}-\theta_n^{(2)} 
\to -(1-\beta),
\end{alignat*}
as \ $n\to\infty$. \ Using the same ideas as before one can easily
show that for sufficiently large \ $n\in{\mathbb N}$ \ and for all \
$k\in\{0,\ldots,[ns_1]+q_1-1\}$ \ and \
$\ell\in\{0,\ldots,[nt_2]+r_2-1\}$ \ 
\begin{align*}
{\mathrm
  P}\big(S_{g_{1,n}+k,\ell}^{(\alpha,1-\beta)}=g_{1,n}+k\big)&\leq 
\exp
\Big(-n(s_2-s_1)I_{\theta_n^{(1)}}\big(\theta_n^{(1)}+(1-\alpha)/2\big)/2
\Big),\\ 
{\mathrm P}\big(S_{k,g_{2,n}+\ell}^{(\alpha,1-\beta)}=k\big)&\leq
\exp
\Big(-n(t_1-t_2)I_{\theta_n^{(2)}}\big(\theta_n^{(2)}-(1-\beta)/2\big)/2
\Big), 
\end{align*}
implying \eqref{eq:eq2.18}.

Case \ $s_1>s_2, \ t_1<t_2$ \ follows by symmetry, too. \ Finally,
let e.g. \ $s_1=s_2, \ t_1>t_2$, \ so for sufficiently large \
$n\in{\mathbb N}$ 
\begin{align*}
\Cov
\big(Y^{(n)}_{q_1,q_2}(s_1,t_1),Y^{(n)}_{r_1,r_2}(s_2,t_2)\big)=\!\! 
\sum_{k=0}^{[ns_1]+r_1\land q_1-1}
\sum_{\ell =0}^{[nt_2]+r_2-1} &{\mathrm P}\big(S_{q_1-r_1\land
  q_1+k,g_{2,n}+\ell}^{(\alpha,1-\beta)}=q_1-r_1\land q_1+k\big) \\
&\times {\mathrm P}\big(S_{r_1-r_1\land
  q_1+k,\ell}^{(\alpha,1-\beta)}=r_1-r_1\land q_1+k\big).
\end{align*}
Hence, this case and the remaining three cases can be handled in the
same way as the earlier ones.

Now, let \ $(\alpha, \beta, \gamma) \in {\mathcal E}_{1+}$. \ Obviously,
\begin{align*}
\Cov
\big(Y^{(n)}_{q_1,q_2}(s_1,t_1),Y^{(n)}_{r_1,r_2}(s_2,t_2)\big)=\Big (
\big ([ns_1]+q_1\big )\land \big([ns_2]+r_1\big) \Big ) &
\beta ^{\left |[nt_1]-[nt_2]+q_2-r_2\right |} \\
&\times\frac {1-\beta ^{([nt_1]+q_2)\land
    ([nt_2]+r_2)}}{1-\beta^2}, 
\end{align*}
that immediately implies the statement of the Proposition. Case \
$(\alpha, \beta, \gamma) \in {\mathcal E}_{2+}$ \ can be handled in
the same way.  

Finally, in case \ $(\alpha, \beta, \gamma) \in {\mathcal V}_+$ \ the
statement directly 
follows from Lemma \ref{covbound}. \proofend

\section{Proof of Proposition \ref{Bn}}
   \label{sec:sec3}

According to the results of the Introduction in the following sections
we may assume 
\ $\alpha \geq 0$, \ $\beta \geq 0$ \ and
\ $\gamma \geq 0$ \ if \ $\alpha\beta\gamma\geq 0$ \ and  \
$\alpha >0$, \ $\beta >0$ \ and \ $\gamma < 0$ \ if \
$\alpha\beta\gamma <0$. \ In this case \ $\Psi_{\alpha,\beta}$ \
equals the three-by-three matrix of ones denoted by \ ${\mathbf 1}$, \
\begin{equation*}
\Sigma_{\alpha,\beta}=\begin{bmatrix} 1&\alpha\beta &\beta \\
                                   \alpha\beta & 1 & \alpha \\
                                   \beta & \alpha & 1
                   \end{bmatrix},
\end{equation*}
and under the conditions of Proposition \ref{Bn} the coefficients
\ $G(k-i,\ell-j;\alpha,\beta,\gamma)$ \ in representation
\eqref{marep} of \ $X_{k,\ell}$ \ are non-negative.

Obviously,
\begin{align*}
&\frac 1{n^2}{\mathsf E}B_n=\frac 1{n^2}\sum_{(k,\ell)\in R_n}
               {\mathsf E}\left(
                \begin{bmatrix}
                  X_{k-1,\ell} \\
                  X_{k,\ell-1}  \\
                  X_{k-1,\ell-1}
                \end{bmatrix}
                \begin{bmatrix}
                  X_{k-1,\ell} \\
                  X_{k,\ell-1}  \\
                  X_{k-1,\ell-1}
                \end{bmatrix}^{\top}\right)
                \\
&=\!\!\int\limits_0^1\!\!\!\int\limits_0^1\!\!
    \begin{bmatrix}
     \!\Var Y^{(n)}_{0,1}(s,t) & \!\!\!\Cov \big
     (Y^{(n)}_{0,1}(s,t),Y^{(n)}_{1,0}(s,t) \big) & \!\!\!\Cov \big
     (Y^{(n)}_{0,1}(s,t), Y^{(n)}_{0,0}(s,t) \big)\! \\
     \!\Cov \big (Y^{(n)}_{0,1}(s,t), Y^{(n)}_{1,0}(s,t)\big) &
     \!\!\!\Var Y^{(n)}_{1,0}(s,t) & \!\!\!\Cov \big
     (Y^{(n)}_{1,0}(s,t),Y^{(n)}_{0,0}s,t) \big)\! \\
     \!\Cov \big (Y^{(n)}_{0,1}(s,t), Y^{(n)}_{0,0}(s,t)\big) &
     \!\!\!\Cov \big (Y^{(n)}_{1,0}(s,t), Y^{(n)}_{0,0}(s,t)\big) &
     \!\!\!\Var Y^{(n)}_{0,0}(s,t)\!
    \end{bmatrix}\! {\mathrm d}s\,{\mathrm d}t.
\end{align*}
By Lemma \ref{covbound} if \ $(\alpha, \beta, \gamma) \in {\mathcal
  F}_{++}$ \ then   
\begin{equation*}
n^{-1/2}\Big |\Cov \big (Y^{(n)}_{q_1,q_2}(s,t),
  Y^{(n)}_{r_1,r_2}(s,t)\big)\Big |\leq \ 
C_{\alpha,\beta} n^{-1/2}\big(2[ns]+2[nt]+2\big)^{1/2} \leq
C_{\alpha,\beta} (2s+2t+2)^{1/2}, 
 \end{equation*}
where \ $C_{\alpha,\beta}$ \ is a positive constant, while in case \
$(\alpha, \beta, \gamma) \in {\mathcal E}_{1+}\cup 
{\mathcal E}_{2+}$ \ we have
\begin{equation*}
n^{-1}\Big |\Cov \big (Y^{(n)}_{q_1,q_2}(s,t),
  Y^{(n)}_{r_1,r_2}(s,t)\big)\Big |\leq \ 
\frac {C_{\alpha,\beta}}n\frac{[ns]+[nt]}{1-\gamma ^2} \leq
\frac{C_{\alpha,\beta}(s+t)}{1-\gamma^2}, 
 \end{equation*}
$(q_1,q_2),(r_1,r_2)\in \big\{ (0,1),(1,0),(0,0)\big \}$. \ As
both upper bounds are integrable on the unit square \ $[0,1]\times
[0,1]$, \ the dominated convergence theorem applies. Hence, if  \
$(\alpha, \beta, \gamma) \in {\mathcal F}_{++}$ 
 \ by Proposition \ref{covlim} we obtain
\begin{equation}
   \label{eq:eq3.1}
\lim _{n\to\infty}\frac 1{n^{5/2}}{\mathsf E}B_n\!=\!\frac 1{\sqrt{\pi
    (\alpha\!+\!\beta)}(1\!-\!\alpha)(1\!-\!\beta)}
\int\limits_0^1\!\!\int\limits_0^1\! \big((1-\alpha)s\big )^{1/2}\land
\big ((1-\beta)t\big )^{1/2} {\mathrm d}s\,{\mathrm d}t\, {\mathbf
  1}\!=\!\sigma ^2_{\alpha,\beta} {\mathbf 1},
\end{equation}
while in case $(\alpha, \beta, \gamma) \in {\mathcal E}_{1+}\cup
{\mathcal E}_{2+}$ \ we have
\begin{equation}
   \label{eq:eq3.2}
\lim _{n\to\infty}\frac 1{n^3}{\mathsf E}B_n=\frac 1{1-\gamma ^2} 
\int\limits_0^1\!\int\limits_0^1 \big (s{\mathbbm 1}_{\{\alpha=1\}}+
t{\mathbbm 1}_{\{\beta=1\}}\big ) {\mathrm d}s\,{\mathrm d}t \,
\Sigma_{\alpha,\beta} =\frac 1{2(1-\gamma ^2)}\Sigma_{\alpha,\beta},
\end{equation}
where \ ${\mathbbm 1}_H$ \ denotes the indicator of a set \ $H$.

Besides \eqref{eq:eq3.1} and \eqref{eq:eq3.2} to prove the first two
statements of Proposition \ref{Bn} one has to show
\begin{align}
\frac 1{n^{\tau}}&\Var \Big (\sum_{(k,\ell)\in R_n}
X_{k-1+q_1,\ell-1+q_2}X_{k-1+r_1,\ell-1+r_2} \Big)  \label{eq:eq3.3} \\
&=\frac 1{n^{\tau}} \!\!\!\!\sum_{(k_1,\ell_1),(k_2,\ell_2)\in R_n}
\!\!\!\! \!\!\!\!  \!\!\!\!
\Cov \Big ( X_{k_1-1+q_1,\ell_1-1+q_2}X_{k_1-1+r_1,\ell_1-1+r_2},
X_{k_2-1+q_1,\ell_2-1+q_2}X_{k_2-1+r_1,\ell_2-1+r_2} \Big ) 
\nonumber \\
&=\frac
1{n^{\tau-4}}\int\limits_0^1\!\!\!\int\limits_0^1\!\!
\int\limits_0^1\!\!\!\int\limits_0^1\Cov \Big
(Y^{(n)}_{q_1,q_2}(s_1,t_1) Y^{(n)}_{r_1,r_2}(s_1,t_1),
Y^{(n)}_{q_1,q_2}(s_2,t_2) Y^{(n)}_{r_1,r_2}(s_2,t_2)\Big )\,{\mathrm
  d}s_1\,{\mathrm d}t_1\,{\mathrm d}s_2\,{\mathrm d}t_2  \!\to\! 
0\nonumber 
\end{align}
as \ $n\to\infty$, \ where \ $ \{q_1,q_2\},\{r_1,r_2\}\in \big\{
(0,1),(1,0),(0,0)\big \}$ \ and
\begin{equation}
   \label{eq:eq3.4}
\tau :=\begin{cases}
         5,& \quad \text{if \  $(\alpha, \beta, \gamma) \in {\mathcal
             F}_{++}$;} \\ 
         6, & \quad \text{if \ $(\alpha, \beta, \gamma) \in {\mathcal
             E}_{1+}\cup {\mathcal E}_{2+}$.} 
         \end{cases}
\end{equation}
Using Lemma 2.8 of \citet{bpz2} we have 
\begin{align*}
\int\limits_0^1\!\!\!&\int\limits_0^1\!\!
\int\limits_0^1\!\!\!\int\limits_0^1\Cov \Big
(Y^{(n)}_{q_1,q_2}(s_1,t_1) Y^{(n)}_{r_1,r_2}(s_1,t_1),
Y^{(n)}_{q_1,q_2}(s_2,t_2) Y^{(n)}_{r_1,r_2}(s_2,t_2)\Big )\,{\mathrm
  d}s_1\,{\mathrm d}t_1\,{\mathrm d}s_2\,{\mathrm d}t_2 \\
&\leq \!M_4\!\int\limits_0^1\!\!\!\int\limits_0^1\!\!
\int\limits_0^1\!\!\!\int\limits_0^1\Cov \Big
(Y^{(n)}_{q_1,q_2}(s_1,t_1),Y^{(n)}_{q_1,q_2}(s_2,t_2)\Big )\Cov \Big
(Y^{(n)}_{r_1,r_2}(s_1,t_1), Y^{(n)}_{r_1,r_2}(s_2,t_2)\Big )\,{\mathrm
  d}s_1\,{\mathrm d}t_1\,{\mathrm d}s_2\,{\mathrm d}t_2 \\
&+\!M_4\!\int\limits_0^1\!\!\!\int\limits_0^1\!\!
\int\limits_0^1\!\!\!\int\limits_0^1\Cov \Big
(Y^{(n)}_{q_1,q_2}(s_1,t_1), Y^{(n)}_{r_1,r_2}(s_2,t_2)\Big)\Cov \Big(
Y^{(n)}_{r_1,r_2}(s_1,t_1),Y^{(n)}_{q_1,q_2}(s_2,t_2) \Big )\,{\mathrm
  d}s_1\,{\mathrm d}t_1\,{\mathrm d}s_2\,{\mathrm d}t_2,
\end{align*}
which by Lemma \ref{covbound}, Proposition \ref{covlim} and by the
dominated convergence theorem implies \eqref{eq:eq3.3}.

Finally, let 
\begin{alignat}{2}
S_{n,1}&:=\sum_{(k,\ell )\in {R_n}}\big(X_{k,\ell -1}-X_{k-1,\ell
  -1}\big)^2, \qquad 
&&S_{n,3}:=\sum_{(k,\ell )\in {R_n}}\big(X_{k,\ell -1}-X_{k-1,\ell
  -1}\big)X_{k-1,\ell -1}, \nonumber \\
S_{n,2}&:=\sum_{(k,\ell )\in {R_n}}\big(X_{k-1,\ell }-X_{k-1,\ell -1}\big)^2, 
&&S_{n,4}:=\sum_{(k,\ell )\in {R_n}}\big(X_{k-1,\ell }-X_{k-1,\ell
  -1}\big)X_{k-1,\ell -1},  \nonumber\\
S_{n,5}&:=\sum_{(k,\ell )\in {R_n}}\big(X_{k,\ell -1}-X_{k-1,\ell
  -1}\big)\big(X_{k-1,\ell }- &&X_{k-1,\ell-1}\big),  \label{eq:eq3.5} \\
 T_n&:=\sum_{(k,\ell )\in {R_n}}
X_{k-1,\ell -1}^2. && \nonumber
\end{alignat}
Observe, that \ $T_n$ \ is exactly entry (3,3) of the matrix \
$B_n$. \ In  case \ $(\alpha, \beta, \gamma) \in {\mathcal V}_+$
\begin{equation}
   \label{eq:eq3.6}
X_{k,\ell}=\sum_{i=1}^k\sum_{j=1}^{\ell}\varepsilon_{i,j},
\end{equation}
so by the continuous mapping theorem (CMT) \citep{billing} 
\begin{equation}
   \label{eq:eq3.7}
\frac 1{n^4}T_n =
\int\limits_0^1\!\!\int\limits_0^1 \bigg (\frac 1n
Y^{(n)}_{0,0}(s,t) \bigg )^2{\mathrm d}s\,{\mathrm d}t \distr
\int\limits_0^1\!\! \int\limits_0^1{\mathcal W}^2(s,t)\,{\mathrm
  d}s\,{\mathrm d} t  \qquad \text{as \ $n\to\infty$}
\end{equation}
follows from  Donsker's theorem \citep{wichura}
\begin{equation}
   \label{eq:eq3.8}
\frac 1n Y^{(n)}_{0,0}(s,t)=\frac 1n \sum _{i=1}^{[ns]}\sum
_{j=1}^{[nt]}\varepsilon_{i,j} \distr {\mathcal W}(s,t) \qquad
\text{as \ $n\to\infty$.} 
\end{equation} 
Further, by \eqref{eq:eq3.6} we have
\begin{equation}
    \label{eq:eq3.9}
X_{k,\ell}-X_{k,\ell-1}=\sum_{i=1}^k\varepsilon_{i,\ell} \qquad \text{and}\qquad
X_{k,\ell}-X_{k-1,\ell}=\sum_{j=1}^{\ell}\varepsilon_{k,j}.
\end{equation}
Using the independence of the error terms \ $\varepsilon _{i,j}$ \
short calculation shows
\begin{equation*}
{\mathsf E} S_{n,1}={\mathsf E} S_{n,2}=n^2(n-1)/2, \qquad \Var \big
(S_{n,1}\big )=\Var \big (S_{n,2}\big )=O(n^5)
\end{equation*}
implying
\begin{equation}
  \label{eq:eq3.10}
n^{-3}S_{n,1} \qmean 1/2 \quad \text{and} \quad n^{-3}S_{n,2} \qmean
1/2 \qquad \text{as \ $n\to\infty$.}
\end{equation}
Applying again the independence of \ $\varepsilon _{i,j}$ \ it is not
difficult to verify 
\begin{equation*}
{\mathsf E} S_{n,3}={\mathsf E} S_{n,4}=0, \qquad \Var \big
(S_{n,3}\big )=\Var \big (S_{n,4}\big )=O(n^6)
\end{equation*}
and 
\begin{equation*}
{\mathsf E} S_{n,5}=0, \qquad \Var \big
(S_{n,5}\big )=O(n^4).
\end{equation*}
Hence, for all \ $\delta>0$ \ we have
\begin{equation}
  \label{eq:eq3.11}
n^{-3-\delta}S_{n,3} \qmean 0, \quad n^{-3-\delta}S_{n,4} \qmean
0 \quad \text{and} \quad n^{-2-\delta}S_{n,5} \qmean 0 \qquad \text{as
  \ $n\to\infty$.} 
\end{equation}
Obviously,
\begin{alignat*}{2}
\sum_{(k,\ell )\in {R_n}}X_{k-1,\ell }^2&=
S_{n,2}+2S_{n,2}+T_n,  \qquad 
&&\sum_{(k,\ell )\in {R_n}}X_{k-1,\ell }X_{k-1,\ell -1}=S_{n,2}+T_n, \\
\sum_{(k,\ell )\in {R_n}}X_{k,\ell-1 }^2&=S_{n,1}+2S_{n,2}+T_n,
&&\sum_{(k,\ell )\in {R_n}}X_{k,\ell -1 }X_{k-1,\ell -1}=S_{n,1}+T_n,\\ 
&\sum_{(k,\ell )\in {R_n}}X_{k,\ell-1 }X_{k-1,\ell }=S_{n,3}&&+S_{n,4}+S_{n,5}+T_n,
\end{alignat*}
so by \eqref{eq:eq3.10} and \eqref{eq:eq3.11} each entry of \
$n^{-4}B_n$ \ has the same limit in distribution, that completes the
proof. \proofend  

\section{Proof of Proposition \ref{An}}
   \label{sec:sec4}
 
To prove the first two statements of Proposition \ref{An} first we
show that \ $(A_n )_{n\geq 1}$ \ is a square integrable three
dimensional martingale with respect to filtration \ $({\mathcal
  F}_n)_{n\geq 1}$, \ where \ ${\mathcal F}_n$ \ denotes the
$\sigma$-algebra generated by the random variables \
$\{\varepsilon _{k,\ell} : (k,\ell)\in R_n \}$. \ 
In order to do this consider the following  decomposition of \
$A_n-A_{n-1}$, \ where \ $A_0:=0$. \ 
Let \ $A_n^{(i)}, \ i=1,2,3,$ \ denote the components of \ $A_n$. \
By representation \eqref{marep},
\begin{align*}
A_n^{(1)}-A_{n-1}^{(1)}=&\sum_{(k,\ell)\in R_n\setminus R_{n-1}}
\varepsilon_{k,\ell} \sum_{(i,j)\in R_{k-1,\ell}}
G(k-1-i,\ell-j;\alpha,\beta,\gamma)\varepsilon_{i,j}, \\ 
A_n^{(2)}-A_{n-1}^{(2)}=&\sum_{(k,\ell)\in R_n\setminus R_{n-1}}
\varepsilon_{k,\ell} \sum_{(i,j)\in R_{k,\ell-1}}
G(k-i,\ell-1-j;\alpha,\beta,\gamma)\varepsilon_{i,j}, \\
A_n^{(3)}-A_{n-1}^{(3)}=&\sum_{(k,\ell)\in R_n\setminus R_{n-1}}
\varepsilon_{k,\ell} \sum_{(i,j)\in R_{k-1,\ell -1}}
G(k-1-i,\ell-1-j;\alpha,\beta,\gamma)\varepsilon_{i,j}.
\end{align*}
Collecting first the terms containing only \ $\varepsilon_{i,j}$ \ with
 \ $(i,j)\in R_n\setminus R_{n-1}$, \ and then the rest, we obtain the
 decomposition
 \begin{equation}
    \label{eq:eq4.1}
  A_n-A_{n-1}
  =A_{n,1}+\sum_{(k,\ell)\in R_n\setminus R_{n-1}}
            \varepsilon_{k,\ell}A_{n,2,k,\ell},  
 \end{equation}
 where \ $A_{n,1}=\big(A_{n,1}^{(1)},A_{n,1}^{(2)},0\big)^{\top}$ \ and \
$A_{n,2,k,\ell}=\big (\widetilde A_{n,2,k-1,\ell},\widetilde
A_{n,2,k,\ell-1 },\widetilde A_{n,2,k-1,\ell-1 } \big )^{\top}$ \ with 
\begin{align}
A_{n,1}^{(1)}:=&\!\!\!\!\!\!\sum_{(k,\ell)\in R_n\setminus R_{n-1}}\!\!\!\!\!\!
\varepsilon_{k,\ell} \!\!\!\sum_{(i,j)\in R_{k-1,\ell}\setminus
  R_{n-1}}\!\!\!\!\!\!\!\!\!
G(k-1-i,\ell-j;\alpha,\beta,\gamma)\varepsilon_{i,j}=\sum_{k=2}^n\sum_{i=1}^{k-1}
\alpha^{k-1-i}\varepsilon_{i,n}\varepsilon_{k,n}, \label{eq:eq4.2}\\    
A_{n,1}^{(2)}:=&\!\!\!\!\!\!\sum_{(k,\ell)\in R_n\setminus R_{n-1}}\!\!\!\!\!\!
\varepsilon_{k,\ell} \!\!\!\sum_{(i,j)\in R_{k,\ell-1}\setminus
  R_{n-1}} \!\!\!\!\!\!\!\!\!
G(k-i,\ell-1-j;\alpha,\beta,\gamma)\varepsilon_{i,j} 
=\sum_{\ell=2}^n\sum_{j=1}^{\ell-1}  
\beta^{\ell-1-j}\varepsilon_{n,j}\varepsilon_{n,\ell}, \label{eq:eq4.3} \\
\widetilde A_{n,2,k,\ell}:=&\!\!\!\!\!\!\sum_{(i,j)\in R_{k,\ell}\cap
  R_{n-1}} \!\!\!\!\!\!\!\!\!
G(k-i,\ell-j;\alpha,\beta,\gamma)\varepsilon_{i,j}. \label{eq:eq4.4}
\end{align}
The first two components of  \ $A_{n,1}$ \ are quadratic forms of the variables
 \ $\{\varepsilon_{i,j}:(i,j)\in R_n\setminus R_{n-1}\}$, \ hence \ 
$A_{n,1}$ \ is independent of \ ${\mathcal F}_{n-1}$.
\ Besides this the terms \ $\widetilde A_{n,2,k,\ell}$ \ are linear 
combinations of the variables \ $\{\varepsilon_{i,j}:(i,j)\in
R_{n-1}\}$, \ thus \ vectors \ $A_{n,2,k,\ell}$  are measurable with respect 
to \ ${\mathcal F}_{n-1}$. \ Consequently,
\begin{equation*}
  {\mathsf E}(A_n-A_{n-1}\mid{\mathcal F}_{n-1})
  ={\mathsf E} A_{n,1}+\sum_{(k,\ell)\in R_n\setminus R_{n-1}}
     A_{n,2,k,\ell}{\mathsf E}(\varepsilon_{k,\ell}\mid {\mathcal F}_{n-1})
  =0.
\end{equation*}
Hence \ $(A_n)_{n\geq1}$ \ is a square integrable martingale with
 respect to the filtration \ $({\mathcal F}_n)_{n\geq1}$.

By the Martingale Central Limit Theorem \citep{jc}, in order to prove
the first two statements of Proposition \ref{An}, it suffices to show
that the conditional variances of the martingale differences converge
in probability and to verify the conditional Lindeberg condition.
To be precise, the statements are consequences of the following two
propositions. 
\begin{Pro}
   \label{CCAn}
\ If \ $(\alpha, \beta, \gamma) \in {\mathcal F}_{++}$ \ then
\begin{equation*}
\frac 1{n^{5/2}}\sum_{m=1}^n
   {\mathsf E}\big((A_m-A_{m-1})(A_m-A_{m-1})^{\top}\mid
     {\mathcal F}_{m-1}\big)
  \stoch \sigma_{\alpha,\beta }^2{\mathbf 1} \qquad \text{as \
    $n\to\infty$.}
\end{equation*}
If \ $(\alpha, \beta, \gamma) \in {\mathcal E}_{1+}\cup{\mathcal
  E}_{2+}$ \ then
\begin{equation*}
\frac 1{n^3}\sum_{m=1}^n
   {\mathsf E}\big((A_m-A_{m-1})(A_m-A_{m-1})^{\top}\mid
     {\mathcal F}_{m-1}\big)
  \stoch \frac 1{2(1-\gamma ^2)}\Sigma_{\alpha,\beta} \qquad \text{as \
    $n\to\infty$.}
\end{equation*}  
\end{Pro}
\begin{Pro}
   \label{LINDAn} \
For all \ $\delta>0$,
\begin{equation*}
  \frac 1{n^{\tau /2}}
  \sum_{m=1}^n {\mathsf E}\left(\Vert A_m-A_{m-1} \Vert ^2
            \bone_{\left\{\Vert A_m-A_{m-1}\Vert 
    \geq\delta n^{\tau /4 }\right\}} \,\Big|\,{\mathcal F}_{m-1}\right)
  \stoch 0
\end{equation*}
as \ $n\to\infty$,\ where \ $\tau$ \ is defined by \eqref{eq:eq3.4}, i.e.
\begin{equation*}
\tau :=\begin{cases}
         5,& \quad \text{if \ $(\alpha, \beta, \gamma) \in {\mathcal
             F}_{++}$;} \\ 
         6, & \quad \text{if \ $(\alpha, \beta, \gamma) \in {\mathcal
             E}_{1+}\cup{\mathcal E}_{2+}$.} 
         \end{cases}
\end{equation*}
\end{Pro}

\noindent
\textbf{Proof of Proposition \ref{CCAn}.} \
Let \ $U_m:={\mathsf E}\big((A_m-A_{m-1})(A_m-A_{m-1})^{\top} \mid
{\mathcal F}_{m-1}\big)$. \ First we show that if \ $(\alpha, \beta,
\gamma) \in {\mathcal F}_{++}$  
\begin{equation}
   \label {eq:eq4.5}
\lim_{n\to\infty}\frac 1{n^{5/2}} \sum _{m=1}^n{\mathsf E} U_m =
\sigma_{\alpha,\beta}^2{\mathbf 1}, 
\end{equation}
while in case \ $(\alpha, \beta, \gamma) \in {\mathcal E}_{1+}\cup{\mathcal
  E}_{2+}$ \ we have
\begin{equation}
  \label {eq:eq4.6}
\lim_{n\to\infty}\frac 1{n^3}\sum _{m=1}^n{\mathsf E} U_m= \frac 1{2(1-\gamma
  ^2)}\Sigma_{\alpha,\beta}.  
\end{equation}
Obviously, 
\begin{equation*}
A_m-A_{m-1}=\sum_{(k,\ell)\in R_m\setminus
  R_{m-1}}\varepsilon_{k,\ell}
      \begin{bmatrix}
        X_{k-1,\ell} \\
        X_{k,\ell-1}\\
        X_{k-1,\ell-1}
       \end{bmatrix}
\end{equation*}
 and by representation \eqref{marep} and independence of the \
 $\varepsilon_{i,j}$, 
 \ the terms in the summation have zero mean  and they are mutually
 uncorrelated. Since for all \ $ \{q_1,q_2\},\{r_1,r_2\}\in \big\{
(0,1),(1,0),(0,0)\big \}$ \ products \
$X_{k-1+q_1,\ell-1+q_2}X_{k-1+r_1,\ell-1+r_2}$ \ and \
$\varepsilon_{k,\ell}$ \ are independent we obtain
\begin{align}
   \label{eq:eq4.7}
{\mathsf E}U_m =&{\mathsf E}(A_m-A_{m-1})(A_m-A_{m-1})^{\top} \\
=&\!\!\!\!\!\!\sum_{(k,\ell)\in
  R_m\setminus R_{m-1} }\!\!\!\!\!\! {\mathsf E}
              \left(
                \begin{bmatrix}
                  X_{k-1,\ell} \\
                  X_{k,\ell-1}\\
                  X_{k-1,\ell-1}
                \end{bmatrix}
                \begin{bmatrix}
                  X_{k-1,\ell} \\
                  X_{k,\ell-1}\\
                  X_{k-1,\ell-1}
                \end{bmatrix}^{\top}
               \right)
{\mathsf E}\varepsilon _{k,\ell}
={\mathsf E}B_m-{\mathsf E}B_{m-1}, \nonumber
\end{align}
where \ $B_0$ \ equals the three-by-three matrix of
zeros. Consequently, \eqref{eq:eq4.5} and \eqref{eq:eq4.6} follow from
\eqref{eq:eq3.1} and \eqref{eq:eq3.2}, respectively.

By decomposition \eqref{eq:eq4.1} and by the measurability of
 \ $A_{m,2,k,\ell}$ \ with respect to \ ${\mathcal F}_{m-1}$ \ one can derive 
 \begin{align*}
  U_m=& \,{\mathsf E}\big(A_{m,1}A_{m,1}^{\top} \mid {\mathcal F}_{m-1}\big) 
   +\!\!\!\!\!\!\!\!\!\!\!\sum_{(k,\ell)\in R_m\setminus
     R_{m-1}}\!\!\!\!\!\!\!\!\!\!\! 
      {\mathsf E}\big(A_{m,1}\varepsilon_{k,\ell} \mid {\mathcal
        F}_{m-1} \big)A_{m,2,k,\ell}^{\top} 
+\!\!\!\!\!\!\!\!\!\!\!\sum_{(k,\ell)\in R_m\setminus R_{m-1}}\!\!\!\!\!\!\!\!\!\!\!
      A_{m,2,k,\ell}{\mathsf E}\big
      (A_{m,1}^{\top}\varepsilon_{k,\ell} \mid {\mathcal F}_{m-1}\big)\\
  &+\sum_{(k_1,\ell_1)\in R_m\setminus R_{m-1}}
     \sum_{(k_2,\ell_2)\in R_m\setminus R_{m-1}}\!\!\!\!\!\!\!\!
      A_{m,2,k_1,\ell_1}A_{m,2,k_2,\ell_2}^{\top}
      {\mathsf E}\big (\varepsilon_{k_1,\ell_1}
      \varepsilon_{k_2,\ell_2} \mid {\mathcal F}_{m-1}\big).
 \end{align*}
By the independence of \ $A_{m,1}$ \ and
 \ $\{\varepsilon_{k,\ell}:(k,\ell)\in R_m\setminus R_{m-1}\}$ \ from
 \ ${\mathcal F}_{m-1}$, \ and by \ ${\mathsf
   E}(A_{m,1}\varepsilon_{k,\ell})=(0,0,0)^{\top}$, \ one obtains 
 \begin{equation}
   \label{eq:eq4.8}
  U_m={\mathsf E} A_{m,1}A_{m,1}^{\top}+\sum_{(k,\ell)\in R_m\setminus R_{m-1}}
  A_{m,2,k,\ell}A_{m,2,k,\ell}^{\top}.
 \end{equation}
This means that to complete the proof of the proposition we have to
show that for all \ $ \{q_1,q_2\},\{r_1,r_2\}\in \big\{
(0,1),(1,0),(0,0)\big \}$
\begin{equation}
   \label{eq:eq4.9}
\lim _{n\to\infty} \frac 1{n^\tau }\Var \Big
(\sum_{m=1}^n\sum_{(k,\ell)\in R_m\setminus 
R_{m-1}}\widetilde A_{m,2,k-1+q_1,\ell-1+q_2}\widetilde
A_{m,2,k-1+r_1,\ell-1+r_2}\Big )=0,   
\end{equation}
where \ $\tau$ \ is defined by \eqref{eq:eq3.4}. Obviously,
\begin{align*}
\Var &\Big
(\sum_{m=1}^n\sum_{(k,\ell)\in R_m\setminus 
R_{m-1}}\widetilde A_{m,2,k-1+q_1,\ell-1+q_2}\widetilde
A_{m,2,k-1+r_1,\ell-1+r_2}\Big ) \\
&=\sum_{m_1=1}^n\sum_{(k_1,\ell_1)\in R_{m_1}\setminus R_{m_1-1}}
\sum_{m_2=1}^n\sum_{(k_2,\ell_2)\in R_{m_2}\setminus R_{m_2-1}} \!\!\!\!\!\!\!\!
\Cov \Big(\widetilde A_{m_1,2,k_1-1+q_1,\ell_1-1+q_2}\widetilde
A_{m_1,2,k_1-1+r_1,\ell_1-1+r_2}, \\
&\phantom{=\sum_{m_1=1}^n\sum_{(k_1,\ell_1)\in R_{m_1}\setminus R_{m_1-1}}
\sum_{m_2=1}^n\sum_{(k_2,\ell_2)\in R_{m_2}\setminus R_{m_2-1}} \!\!\!\!\!\!\!\!
\Cov \Big(}
\widetilde A_{m_2,2,k_2-1+q_1,\ell_2-1+q_2}\widetilde
A_{m_2,2,k_2-1+r_1,\ell_2-1+r_2} \Big),
\end{align*}
 and using Lemma 2.8 of \citet{bpz2} we have 
\begin{align*}
\Cov \Big(\widetilde A_{m_1,2,k_1-1+q_1,\ell_1-1+q_2}&\widetilde
A_{m_1,2,k_1-1+r_1,\ell_1-1+r_2},\widetilde A_{m_2,2,k_2-1+q_1,\ell_2-1+q_2}\widetilde
A_{m_2,2,k_2-1+r_1,\ell_2-1+r_2}\Big) \\
\leq &M_4 \Cov\Big (\widetilde A_{m_1,2,k_1-1+q_1,\ell_1-1+q_2},
\widetilde A_{m_2,2,k_2-1+q_1,\ell_2-1+q_2}  \Big) \\
&\times\Cov\Big (\widetilde
A_{m_1,2,k_1-1+r_1,\ell_1-1+r_2},\widetilde
A_{m_2,2,k_2-1+r_1,\ell_2-1+r_2} \Big) \\
+ &M_4 \Cov\Big (\widetilde A_{m_1,2,k_1-1+q_1,\ell_1-1+q_2},
\widetilde A_{m_2,2,k_2-1+r_1,\ell_2-1+r_2}  \Big) \\ 
&\times\Cov\Big (\widetilde
A_{m_1,2,k_1-1+r_1,\ell_1-1+r_2},\widetilde
A_{m_2,2,k_2-1+q_1,\ell_2-1+q_2} \Big).
\end{align*}
Moreover, by \eqref{eq:eq4.4} and representation \eqref{marep}
\begin{align*}
\Cov\Big (\widetilde A_{m_1,2,k_1,\ell_1},
\widetilde A_{m_2,2,k_2,\ell_2}  \Big) = \!\!\!\!\!\!\!\!\!\!\!\!\sum _{(i,j)\in
  R_{k_1\land k_2,\ell_1\land \ell_2} \cap R_{m_1\land m_2-1} }
\!\!\!\!\!\!\!\!\!\!\!\!\!\!
G(k_1-i,\ell_1-j;\alpha,\beta,\gamma)G&(k_2-i,\ell_2-j;\alpha,\beta,\gamma)\\
&\leq \Cov\Big (X_{k_1,\ell_1},X_{k_2,\ell_2}\Big).
\end{align*}
Furthermore,
\begin{align*}
&\sum_{m_1=1}^n\sum_{(k_1,\ell_1)\in R_{m_1}\setminus R_{m_1-1}}
\sum_{m_2=1}^n\sum_{(k_2,\ell_2)\in R_{m_2}\setminus R_{m_2-1}}
\Cov\Big (X_{k_1-1+q_1,\ell_1-1+q_2},X_{k_2-1+q_1,\ell_2-1+q_2} \Big)\\
&\phantom{\sum_{m_1=1}^n\sum_{(k_1,\ell_1)\in R_{m_1}\setminus R_{m_1-1}}
\sum_{m_2=1}^n\sum_{(k_2,\ell_2)\in R_{m_2}\setminus R_{m_2-1}}}
\times \Cov\Big (X_{k_1-1+r_1,\ell_1-1+r_2},X_{k_2-1+r_1,\ell_2-1+r_2}
\Big) \\
&=\!\!\!\!\!\!\!\sum_{(k_1,\ell_1),(k_2,\ell_2)\in R_n}\!\!\!\!\!\!\!
\Cov\Big (X_{k_1-1+q_1,\ell_1-1+q_2},X_{k_2-1+q_1,\ell_2-1+q_2} \Big)
\Cov\Big (X_{k_1-1+r_1,\ell_1-1+r_2},X_{k_2-1+r_1,\ell_2-1+r_2} \Big).  
\end{align*}
Hence,
\begin{align*}
&\Var \Big
(\sum_{m=1}^n\sum_{(k,\ell)\in R_m\setminus 
R_{m-1}}\widetilde A_{m,2,k-1+q_1,\ell-1+q_2}\widetilde
A_{m,2,k-1+r_1,\ell-1+r_2}\Big ) \leq M_4 \\ 
&\times\!\!\!\!\!\!\!\!\sum_{(k_1,\ell_1),(k_2,\ell_2)\in R_n}\!\!\!\!
\bigg (\!\Cov\Big (X_{k_1-1+q_1,\ell_1-1+q_2},X_{k_2-1+q_1,\ell_2-1+q_2} \Big)
\Cov\Big (X_{k_1-1+r_1,\ell_1-1+r_2},X_{k_2-1+r_1,\ell_2-1+r_2} \Big)
\\
&\phantom{=====}+\Cov\Big
(X_{k_1-1+q_1,\ell_1-1+q_2},X_{k_2-1+r_1,\ell_2-1+r_2} \Big) 
\Cov\Big (X_{k_1-1+r_1,\ell_1-1+r_2},X_{k_2-1+q_1,\ell_2-1+q_2} \Big)\!\bigg),
\end{align*}
so \eqref{eq:eq4.9} can be proved in a similar way as
\eqref{eq:eq3.3}. \proofend

\medskip
\noindent
\textbf{Proof of Proposition \ref{LINDAn}.}  \ Since 
 \begin{equation*}
  \bone_{\left\{\Vert A_m-A_{m-1}\Vert \geq\delta n^{\tau /4 }\right\}}
   \leq\delta^{-2}n^{-\tau /2 }\Vert A_m-A_{m-1}\Vert ^2,
 \end{equation*}
 to prove Proposition \ref{LINDAn} it suffices to show 
 \begin{equation}
    \label{eq:eq4.10}
  \frac 1{n^{\tau }}\sum_{m=1}^n
         {\mathsf E}\big(\Vert A_m-A_{m-1}\Vert ^4 \mid {\mathcal F}_{m-1}\big)
  \stoch0\qquad\text{as \ $n\to\infty$,}
 \end{equation}
where \ $\tau$ \ is defined by \eqref{eq:eq3.4}. By the
decomposition \eqref{eq:eq4.1} of \ $A_m-A_{m-1}$ \ and by the inequality 
 \ $(x+y)^4\leq2^3(x^4+y^4)$ \ for \ $x,y\in{\mathbb R}$, 
 \begin{equation*}
  \Vert A_m-A_{m-1}\Vert ^4
  \leq2^3\Vert A_{m,1}\Vert ^4
      +2^3\bigg\Vert\sum_{(k,\ell)\in R_m\setminus R_{m-1}}
                 \varepsilon_{k,\ell}A_{m,2,k,\ell}\bigg\Vert^4.
\end{equation*}
By the independence of \ $A_{m,1}$ \ and \ ${\mathcal F}_{m-1}$, \ we have
 \ ${\mathsf E}\big(\Vert A_{m,1}\Vert ^4 \mid {\mathcal F}_{m-1}\big)=
{\mathsf E}\Vert A_{m,1}\Vert ^4$. \ Applying the measurability of \
$A_{m,2,k,\ell}$ \ with respect to 
 \ ${\mathcal F}_{m-1}$, \ we obtain  
\begin{equation*}
  {\mathsf E}\left(\bigg\Vert \sum_{(k,\ell)\in R_m\setminus R_{m-1}}
                  \varepsilon_{k,\ell}A_{m,2,k,\ell}\bigg\Vert ^4
           \,\bigg|\,{\mathcal F}_{m-1}\right)
  \leq\big((M_4-3)^++3\big)
      \left(\sum_{(k,\ell)\in R_m\setminus R_{m-1}}\Vert 
   A_{m,2,k,\ell}\Vert ^2\right)^2.
 \end{equation*}
Hence, in order to prove \eqref{eq:eq4.10}, it suffices to show 
 \begin{align}
  \lim_{n\to\infty}\frac 1{n^{\tau}}\sum_{m=1}^n{\mathsf E} 
   \Vert A_{m,1}\Vert ^4&=0,\label{eq:eq4.11}\\
  \lim_{n\to\infty}\frac 1{n^{\tau}}
   \sum_{m=1}^n
    {\mathsf E}\left(\sum_{(k,\ell)\in R_m\setminus R_{m-1}}
              \Vert A_{m,2,k,\ell}\Vert ^2\right)^2&=0.\label{eq:eq4.12}
 \end{align}
Using \eqref{eq:eq4.2} and \eqref{eq:eq4.3} it is easy to see that
\begin{equation*}
 \Vert A_{m,1}\Vert ^4\leq 2\left(\sum_{k=2}^m\sum_{i=1}^{k-1}
           \alpha ^{k-1-i}\varepsilon_{i,m}\varepsilon_{k,m}\right)^4
               +2\left(\sum_{\ell=2}^m\sum_{j=1}^{\ell-1}
            \beta^{\ell-1-j}\varepsilon_{m,j}\varepsilon_{m,\ell}\right)^4.
\end{equation*}
If \ $0<\alpha,\beta <1$ \ then by Lemma 12 of \citet{bpz1} we have \ 
${\mathsf E} \Vert A_{m,1}\Vert ^4=O(m^2)$, \ while for \ $\alpha=1$ \
or \ $\beta=1$ \ a short calculation shows that \ ${\mathsf E}\Vert
A_{m,1}\Vert ^4=O(m^4)$. \ Hence,
\eqref{eq:eq4.11} is satisfied for both possible values of \
$\tau$. 

Furthermore, 
\begin{align*}
 {\mathsf E}\left(\sum_{(k,\ell)\in R_m\setminus R_{m-1}}
              \Vert A_{m,2,k,\ell}\Vert ^2\right)^2\!\!\!\!
= \!\!\!\!\!\!\!\!\sum_{(k_1,\ell_1)(k_2,\ell_2)\in R_m\setminus
  R_{m-1}}
   \!\!\!\!\!\!\!\!\!\!
     &{\mathsf E} \Big (\big(\widetilde A_{m,2,k_1-1,\ell_1}^2+
\widetilde A_{m,2,k_1,\ell_1-1}^2+\widetilde A_{m,2,k_1-1,\ell_1-1}^2\big) \\
&\times \big(\widetilde A_{m,2,k_2-1,\ell_2}^2+
 \widetilde A_{m,2,k_2,\ell_2-1}^2+\widetilde A_{m,2,k_2-1,\ell_2-1}^2\big)\Big ).
\end{align*}
From Lemma 2.8 of \citet{bpz2} follows
\begin{equation*}
{\mathsf E}\big(\widetilde A_{m,2,k_1,\ell_1}^2\widetilde A_{m,2,k_2,\ell_2}^2\big)
  \leq 3M_4{\mathsf E} \widetilde A_{m,2,k_1,\ell_1}^2
{\mathsf E} \widetilde A_{m,2,k_2,\ell_2}^2,
\end{equation*}
while using \eqref{eq:eq4.4} and representation \eqref{marep} one can
easily see \ ${\mathsf E}\widetilde A_{m,2,k,\ell}^2 \leq \Var\,
X_{k,\ell}$. \ As by Lemma \ref{covbound} there exists a positive constant \
$C_{\alpha,\beta}$ \ such that 
\begin{equation*}
\Var \,X_{k,\ell} \leq\begin{cases}
         C_{\alpha,\beta}\sqrt {k+\ell},& \quad \text{if \ $(\alpha,
           \beta, \gamma) 
           \in {\mathcal F}_{++}$;} \\ 
         C_{\alpha,\beta}(k+\ell), & \quad \text{if \ $(\alpha, \beta,
           \gamma) \in 
           {\mathcal E}_{1+}\cup{\mathcal E}_{2+}$,}  
         \end{cases}
\end{equation*}
short calculation shows
\begin{equation*}
{\mathsf E}\left(\sum_{(k,\ell)\in R_m\setminus R_{m-1}}
              \Vert A_{m,2,k,\ell}\Vert ^2\right)^2=O\big(m^{\tau -2}\big),
\end{equation*}
which implies \eqref{eq:eq4.12}. \proofend

Now, consider the case \ $(\alpha, \beta, \gamma) \in {\mathcal V}_+$. \ Let 
\begin{equation*}
Z_n:=\sum_{(k,\ell )\in {R_n}} X_{k-1,\ell-1}\varepsilon_{k,\ell},
\end{equation*}
and \ $C_n^{(1)}$ \ and \ $C_n^{(2)}$ \ be the random sequences
defined by \eqref{Cn1} and 
\eqref{Cn2}, respectively.
Using equation \eqref{model} which in this case takes form \
$\Delta_1\Delta_2X_{k,\ell}=\varepsilon_{k,\ell}$ \ with
$\Delta_1X_{k,\ell}:=X_{k,\ell}-X_{k-1,\ell}, \
\Delta_2X_{k,\ell}:=X_{k,\ell}-X_{k,\ell-1}$, \ from \eqref{eq:eq3.8}
and CMT we obtain 
\begin{align*}
\frac 1{n^2}Z_n&=\sum_{(k,\ell )\in {R_n}} \bigg(\frac 1n
Y^{(n)}\Big(\frac{k-1}n,\frac{\ell -1}n\Big)\bigg)\Delta_1\Delta_2\bigg(\frac 1n
Y^{(n)}\Big(\frac kn,\frac{\ell}n\Big) \bigg)\\
&=\int\limits_0^1\!\!\int\limits_0^1 \bigg (\frac 1n 
Y^{(n)}_{0,0}(s,t) \bigg ) \bigg (\frac 1n 
Y^{(n)}_{0,0}({\mathrm d}s,{\mathrm d}t)\bigg)
\distr\int\limits_0^1\!\! \int\limits_0^1{\mathcal W}(s,t) {\mathcal
  W}({\mathrm d}s,{\mathrm d} t)  \qquad \text{as \ $n\to\infty$.} \nonumber
\end{align*}
Further, using the independence of the error terms \
$\varepsilon_{i,j}$  \ and \eqref{eq:eq3.9} short calculation shows
\begin{equation*}
{\mathsf E}C_n^{(1)}={\mathsf E}C_n^{(2)}=0, \qquad \Var
\big(C_n^{(1)}\big)=\Var \big(C_n^{(2)}\big)=O(n^3). 
\end{equation*} 
Hence, for all \ $\delta>0$ \ we have
\begin{equation}
   \label{eq:eq4.13}
n^{-3/2-\delta}C_n^{(1)} \qmean 0 \quad \text{and} \quad
n^{-3/2-\delta}C_n^{(2)} \qmean 0 \qquad \text{as \ $n\to\infty$.}
\end{equation}
Obviously, 
\begin{equation*}
A_n-Z_n(1, \ 1, \ 1)^{\top}=\big(C_n^{(1)}, \ C_n^{(2)}, \ 0\big)^{\top},
\end{equation*}
that together with \eqref{eq:eq4.13} completes the proof. \proofend

\section{Proof of Proposition \ref{DetB}}
   \label{sec:sec5}

According to the results of the Introduction it suffices to consider
the case \ $\alpha \geq 0$, \ $\beta \geq 0$ \ and
\ $\gamma \geq 0$ \ if \ $\alpha\beta\gamma\geq 0$ \ and  \
$\alpha >0$, \ $\beta >0$ \ and \ $\gamma < 0$ \ if \
$\alpha\beta\gamma <0$. \

Consider the following expression of \ $\det (B_n)$
\begin{equation*}
\det (B_n)=\sum_{(k_1,\ell _1)\in {R_n}}\sum_{(k_2,\ell _2)\in
  {R_n}}\sum_{(k_3,\ell _3)\in {R_n}} W_{k_1,\ell_1,k_2,\ell_2,k_3,\ell_3},
\end{equation*}
where
\begin{align*}
&W_{k_1,\ell_1,k_2,\ell_2,k_3,\ell_3}:=2X_{k_1-1,\ell_1}X_{k_1-1,\ell_1-1}
  X_{k_2,\ell_2-1}X_{k_2-1,\ell_2-1}X_{k_3-1,\ell_3}X_{k_3,\ell_3-1}
  \\ &\phantom{=}+\! 
X^2_{k_1-1,\ell_1}X^2_{k_2,\ell_2-1}X^2_{k_3-1,\ell_3-1} \!-\!
X^2_{k_1,\ell_1-1}
  X_{k_2-1,\ell_2}X_{k_2-1,\ell_2-1}X_{k_3-1,\ell_3}X_{k_3-1,\ell_3-1}
  \\ &\phantom{=}-\!X^2_{k_1-1,\ell_1}
  X_{k_2,\ell_2-1}X_{k_2-1,\ell_2-1}X_{k_3,\ell_3-1}X_{k_3-1,\ell_3-1}
  \!-\!X^2_{k_1-1,\ell_1-1} X_{k_2,\ell_2-1}X_{k_2-1,\ell_2}X_{k_3,\ell_3-1}X_{k_3,\ell_3-1}.
\end{align*}
Short calculation shows that
\begin{align}
&W_{k_1,\ell_1,k_2,\ell_2,k_3,\ell_3}=\big(
X_{k_1,\ell_1-1}\!-\!X_{k_1-1,\ell_1-1}\big )^2
\big
(X_{k_2-1,\ell_2}\!-\!X_{k_2-1,\ell_2-1}\big)^2X^2_{k_3-1,\ell_3-1} \nonumber
\\ 
&\phantom{=}+2\big (X_{k_1,\ell_1-1}\!-\!X_{k_1-1,\ell_1-1}\big)
\big (X_{k_1-1,\ell_1}\!-\!X_{k_1-1,\ell_1-1}\big)
\big (X_{k_2,\ell_2-1}\!-\!X_{k_2-1,\ell_2-1}\big) \nonumber \\
&\phantom{==}\times \big (X_{k_3-1,\ell_3}\!-\!X_{k_3-1,\ell_3-1}\big)X_{k_2-1,\ell_2-1}
X_{k_3-1,\ell_3-1}  \label{eq:eq5.1}\\ 
&\phantom{=}-\big (X_{k_1,\ell_1-1}\!-\!X_{k_1-1,\ell_1-1}\big)
\big (X_{k_1-1,\ell_1}\!-\!X_{k_1-1,\ell_1-1}\big)
\big (X_{k_2,\ell_2-1}\!-\!X_{k_2-1,\ell_2-1}\big)  \nonumber \\
&\phantom{==}\times \big (X_{k_2-1,\ell_2}\!-\!X_{k_2-1,\ell_2-1}\big)
X^2_{k_3-1,\ell_3-1}  \nonumber \\
&\phantom{=}-\big (X_{k_1,\ell_1-1}\!-\!X_{k_1-1,\ell_1-1}\big)^2
\big (X_{k_2-1,\ell_2}\!-\!X_{k_2-1,\ell_2-1}\big)
\big (X_{k_3-1,\ell_3}\!-\!X_{k_3-1,\ell_3-1}\big)X_{k_2-1,\ell_2-1}
X_{k_3-1,\ell_3-1} \nonumber \\
&\phantom{=}-\big (X_{k_1-1,\ell_1}\!-\!X_{k_1-1,\ell_1-1}\big)^2
\big (X_{k_2,\ell_2-1}\!-\!X_{k_2-1,\ell_2-1}\big)
\big (X_{k_3,\ell_3-1}\!-\!X_{k_3-1,\ell_3-1}\big)X_{k_2-1,\ell_2-1}
X_{k_3-1,\ell_3-1}. \nonumber
\end{align}

First let \ $(\alpha, \beta, \gamma) \in {\mathcal F}_{++}$. \ Using 
notations \eqref{eq:eq3.5} introduced in Section \ref{sec:sec3}, by
representation  \eqref{eq:eq5.1} we have
\begin{align}
   \label{eq:eq5.2}
n^{-13/2}\det (B_n)=&\big (n^{-2}S_{n,1}\big)\big (n^{-2}S_{n,2}\big)\big
(n^{-5/2}\,T_n\big)+2\big
(n^{-2}S_{n,5}\big)\big (n^{-9/4}S_{n,3}\big)\big
(n^{-9/4}S_{n,4}\big)\\
&-\big (n^{-2}S_{n,5}\big)^2\big (n^{-5/2}\,T_n\big)-\big
(n^{-2}S_{n,1}\big)\big (n^{-9/4}S_{n,4}\big)^2-\big
(n^{-2}S_{n,2}\big)\big (n^{-9/4}S_{n,3}\big)^2. \nonumber
\end{align}
Obviously,
\begin{equation*}
{\mathsf E}S_{n,3}=
n^2\int\limits_0^1\!\!\int\limits_0^1 \Cov
\big(Y^{(n)}_{1,0}(s,t),Y^{(n)}_{0,0}(s,t)\big)- \Var
\big(Y^{(n)}_{0,0}(s,t)\big)  {\mathrm d}s\,{\mathrm d}t, 
\end{equation*}
and with the help of \eqref{eq:eq2.1}, Lemmas \ref{binorep} and
\ref{lclt}, Corollary \ref{probdiff}  and Lemma
2.8 of \cite{bpz2} one can show
\begin{align}
   \label{eq:eq5.3}
\Var \big(S_{n,3}\big)\!\leq &n^4\!\!\int\limits_0^1\!\!
\int\limits_0^1\!\!\int\limits_0^1\!\!\int\limits_0^1 
\Big |\Cov \big(Y^{(n)}_{1,0}(s_1,t_1),Y^{(n)}_{0,0}(s_2,t_2)\big)\!-\! 
\Cov \big(Y^{(n)}_{0,0}(s_1,t_1),Y^{(n)}_{0,0}(s_2,t_2)\big)\Big | \\
&\phantom{=}\times \Big |\Cov \big(Y^{(n)}_{1,0}(s_2,t_2),
Y^{(n)}_{0,0}(s_1,t_1)\big)\!-\!  
\Cov \big(Y^{(n)}_{0,0}(s_2,t_2),Y^{(n)}_{0,0}(s_1,t_1)\big)\Big |
{\mathrm d}s_1\,{\mathrm d}t_1\,{\mathrm d}s_2\,{\mathrm d}t_2 \nonumber\\
+&n^4\!\!\int\limits_0^1\!\!
\int\limits_0^1\!\!\int\limits_0^1\!\!\int\limits_0^1 \bigg(\Big |\Cov
\big(Y^{(n)}_{1,0}(s_1,t_1),Y^{(n)}_{0,0}(s_2,t_2)\big)\!-\!  
\Cov \big(Y^{(n)}_{0,0}(s_1,t_1),Y^{(n)}_{0,0}(s_2,t_2)\big)\Big
| \nonumber \\
&\phantom{=}+\Big |\Cov
\big(Y^{(n)}_{1,0}(s_1,t_1),Y^{(n)}_{0,0}(s_2\!+\!1/n,t_2)\big)\!-\! 
\Cov \big(Y^{(n)}_{0,0}(s_1,t_1),Y^{(n)}_{0,0}(s_2\!+\!1/n,t_2)\big)\Big |
\bigg ) \nonumber \\
&\phantom{=}\times \Big |\Cov
\big(Y^{(n)}_{0,0}(s_1,t_1),Y^{(n)}_{0,0}(s_2,t_2)\big)\Big |
{\mathrm d}s_1\,{\mathrm d}t_1\,{\mathrm d}s_2\,{\mathrm d}t_2 
+n^4(M_4-3)^+C_{\alpha,\beta}, \nonumber
\end{align}
where \ $C_{\alpha,\beta}$ \ is a positive constant.
In this way Propositions \ref{covlim} and \ref{covdiff} and the dominated
convergence theorem imply
\begin{equation*}
\lim_{n\to\infty} n^{-9/4}{\mathsf E}S_{n,3} =0 \qquad \text{and} \qquad 
\lim_{n\to\infty} n^{-9/2}\Var \big(S_{n,3}\big) =0,
\end{equation*}
and the same result can be proved for \ $S_{n,4}$. \ Hence,
\begin{equation}
  \label{eq:eq5.4}
n^{-9/4}S_{n,3} \qmean 0 \qquad \text{and} \qquad n^{-9/4}S_{n,4} \qmean
0 \quad\quad \text{as \ $n\to \infty$.}
\end{equation}

Further,
\begin{equation*}
{\mathsf E}S_{n,1}=
n^2\int\limits_0^1\!\!\int\limits_0^1 {\mathsf E}\big
(Y^{(n)}_{1,0}(s,t) - Y^{(n)}_{0,0}(s,t) \big)^2{\mathrm d}s\,{\mathrm
  d}t,
\end{equation*}
and using representation \eqref{marep} and Lemma \ref{binorep} with
notations of Lemma \ref{lclt} we obtain
\begin{align*}
{\mathsf E}\big
(Y^{(n)}_{1,0}(s,t) - Y^{(n)}_{0,0}(s,t) \big)^2=&\frac
{1-\beta^{2[nt]}}{1-\beta^2} \\
+(1-\alpha )^2&\sum_{k=0}^{[ns]-1}\sum_{\ell
  =0}^{[nt]-1}\Big ({\mathsf
  P}\big(S_{k,\ell}^{(\alpha,1-\beta)}=k+1\big)-{\mathsf
  P}\big(S_{k,\ell}^{(\alpha,1-\beta)}=k\big) \Big)^2. 
\end{align*}
Obviously, \ ${\mathsf E}\big (Y^{(n)}_{1,0}(s,t) - Y^{(n)}_{0,0}(s,t)
\big)^2$ \ is a monotone increasing sequence and by Proposition
\ref{covdiff} it has an upper bound independent of \ $s,t$ \ and \
$n$. \ Hence, 
\begin{equation}
   \label{eq:eq5.5}
\lim_{n\to\infty} {\mathsf E}\big (Y^{(n)}_{1,0}(s,t) - Y^{(n)}_{0,0}(s,t)
\big)^2=\frac 1{1-\beta^2}+(1-\alpha )^2\varrho_{\alpha,\beta}^{(1)}>0. 
\end{equation}
Similarly to \eqref{eq:eq5.3} one can show
\begin{align}
  \label{eq:eq5.6}
\Var &\big(S_{n,1}\big)\!\leq n^4(M_4-3)^+C_{\alpha,\beta}\\
&+2n^4\!\!\int\limits_0^1\!\!
\int\limits_0^1\!\!\int\limits_0^1\!\!\int\limits_0^1 \!\!\bigg(
\Big |\Cov
\big(Y^{(n)}_{1,0}(s_1,t_1),Y^{(n)}_{0,0}(s_2\!+\!1/n,t_2)\big)\!-\! 
\Cov
\big(Y^{(n)}_{0,0}(s_1,t_1),Y^{(n)}_{0,0}(s_2\!+\!1/n,t_2)\big)\Big |
\nonumber \\
&\phantom{====}+\Big |\Cov
\big(Y^{(n)}_{1,0}(s_1,t_1),Y^{(n)}_{0,0}(s_2,t_2)\big)\!-\!  
\Cov \big(Y^{(n)}_{0,0}(s_1,t_1),Y^{(n)}_{0,0}(s_2,t_2)\big)\Big |
\bigg )^2 \! 
{\mathrm d}s_1\,{\mathrm d}t_1\,{\mathrm d}s_2\,{\mathrm d}t_2, \nonumber
\end{align}
where \ $C_{\alpha,\beta}$ \ is a positive constant. Again, 
Propositions \ref{covlim} 
and \ref{covdiff}, the dominated 
convergence theorem and \eqref{eq:eq5.5} imply
\begin{equation*}
\lim_{n\to\infty} n^{-2}{\mathsf E}S_{n,1} =\frac
1{1-\beta^2}+(1-\alpha )^2\varrho_{\alpha,\beta}^{(1)} \qquad \text{and} \qquad  
\lim_{n\to\infty} n^{-4}\Var \big(S_{n,1}\big) =0,
\end{equation*}
and a similar result can be proved for \ $S_{n,2}$. \ Hence,
\begin{equation}
  \label{eq:eq5.7}
n^{-2}S_{n,1} \qmean \frac 1{1\!-\!\beta^2}+(1-\alpha
)^2\varrho_{\alpha,\beta}^{(1)}= \kappa_{\alpha,\beta}^{(1)}
\ \ \text{and} \ \  n^{-2}S_{n,2} \qmean 
\frac 1{1\!-\!\alpha^2}+(1-\beta
)^2\varrho_{\beta,\alpha}^{(1)}=\kappa_{\beta,\alpha}^{(1)} 
\end{equation}
as \ $n\to \infty$.

Finally, 
\begin{equation*}
{\mathsf E}S_{n,5}=
n^2\int\limits_0^1\!\!\int\limits_0^1 {\mathsf E}\big
(Y^{(n)}_{1,0}(s,t) - Y^{(n)}_{0,0}(s,t) \big)\big
(Y^{(n)}_{0,1}(s,t) - Y^{(n)}_{0,0}(s,t) \big){\mathrm d}s\,{\mathrm
  d}t,
\end{equation*}
while for \ $\Var \big(S_{n,5}\big)$ \ one can find a result similar
to \eqref{eq:eq5.6}.
Using again representation \eqref{marep} and Lemma \ref{binorep} we
obtain
\begin{equation*}
{\mathsf E}\big (Y^{(n)}_{1,0}(s,t) - Y^{(n)}_{0,0}(s,t) \big)\big
(Y^{(n)}_{0,1}(s,t) - Y^{(n)}_{0,0}(s,t) \big) =
(1-\alpha )(1-\beta){\mathcal V}^{(n)}_{\alpha,\beta}(s,t),
\end{equation*}
where
\begin{align*}
{\mathcal V}^{(n)}_{\alpha,\beta}(s,t):=\sum_{k=0}^{[ns]-1}\sum_{\ell
  =0}^{[nt]-1} & \,\Big ({\mathsf
  P}\big(S_{k,\ell}^{(\alpha,1-\beta)}=k+1\big)-
{\mathsf P}\big(S_{k,\ell}^{(\alpha,1-\beta)}=k\big) \Big) \\
&\times \Big ({\mathsf
  P}\big(S_{\ell,k}^{(\beta,1-\alpha)}=\ell+1\big)-{\mathsf 
    P}\big(S_{\ell,k}^{(\beta,1-\alpha)}=\ell\big) \Big).
\end{align*}
By  Proposition
\ref{covdiff} \ ${\mathsf E}\big (Y^{(n)}_{1,0}(s,t) - Y^{(n)}_{0,0}(s,t) \big)\big
(Y^{(n)}_{0,1}(s,t) - Y^{(n)}_{0,0}(s,t) \big)$ \ is bounded with a
bound independent of \ $s,t$ \ and \ $n$, \ but one
has to show that \ ${\mathcal V}^{(n)}_{\alpha,\beta}(s,t)$ \ has a limit as \
$n\to\infty$. \ In order to prove this we show that for fixed
\ $s$ \ and \ 
$t$ \ values \ ${\mathcal V}^{(n)}_{\alpha,\beta}(s,t)$ \ is a Cauchy sequence.

Let \ $n,m\in {\mathbb N}, \ n>m,\ 0<s,t<1$, \ and without loss of
generality we may assume \ $[ms]\geq 1$ \ and \ $[mt]\geq 1$. \ The
local version of the CLT 
given in Lemma \ref{lcltdiff} yields approximation
\begin{align*}
{\mathcal V}^{(n)}_{\alpha,\beta}(s,t)\!-\! {\mathcal
  V}^{(m)}_{\alpha,\beta}(s,t)\approx 
{\mathcal D}^{(n,m)}_{\alpha,\beta}(s,t)\!:=-\frac
1{2\pi}\!\sum_{k=1}^{[ms]-1} \sum_{\ell=[mt]}^{[nt]-1}\!f(b_{k,\ell},a_{k,\ell}) 
-\frac 1{2\pi}\!\sum_{k=[ms]}^{[ns]-1}
\sum_{\ell=1}^{[nt]-1}\!f(b_{k,\ell},a_{k,\ell}),  
\end{align*}
where now
\begin{equation*}
f(u,v)=\frac {v^2}{u^3}\exp \Big (-\frac {v^2}{2u}\Big),
\end{equation*}
while \ $a_{k,l}$ \ and \ $b_{k,l}$ \ are defined in the same way as
in the proof of Proposition \ref{covdiff}. Using  Lemma \ref{lcltdiff}
one can easily show that for the error of the approximation we have
\begin{align}
   \label{eq:eq5.8}
\Big |&{\mathcal V}^{(n)}_{\alpha,\beta}(s,t)- {\mathcal
  V}^{(m)}_{\alpha,\beta}(s,t)-{\mathcal
  D}^{(n,m)}_{\alpha,\beta}(s,t)\Big|\leq C_{\alpha,\beta}
\bigg(\sum_{k=1}^{[ms]-1} \sum_{\ell=[mt]}^{[nt]-1} \frac
1{b_{k,\ell}^{5/2}} +\sum_{k=[ms]}^{[ns]-1}\sum_{\ell=1}^{[nt]-1}\frac
1{b_{k,\ell}^{5/2}}\\ &+\sum_{\ell=[mt]}^{[nt]-1} \frac
1{b_{0,\ell}^2}+ \sum_{k=[ms]}^{[ns]-1}\frac 1{b_{k,0}^2}\bigg)\leq
\frac {8C_{\alpha,\beta}}{\alpha\beta(1-\alpha)(1-\beta)}\bigg(\frac
1{b_{1,[mt]}^{1/2}} +\frac 1{b_{[ms],1}^{1/2}}+\frac 1{b_{0,[mt]}}
+\frac 1{b_{[ms],0}}\bigg)  \to 0 \nonumber
\end{align}
as \ $m,n\to\infty$, \ where \ $C_{\alpha,\beta}$ \ is a positive
constant. Further, as 
\begin{equation*}
F(u,v)\!:=\!\int\!\! f(u,v)\,{\mathrm d}v =-\frac v{2u^2}\exp \Big(\!-\frac
{v^2}u\Big ) +\frac {\sqrt{\pi}}{4u^{3/2}}\widetilde \Phi\Big(\frac
v{u^{1/2}}\Big ),\quad  \text{so} \quad |F(u,v)|\leq \frac 1{u^{3/2}}, 
\end{equation*}
using the notations of the proof of Proposition \ref{covdiff} we obtain
\begin{align*}
\big|{\mathcal D}^{(n,m)}_{\alpha,\beta}(s,t)\big|&\leq \frac 2{\pi}
\int\limits_1^{[ms]}\int\limits_{[mt]}^{[nt]}f(b_{y,z},a_{y,z})\,{\mathrm
    d}z\,{\mathrm d}y +\frac 2{\pi} \int\limits_{[ms]}^{[ns]}
  \int\limits_1^{[nt]} 
  f(b_{y,z},a_{y,z})\,{\mathrm d}z\,{\mathrm d}y \\
&\leq \frac {2H_{\alpha,\beta}}{\pi b_{1,1}^{-1/2}}\bigg(
\int\limits_{b_{1,[mt]}}^{b_{[ms],[nt]}}\int\limits_{a_{1,[nt]}}^{a_{[ms],[mt]}}
\!\!f(u,v)\,{\mathrm d}v\,{\mathrm d}u +
\int\limits_{b_{[ms],1}}^{b_{[ns],[nt]}}\int\limits_{a_{[ms],[nt]}}^{a_{[ns],1}}
\!\!f(u,v) \,{\mathrm d}v\,{\mathrm d}u 
\bigg) \\
&\leq  \frac {4H_{\alpha,\beta}}{\pi b_{1,1}^{-1/2}}\bigg(
\int\limits_{b_{1,[mt]}}^{b_{[ms],[nt]}}\frac 1{u^{3/2}}\,{\mathrm d}u+
\int\limits_{b_{[ms],1}}^{b_{[ns],[nt]}}\frac 1{u^{3/2}}\,{\mathrm
  d}u\bigg)\leq  \frac {8H_{\alpha,\beta}}{\pi b_{1,1}^{-1/2}}\bigg(
\frac 1{b_{1,[mt]}^{1/2}} +\frac 1{b_{[ms],1}^{1/2}}\bigg) \to 0
\end{align*}
as \ $m,n \to\infty$, \ which together with \eqref{eq:eq5.8} proves
that \ ${\mathcal V}^{(n)}_{\alpha,\beta}(s,t)$ \ is a Cauchy.

In this way 
\begin{equation*}
\lim _{n\to\infty}{\mathsf E}\big (Y^{(n)}_{1,0}(s,t) - Y^{(n)}_{0,0}(s,t) \big)\big
(Y^{(n)}_{0,1}(s,t) - Y^{(n)}_{0,0}(s,t) \big)=\varrho_{\alpha,\beta}^{(2)},
\end{equation*}
so by Propositions \ref{covlim} 
and \ref{covdiff} and the dominated 
convergence theorem we have
\begin{equation*}
\lim_{n\to\infty} n^{-2}{\mathsf E}S_{n,5}
=(1-\alpha)(1-\beta)\varrho_{\alpha,\beta}^{(2)} \qquad \text{and}
\qquad   
\lim_{n\to\infty} n^{-4}\Var \big(S_{n,5}\big) =0.
\end{equation*}
Hence,
\begin{equation}
   \label{eq:eq5.9}
n^{-2}S_{n,5} \qmean (1-\alpha )(1-\beta
)\varrho_{\alpha,\beta}^{(2)}=\kappa_{\alpha,\beta}^{(2)} 
\qquad \text{as \ $n\to \infty$,}
\end{equation}
that together with \eqref{eq:eq5.4}, \eqref{eq:eq5.7} and Proposition
\ref{Bn} implies the first statement of Proposition
\ref{DetB}. Observe, the positivity of the limit of \ $n^{-13/2}\det
(B_n)$ \ follows from the non negativity of \
$\varrho_{\alpha,\beta}^{(1)}\varrho_{\beta,\alpha}^{(1)}- \big
(\varrho_{\alpha,\beta}^{(2)}\big)^2$ \ that is a trivial consequence
of the Cauchy-Schwarz inequality.

Now, let \ $(\alpha, \beta, \gamma) \in {\mathcal E}_{1+}$, \ so
\begin{equation*}
X_{k,\ell}=\sum _{i=1}^k\sum _{j=1}^{\ell}\beta ^{\ell -j}\varepsilon _{i,j}.
\end{equation*}
In this way
\begin{equation*}
X_{k,\ell -1}-X_{k-1,\ell -1}=\sum_{j=1}^{\ell -1}\beta
^{\ell-1-j}\varepsilon_{k,j} \quad \text{and} \quad X_{k-1,\ell}-X_{k-1,\ell-1}
=\sum_{i=1}^{k-1}\varepsilon_{i,\ell} -(1-\beta)X_{k-1,\ell-1}.
\end{equation*}
Using the independence of the error terms \ $\varepsilon _{i,j}$ \
short straightforward calculations show
\begin{alignat*}{2}
{\mathsf E}S_{n,1}&=\frac {n^2}{1-\gamma^2}+\frac
{n(1-\gamma^{2n})}{(1-\gamma^2)^2},  \qquad \qquad &&{\mathsf E}S_{n,3}
={\mathsf E}S_{n,5}=0, \\ 
{\mathsf E}S_{n,2}&=\frac {n^2(n-1)}2+(1+\gamma)^2{\mathsf E}T_n,
 \qquad \qquad &&{\mathsf E}S_{n,4}=-(1+\gamma){\mathsf E}T_n,
\end{alignat*}
and
\begin{alignat*}{3}
\Var \big(S_{n,1}\big)&=O(n^2), \qquad \Var \big(S_{n,3}\big)&=O(n^3),  
\qquad \Var \big(S_{n,5}\big)&=O(n^4),
\\
&\Var \big(S_{n,2}\big)=O(n^5), \quad &\Var \big(S_{n,4}\big)=O(n^5). 
\end{alignat*}
Hence
\begin{equation}
   \label{eq:eq5.10}
n^{-2}S_{n,1} \qmean \frac 1{1-\gamma ^2}, \qquad 
n^{-3}S_{n,2}\qmean \frac 1{1-\gamma}, \qquad
n^{-3}S_{n,4}\qmean - \frac 1{2(1-\gamma)}
\end{equation}
and for all \ $\delta>0$ \
\begin{equation}
    \label{eq:eq5.11}
n^{-3/2-\delta}S_{n,3} \qmean 0, \qquad 
n^{-2-\delta}S_{n,5}\qmean 0
\end{equation}
as \ $n\to\infty$. \ As
by representation 
\eqref{eq:eq5.1} we have
\begin{align*}
n^{-8}\det (B_n)=&\big (n^{-2}S_{n,1}\big)\big (n^{-3}S_{n,2}\big)\big
(n^{-3}\,T_n\big)+2\big
(n^{-5/2}S_{n,5}\big)\big (n^{-5/2}S_{n,3}\big)\big
(n^{-3}S_{n,4}\big)\\
&-\big (n^{-5/2}S_{n,5}\big)^2\big (n^{-3}\,T_n\big)-\big
(n^{-2}S_{n,1}\big)\big (n^{-3}S_{n,4}\big)^2-\big
(n^{-3}S_{n,2}\big)\big (n^{-5/2}S_{n,3}\big)^2, \nonumber
\end{align*}
Proposition \ref{Bn} and limits \eqref{eq:eq5.10} and \eqref{eq:eq5.11}
imply the second statement of Proposition \ref{DetB}. Case \ $(\alpha,
\beta, \gamma) \in {\mathcal E}_{2+}$ \ can be handled in the same
way.   

Finally, consider the case \ $(\alpha, \beta, \gamma) \in {\mathcal V}_+$. \ 
As by representation
\eqref{eq:eq5.1} we have
\begin{align*}
n^{-10}\det (B_n)=&\big (n^{-3}S_{n,1}\big)\big (n^{-3}S_{n,2}\big)\big
(n^{-4}\,T_n\big)+2\big
(n^{-3}S_{n,5}\big)\big (n^{-7/2}S_{n,3}\big)\big
(n^{-7/2}S_{n,4}\big)\\
&-\big (n^{-3}S_{n,5}\big)^2\big (n^{-4}\,T_n\big)-\big
(n^{-3}S_{n,1}\big)\big (n^{-7/2}S_{n,4}\big)^2-\big
(n^{-3}S_{n,2}\big)\big (n^{-7/2}S_{n,3}\big)^2, \nonumber
\end{align*}
the last statement of Proposition \ref{DetB} is a direct consequence of
Slutsky's lemma, \eqref{eq:eq3.7},  \eqref{eq:eq3.10} and
\eqref{eq:eq3.11}. \proofend

\section{Proof of Proposition \ref{Cnlim}}
   \label{sec:sec6}

To prove Proposition \ref{Cnlim} we are going to apply the same ideas
as in the proof of Proposition \ref{An}. Consider first the cases \
$(\alpha, \beta, \gamma) \in {\mathcal F}_{++}$ \
and \ $(\alpha, \beta, \gamma) \in {\mathcal E}_{1+}\cup{\mathcal
  E}_{2+}$. \ As \ $(A_n)_{n\geq1}$ \ is a 
three dimensional 
square integrable martingale with respect to the filtration \
$({\mathcal F}_n)_{n\geq1}$, \ 
random sequence
\begin{equation}
   \label{eq:eq6.1}
C_n-C_{n-1}=C_{n,1}+\sum_{(k,\ell)\in R_n\setminus R_{n-1}}
            \varepsilon_{k,\ell}C_{n,2,k,\ell}
\end{equation}
is a two dimensional martingale difference with respect to the same
filtration,  where 
\begin{equation}
  \label{eq:eq6.2}
C_{n,1}:=\begin{bmatrix}
                 A_{n,1}^{(1)} \\ A_{n,1}^{(2)}        
                 \end{bmatrix}, \quad
C_{n,2,k,\ell}=\begin{bmatrix}
                 C_{n,2,k,\ell}^{(1)} \\ C_{n,2,k,\ell}^{(2)}        
                 \end{bmatrix} := 
                 \begin{bmatrix}
                 \widetilde A_{n,2,k-1,\ell}-\widetilde
                 A_{n,2,k-1,\ell-1} \\ \widetilde A_{n,2,k,\ell
                   -1}-\widetilde A_{n,2,k-1,\ell-1}          
                 \end{bmatrix},
\end{equation}
with \ $A_{n,1}^{(1)}, \ A_{n,1}^{(2)}$ \ and \ $\widetilde
A_{n,2,k,\ell}$ \ defined by \eqref{eq:eq4.2}, \eqref{eq:eq4.3} and
\eqref{eq:eq4.4}, respectively. Here \ $C_{n,1}$ \ is independent of \
${\mathcal F}_{n-1}$, \ while \ $C_{n,2,k,\ell}$ \ is measurable with
respect to it. However, representation \eqref{eq:eq6.1} is also valid
in the case \ $(\alpha, \beta, \gamma) \in {\mathcal V}_+$, \ when
\begin{equation*}
\sum_{(k,\ell)\in R_n\setminus R_{n-1}}
            \varepsilon_{k,\ell}C_{n,2,k,\ell}^{(1)}=\sum_{\ell=1}^{n-1}
            \sum_{i=1}^{n-1} 
            \varepsilon_{i,\ell}\varepsilon_{n,\ell} , \qquad
\sum_{(k,\ell)\in R_n\setminus R_{n-1}}
            \varepsilon_{k,\ell}C_{n,2,k,\ell}^{(2)}=\sum_{k=1}^{n-1}\sum_{j=1}^{n-1}
            \varepsilon_{k,j}\varepsilon_{k,n}. 
\end{equation*}
Hence, \ $C_n-C_{n-1}$ \ is a martingale difference in this case, too.
This means that according to the Martingale Central Limit Theorem the
statement of Proposition \ref{Cnlim} follows from the propositions
below.

\begin{Pro}
   \label{CCCn}
\ If \ $(\alpha, \beta, \gamma) \in {\mathcal F}_{++}$ \ then
\begin{equation*}
\frac 1{n^2}\sum_{m=1}^n
   {\mathsf E}\big((C_m-C_{m-1})(C_m-C_{m-1})^{\top}\mid
     {\mathcal F}_{m-1}\big)
  \stoch {\mathcal K}_{\alpha,\beta}  \qquad \text{as \
    $n\to\infty$.}
\end{equation*}
If \ $(\alpha, \beta, \gamma) \in {\mathcal E}_{1+}$ \ 
\begin{equation*}
\frac 1{n^3}\sum_{m=1}^n
   {\mathsf E}\Big(\big (C_m^{(1)}-C_{m-1}^{(1)}\big)^2\,\big |\,
     {\mathcal F}_{m-1}\Big)
  \stoch \frac 1{1-\gamma },  \ \ \frac 1{n^2}\sum_{m=1}^n
   {\mathsf E}\Big(\big(C_m^{(2)}-C_{m-1}^{(2)}\big)^2 \,\big |\,
     {\mathcal F}_{m-1}\Big)
  \stoch \frac 1{1-\gamma ^2}
\end{equation*} 
 as \ $n\to\infty$.

\noindent
If \ $(\alpha, \beta, \gamma) \in {\mathcal E}_{2+}$ \ then
\begin{equation*}
\frac 1{n^2}\sum_{m=1}^n
   {\mathsf E}\Big(\big (C_m^{(1)}-C_{m-1}^{(1)}\big)^2\,\big |\,
     {\mathcal F}_{m-1}\Big)
  \stoch \frac 1{1-\gamma ^2},  \ \ \frac 1{n^3}\sum_{m=1}^n
   {\mathsf E}\Big(\big(C_m^{(2)}-C_{m-1}^{(2)}\big)^2 \,\big |\,
     {\mathcal F}_{m-1}\Big)
  \stoch \frac 1{1-\gamma}
\end{equation*} 
 as \ $n\to\infty$.

\noindent
If \ $(\alpha, \beta, \gamma) \in {\mathcal V}_+$  \ then
\begin{equation*}
\frac 1{n^3}\sum_{m=1}^n
   {\mathsf E}\big((C_m-C_{m-1})(C_m-C_{m-1})^{\top}\mid
     {\mathcal F}_{m-1}\big)
  \stoch \frac 12\, {\mathcal I}_2 \qquad \text{as \
    $n\to\infty$.}
\end{equation*}
\end{Pro}

\begin{Pro}
   \label{LINDCn} \
If \ $(\alpha, \beta, \gamma) \in {\mathcal F}_{++}$ \ then 
for all \ $\delta>0$
\begin{equation*}
  \frac 1{n^2}
  \sum_{m=1}^n {\mathsf E}\left(\Vert C_m-C_{m-1} \Vert ^2
            \bone_{\left\{\Vert C_m-C_{m-1}\Vert 
    \geq\delta n\right\}} \,\Big|\,{\mathcal F}_{m-1}\right)
  \stoch 0 \qquad \text{as \ $n\to\infty$.}
\end{equation*}
If \ $(\alpha, \beta, \gamma) \in {\mathcal E}_{1+}\cup{\mathcal
  E}_{2+}$ \ then for all \ $\delta>0$
\begin{align*}
  \frac 1{n^3}
  &\sum_{m=1}^n {\mathsf E}\left(\vert C_m^{(i)}-C_{m-1}^{(i)} \vert ^2
            \bone_{\left\{\vert C_m^{(i)}-C_{m-1}^{(i)}\vert 
    \geq\delta n^{3/2}\right\}} \,\Big|\,{\mathcal F}_{m-1}\right)
  \stoch 0\\
  \frac 1{n^2}
  &\sum_{m=1}^n {\mathsf E}\left(\vert C_m^{(j)}-C_{m-1}^{(j)} \vert ^2
            \bone_{\left\{\vert C_m^{(j)}-C_{m-1}^{(j)}\vert 
    \geq\delta n\right\}} \,\Big|\,{\mathcal F}_{m-1}\right)
  \stoch 0
\end{align*}
as \ $n\to\infty$,\ where 
\begin{equation*}
(i,j):=\begin{cases}
       (1,2), & \text{if \ $(\alpha, \beta, \gamma) \in {\mathcal E}_{1+}$;} \\
       (2,1), & \text{if \ $(\alpha, \beta, \gamma) \in {\mathcal E}_{2+}$.}
      \end{cases}
\end{equation*}
If \  $(\alpha, \beta, \gamma) \in {\mathcal V}_+$ \ then 
for all \ $\delta>0$
\begin{equation*}
  \frac 1{n^3}
  \sum_{m=1}^n {\mathsf E}\left(\Vert C_m-C_{m-1} \Vert ^2
            \bone_{\left\{\Vert C_m-C_{m-1}\Vert 
    \geq\delta n^{3/2}\right\}} \,\Big|\,{\mathcal F}_{m-1}\right)
  \stoch 0 \qquad \text{as \ $n\to\infty$.}
\end{equation*}
\end{Pro}

\medskip
\noindent
\textbf{Proof of Proposition \ref{CCCn}.}  \ The proof is very similar
to the proof of Proposition \ref{CCAn}. Let \ $V_m:={\mathsf
  E}\big((C_m-C_{m-1})(C_m-C_{m-1})^{\top} \mid 
{\mathcal F}_{m-1}\big)$. \
First we show that if \ $(\alpha, \beta, \gamma) \in {\mathcal F}_{++}$ 
\begin{equation}
   \label {eq:eq6.3}
\lim_{n\to\infty}\frac 1{n^2} \sum _{m=1}^n{\mathsf E} V_m ={\mathcal
  K}_{\alpha,\beta},  
\end{equation}
if  \ $(\alpha, \beta, \gamma) \in {\mathcal E}_{1+}$
\begin{equation}
  \label {eq:eq6.4}
\lim_{n\to\infty}\frac 1{n^3}\sum _{m=1}^n (1, \ 0)\,{\mathsf E} V_m \,(1, \
0)^{\top} = \frac 1{1-\gamma}, \qquad \lim_{n\to\infty}\frac 1{n^2}\sum
_{m=1}^n(0, \ 1)\,{\mathsf E} V_m \,(0, \ 1)^{\top} =\frac 1{1-\gamma
  ^2},
\end{equation}
if \ $(\alpha, \beta, \gamma) \in {\mathcal E}_{2+}$
\begin{equation}
  \label {eq:eq6.5}
\lim_{n\to\infty}\frac 1{n^2}\sum _{m=1}^n(1, \ 0)\, {\mathsf E}V_m \,(1, \
0)^{\top} = \frac 1{1-\gamma^2}, \quad \lim_{n\to\infty}\frac 1{n^3}\sum
_{m=1}^n(0, \ 1)\,{\mathsf E}V_m \, (0, \ 1)^{\top} =\frac 1{1-\gamma}, 
\end{equation}
while in case \ $(\alpha, \beta, \gamma) \in {\mathcal V}_+$ \ we have
\begin{equation}
   \label{eq:eq6.6}
\lim_{n\to\infty}\frac 1{n^3} \sum _{m=1}^n{\mathsf E} V_m =\frac 12\,{\mathcal
  I}_2,  
\end{equation}

As \ $C_n={\mathcal H}A_n$,  \ we obviously have \
$V_m={\mathcal H}U_m{\mathcal H}^{\top}$, where \ $U_m$ \
is the conditional expectation defined in the proof of Proposition
\ref{CCAn}. Hence, using notations \eqref{eq:eq3.5} introduced in
Section \ref{sec:sec3} and \eqref{eq:eq4.7} we obtain
\begin{equation*}
\sum_{m=1}^n{\mathsf E}V_m={\mathcal H}\,{\mathsf E}B_m{\mathcal
  H}^{\top}={\mathsf E} 
\begin{bmatrix}
S_{n,2} & S_{n,5} \\
S_{n,5} & S_{n,1}
\end{bmatrix}.
\end{equation*}  
In this way \eqref{eq:eq6.3} directly follows from \eqref{eq:eq5.7} and
\eqref{eq:eq5.9},  \eqref{eq:eq6.4} is a consequence of
\eqref{eq:eq5.10} and \eqref{eq:eq5.11},  \eqref{eq:eq6.5} can be
proved in the same way as \eqref{eq:eq6.4}, while \eqref{eq:eq6.6}
follows from \eqref{eq:eq3.10} and \eqref{eq:eq3.11}.

Further, \eqref{eq:eq4.8} implies
\begin{equation*}
  V_m=
 {\mathsf E}C_{m,1}C_{m,1}^{\top}+\sum_{(k,\ell)\in R_m\setminus R_{m-1}} 
  C_{m,2,k,\ell}C_{m,2,k,\ell}^{\top}.
 \end{equation*}
This means that to complete the proof one has to show that 
if \  $(\alpha, \beta, \gamma) \in {\mathcal F}_{++}$ 
\begin{equation}
   \label {eq:eq6.7}
\frac 1{n^4} \Var \Big
(\sum_{m=1}^n\sum_{(k,\ell)\in R_m\setminus 
R_{m-1}}C_{m,2,k,\ell}^{(i)}C_{m,2,k,\ell}^{(j)}\Big )\to 0, \quad
\text{for all \ $i,j=1,2$,}
\end{equation}
if  \ $(\alpha, \beta, \gamma) \in {\mathcal E}_{1+}$
\begin{equation}
  \label {eq:eq6.8}
\frac 1{n^6} \Var \Big
(\sum_{m=1}^n\sum_{(k,\ell)\in R_m\setminus 
R_{m-1}}\!\!\!\! \big (C_{m,2,k,\ell}^{(1)}\big )^2\Big )\to 0, \quad
\frac 1{n^4} \Var \Big
(\sum_{m=1}^n\sum_{(k,\ell)\in R_m\setminus 
R_{m-1}}\!\!\!\! \big (C_{m,2,k,\ell}^{(2)}\big )^2\Big )\to 0,
\end{equation}
if \ $(\alpha, \beta, \gamma) \in {\mathcal E}_{2+}$
\begin{equation}
  \label {eq:eq6.9}
\frac 1{n^4} \Var \Big
(\sum_{m=1}^n\sum_{(k,\ell)\in R_m\setminus 
R_{m-1}}\!\!\!\! \big (C_{m,2,k,\ell}^{(1)}\big )^2\Big ) \to 0, \quad
\frac 1{n^6} \Var \Big
(\sum_{m=1}^n\sum_{(k,\ell)\in R_m\setminus 
R_{m-1}}\!\!\!\! \big (C_{m,2,k,\ell}^{(2)}\big )^2\Big )\to 0,
\end{equation}
while in case \ $(\alpha, \beta, \gamma) \in {\mathcal V}_+$ \ we have
\begin{equation}
   \label{eq:eq6.10}
\frac 1{n^6} \Var \Big
(\sum_{m=1}^n\sum_{(k,\ell)\in R_m\setminus 
R_{m-1}}C_{m,2,k,\ell}^{(i)}C_{m,2,k,\ell}^{(j)}\Big )=0, \quad
\text{for all \ $i,j=1,2$,} 
\end{equation}
as \ $n\to\infty$.

By representation \eqref{marep} of \ $X_{k,\ell}$ \ and  by definition
\eqref{eq:eq4.4} of \ ${\widetilde A}_{m,2,k,\ell}$ \ we have
\begin{alignat*}{2}
{\widetilde A}_{m,2,k-1,m}&=X_{k-1,m}-\sum_{i=1}^{k-1}\alpha
^{k-1-i}\varepsilon _{i,m}, \quad {\widetilde
  A}_{m,2,k-1,m-1}=X_{k-1,m-1}\qquad &&k=1,2,\ldots, m,\\ 
{\widetilde A}_{m,2,m,\ell-1}&=X_{m,\ell-1}-\sum_{j=1}^{\ell-1}\beta
^{\ell-1-j}\varepsilon _{m,j}, \quad
{\widetilde A}_{m,2,m-1,\ell-1}=X_{m-1,\ell-1}, \qquad
&&\ell=1,2,\ldots, m \\
{\widetilde A}_{m,2,m-1,k}&=X_{m-1,k}, \qquad  \qquad \qquad \qquad \
{\widetilde A}_{m,2,k,m-1}=X_{k,m-1}, \qquad && k=1,2,\ldots, m-1, 
\end{alignat*}
so according to \eqref{eq:eq6.2} it is not difficult to see that
\begin{align*}
\sum_{m=1}^n\sum_{(k,\ell)\in R_m\setminus R_{m-1}}\!\!\!\! \big
(C_{m,2,k,\ell}^{(1)}\big )^2&=S_{n,2}+{\mathcal C}^{(1)}_{n,1}-2{\mathcal
  C}^{(1)}_{n,2}, \\
\sum_{m=1}^n\sum_{(k,\ell)\in R_m\setminus R_{m-1}}\!\!\!\! \big
(C_{m,2,k,\ell}^{(2)}\big )^2&=S_{n,1}+{\mathcal
  C}^{(2)}_{n,1}-2{\mathcal C}^{(2)}_{n,2}, \\  
\sum_{m=1}^n\sum_{(k,\ell)\in R_m\setminus R_{m-1}}\!\!\!\! 
C_{m,2,k,\ell}^{(1)}C_{m,2,k,\ell}^{(2)}&=S_{n,5}-{\mathcal
  C}^{(3)}_{n,1}-{\mathcal C}^{(3)}_{n,2}+{\mathcal
  C}^{(3)}_{n,3}, 
\end{align*}
where
\begin{alignat*}{2}
{\mathcal
  C}^{(1)}_{n,1}&:=\!\!\sum_{m=2}^n\sum_{k=2}^m\!\!\bigg(\sum_{i=1}^{k-1}\alpha^{k-1-i}
\varepsilon _{i,m}\bigg)^2, \qquad 
{\mathcal
  C}^{(1)}_{n,2}&:=\!\!\sum_{m=2}^n\sum_{k=2}^m\!\!\bigg(\sum_{i=1}^{k-1}\alpha^{k-1-i}
\varepsilon _{i,m}\bigg)\big(X_{k-1,m}-X_{k-1,m-1}\big),\\
{\mathcal
  C}^{(2)}_{n,1}&:=\!\!\sum_{m=2}^n\sum_{\ell=2}^m\!\!
\bigg(\sum_{j=1}^{\ell-1}\beta^{\ell-1-j}  \varepsilon _{m,j}\bigg)^2, \qquad 
{\mathcal
  C}^{(2)}_{n,2}&:=\!\!\sum_{m=2}^n\sum_{\ell=2}^m\!\!
\bigg(\sum_{j=1}^{\ell-1}\beta^{\ell-1-j} 
\varepsilon _{m,j}\bigg)\big(X_{m,\ell -1}-X_{m-1,\ell -1}\big),
\end{alignat*}
and
\begin{align*}
{\mathcal
  C}^{(3)}_{n,1}&:=\!\!\sum_{m=2}^n\sum_{k=2}^m\!\!
\bigg(\sum_{i=1}^{k-1}\alpha^{k-1-i} 
\varepsilon _{i,m}\bigg)\big(X_{k,m-1}-X_{k-1,m-1}\big),\\
{\mathcal
  C}^{(3)}_{n,2}&:=\!\!\sum_{m=2}^n\sum_{\ell=2}^m\!\!
\bigg(\sum_{j=1}^{\ell-1}\beta^{\ell-1-j} 
\varepsilon _{m,j}\bigg)\big(X_{m-1,\ell}-X_{m-1,\ell -1}\big),\\
{\mathcal
  C}^{(3)}_{n,3}&:=\!\!\sum_{m=2}^n\bigg(\sum_{i=1}^{m-1}\alpha^{m-1-i} 
\varepsilon _{i,m}\bigg)\bigg(\sum_{j=1}^{m-1}\beta^{m-1-j} 
\varepsilon _{m,j}\bigg).
\end{align*}
Using the independence of the error terms \ $\varepsilon_{i,j}$, \
Lemma 2.8 of \citet{bpz2} and Proposition \ref{covdiff},
after long and straightforward calculations we obtain that
if \ $(\alpha, \beta, \gamma) \in {\mathcal F}_{++}$ 
\begin{equation}
  \label{eq:eq6.11}
\Var \Big({\mathcal C}^{(i)}_{n,j}\Big)=O(n^3), \qquad i=1,2,3, \ j=1,2, \quad
\text{and} \quad \Var \Big( {\mathcal C}^{(3)}_{n,3} \Big)=O(n),
\end{equation}
if  \ $(\alpha, \beta, \gamma) \in {\mathcal E}_{1+}$
\begin{equation}
  \label{eq:eq6.12}
\Var \Big( {\mathcal C}^{(1)}_{n,j}\Big)=O(n^5) \quad  \text{and}
\quad  \Var \Big({\mathcal C}^{(2)}_{n,j}\Big)=O(n^3), \qquad  j=1,2,
\end{equation}
if  \ $(\alpha, \beta, \gamma) \in {\mathcal E}_{2+}$
\begin{equation}
  \label{eq:eq6.13}
\Var \Big({\mathcal C}^{(1)}_{n,j}\Big)=O(n^3) \quad  \text{and} \quad 
\Var \Big({\mathcal C}^{(2)}_{n,j}\Big)=O(n^5), \qquad  j=1,2,
\end{equation}
while in case \ $(\alpha, \beta, \gamma) \in {\mathcal V}_+$ \ we have
\begin{equation}
  \label{eq:eq6.14}
\Var \Big({\mathcal C}^{(i)}_{n,j}\Big)=O(n^5), \qquad i=1,2,3, \ j=1,2, \quad
\text{and} \quad \Var \Big({\mathcal C}^{(3)}_{n,3}\Big)=O(n^3).
\end{equation}
In this way, \eqref{eq:eq6.7} -- \eqref{eq:eq6.10} follow from
\eqref{eq:eq3.10}, \eqref{eq:eq3.11},
\eqref{eq:eq5.7}, \eqref{eq:eq5.9}, \eqref{eq:eq5.10} and it's pair
for the case \ $(\alpha, \beta, \gamma) \in {\mathcal E}_{2+}$ \ and
\eqref{eq:eq6.11} -- \eqref{eq:eq6.14}.  \proofend

\medskip
\noindent
\textbf{Proof of Proposition \ref{LINDCn}.} \
Similarly to the proof of Proposition \ref{LINDAn} it suffices to show
that if \ $(\alpha, \beta, \gamma) \in {\mathcal F}_{++}$ \ then  
\begin{equation}
   \label{eq:eq6.15}
  \frac 1{n^4}
  \sum_{m=1}^n {\mathsf E}\left(\Vert C_m-C_{m-1} \Vert ^4
             \,\Big|\,{\mathcal F}_{m-1}\right)
  \stoch 0 \qquad \text{as \ $n\to\infty$,}
\end{equation}
if \ $(\alpha, \beta, \gamma) \in {\mathcal E}_{1+}\cup {\mathcal
  E}_{2+}$ \ then  
\begin{equation}
    \label{eq:eq6.16}
  \frac 1{n^6}
  \sum_{m=1}^n {\mathsf E}\left(\vert C_m^{(i)}-C_{m-1}^{(i)} \vert ^4
             \,\Big|\,{\mathcal F}_{m-1}\right)
  \stoch 0, \quad
  \frac 1{n^4}
  \sum_{m=1}^n {\mathsf E}\left(\vert C_m^{(j)}-C_{m-1}^{(j)} \vert ^4
            \,\Big|\,{\mathcal F}_{m-1}\right)
  \stoch 0
\end{equation}
as \ $n\to\infty$, \ where 
\begin{equation*}
(i,j):=\begin{cases}
       (1,2), & \text{if \ $(\alpha, \beta, \gamma) \in {\mathcal E}_{1+}$;} \\
       (2,1), & \text{if \ $(\alpha, \beta, \gamma) \in {\mathcal E}_{2+}$,}
      \end{cases}
\end{equation*}
while for \  $(\alpha, \beta, \gamma) \in {\mathcal V}_+$ \ 
\begin{equation}
    \label{eq:eq6.17}
  \frac 1{n^6}
  \sum_{m=1}^n {\mathsf E}\left(\Vert C_m-C_{m-1} \Vert ^4
            \,\Big|\,{\mathcal F}_{m-1}\right)
  \stoch 0 \qquad \text{as \ $n\to\infty$.}
\end{equation}
Using decomposition \eqref{eq:eq6.1} we have 
 \begin{equation*}
  \Vert C_m-C_{m-1}\Vert ^4 \leq 2^4\,
      \sum _{i=1}^2\left (\big( A_{m,1}^{(i)}\big )^4+
        \bigg(\sum_{(k,\ell)\in R_m\setminus R_{m-1}} 
                 \varepsilon_{k,\ell}C_{m,2,k,\ell}^{(i)}\bigg)^4
                \right).
\end{equation*}
Further, independence of \ $A_{m,1}$ \ and \ ${\mathcal F}_m$ \
and  measurability of \ $C_{m,2,k,\ell}$ \ with respect to \
${\mathcal F}_m$ \ imply
\begin{equation*}
{\mathsf E}\Big(\big( A_{m,1}^{(i)}\big)^4 \,\Big
|\,{\mathcal F}_{m-1}\Big)={\mathsf E}\big( A_{m,1}^{(i)}\big)^4 
\end{equation*}
and
\begin{equation*}
  {\mathsf E}\left(\bigg ( \sum_{(k,\ell)\in R_m\setminus R_{m-1}}\!\!\!\!\!\!\!\!
                  \varepsilon_{k,\ell}C_{m,2,k,\ell}^{(i)}\bigg) ^4
           \,\bigg|\,{\mathcal F}_{m-1}\right)
  \leq\big((M_4-3)^++3\big)
      \left(\sum_{(k,\ell)\in R_m\setminus R_{m-1}}\!\!\!\!\!\!\!\! \big( 
   C_{m,2,k,\ell}^{(i)}\big)^2\right)^2, \quad i=1,2,
 \end{equation*}
respectively. Hence, to prove \eqref{eq:eq6.15} -- \eqref{eq:eq6.17}
one has to show
\begin{align}
   \label{eq:eq6.18}
\lim_{n\to\infty}\frac 1{n^{\tau_i}}\sum_{m=1}^n{\mathsf E}\big(
A_{m,1}^{(i)}\big)^4 &=0, \\
\lim_{n\to\infty}\frac 1{n^{\tau_i}}\sum_{m=1}^n {\mathsf E}
\left(\sum_{(k,\ell)\in R_m\setminus R_{m-1}}\!\!\!\!\! \big(  
   C_{m,2,k,\ell}^{(i)}\big)^2\right)^2&=0, \quad i=1,2,  \label{eq:eq6.19}
\end{align}
where
\begin{equation*}
\tau_1:=\begin{cases}
       4& \text{if \  $(\alpha, \beta, \gamma) \in {\mathcal
           F}_{++}\cup {\mathcal E}_{2+}$;}
\\
       6,& \text{if \ $(\alpha, \beta, \gamma) \in {\mathcal V}_+\cup
         {\mathcal E}_{1+}$,}  
       \end{cases}  \qquad 
\tau_2:=\begin{cases}
       4& \text{if \  $(\alpha, \beta, \gamma) \in {\mathcal
           F}_{++}\cup {\mathcal E}_{1+}$;}
\\
       6,& \text{if \ $(\alpha, \beta, \gamma) \in {\mathcal V}_+\cup
         {\mathcal E}_{2+}$.}  
       \end{cases}
\end{equation*}
Obviously, using the same arguments as in the case of
\eqref{eq:eq4.11} one can easily see that the limit \eqref{eq:eq6.18}
directly follows from 
\eqref{eq:eq4.2} and \eqref{eq:eq4.3}. Further, similarly to the proof
of \eqref{eq:eq6.7} -- \eqref{eq:eq6.10}, after straightforward
calculations we obtain
\begin{equation}
   \label{eq:eq6.20}
{\mathsf E}
\left(\sum_{(k,\ell)\in R_m\setminus R_{m-1}}\!\!\!\!\! \big(  
   C_{m,2,k,\ell}^{(1)}\big)^2\right)^2 \leq \sum_{(k_1,\ell_1)\in
   R_m\setminus R_{m-1}}\sum_{(k_2,\ell_2)\in
   R_m\setminus R_{m-1}}\!\!\!\!\!{\mathcal G}_{k_1,\ell_1,k_2,\ell_2}
 +{\mathcal R}_m 
\end{equation}
where
\begin{equation*}
{\mathcal G}_{k_1,\ell_1,k_2,\ell_2}:={\mathsf E}(
 X_{k_1-1,\ell_1}-X_{k_1-1,\ell_1-1})^2 (
 X_{k_2-1,\ell_2}-X_{k_2-1,\ell_2-1})^2,
\end{equation*}
and
\begin{equation}
   \label{eq:eq6.21}
\lim_{n\to\infty}\frac 1{n^{\tau_1}}\sum_{m=1}^n{\mathcal R}_m =0.
\end{equation}
Now, in the case \ $(\alpha, \beta, \gamma) \in {\mathcal F}_{++}$ \
with the help of Proposition \ref{covdiff} 
and Lemma 2.8 of \citet{bpz2}, while in the remaining cases by
direct calculations one can show that
\begin{equation*}
{\mathcal G}_{k_1,\ell_1,k_2,\ell_2}\leq 
      \begin{cases}
      C_{\alpha,\beta},& \text{if \ $(\alpha, \beta, \gamma) \in
        {\mathcal F}_{++}\cup {\mathcal E}_{2+} $;}
\\
       C_{\alpha,\beta}(k_1-1)(k_2-1),& \text{if \
        $(\alpha, \beta, \gamma) \in {\mathcal V}_+\cup {\mathcal
          E}_{1+}$,}                                
      \end{cases}
\end{equation*}  
where \ $C_{\alpha,\beta}$ \ is a positive constant, that together
with \eqref{eq:eq6.20} and \eqref{eq:eq6.21} implies 
\eqref{eq:eq6.19}. For \ $i=2$ \ \eqref{eq:eq6.19} can be proved in the
same way. \proofend

\section{Proof of Proposition \ref{LimAn}}
   \label{sec:sec7}

Let \ $(\alpha, \beta, \gamma) \in {\mathcal F}_{++}$. \ Using  
 notations \eqref{eq:eq3.5} introduced in Section \ref{sec:sec3} and
 definitions \eqref{Cn1} and 
\eqref{Cn2} after short calculation we obtain
\begin{equation}
   \label{eq:eq7.1} 
n^{-11/2}{\overline B}_nA_n=\big(n^{-9/2}Q_n^{(1)}\big)\big(n^{-1}C_n\big)+
\big(n^{-17/4}Q_n^{(2)}\big)(0, \ 0, \ 1)\big(n^{-5/4}A_n\big), 
\end{equation} 
where \ $Q_n^{(1)}$ \ is a three-by-two matrix with entries
\begin{align*}
Q_{n,1,1}^{(1)}&:=S_{n,1}T_n-S_{n,3}^2, \qquad \qquad \qquad 
Q_{n,2,2}^{(1)}:=S_{n,2}T_n-S_{n,4}^2,\\
Q_{n,1,2}^{(1)}&=Q_{n,2,1}^{(1)}:=S_{n,3}S_{n,4}-S_{n,5}T_n,\\
Q_{n,3,1}^{(1)}&:=\big(S_{n,3}+S_{n,5}\big)S_{n,3}-\big(S_{n,1}+S_{n,3}\big)S_{n,4}
+\big(S_{n,5}-S_{n,1}\big)T_n ,\\
Q_{n,3,2}^{(1)}&:=\big(S_{n,4}+S_{n,5}\big)S_{n,4}-\big(S_{n,2}+S_{n,4}\big)S_{n,3}
+\big(S_{n,5}-S_{n,2}\big)T_n, 
\end{align*} 
and \ $Q_n^{(2)}=\big(Q_{n,1}^{(2)} \ Q_{n,2}^{(2)} \ Q_{n,3}^{(2)}
\big)^{\top}$ \ with
\begin{align*}
Q_{n,1}^{(2)}&:=S_{n,3}S_{n,5}-S_{n,4}S_{n,1}, \qquad \qquad
Q_{n,2}^{(2)}:=S_{n,4}S_{n,5}-S_{n,3}S_{n,2},\\
Q_{n,3}^{(2)}&:=S_{n,1}S_{n,2}-S_{n,5}^2+\big(S_{n,2}-S_{n,5}\big)S_{n,3}+
\big(S_{n,1}-S_{n,5}\big)S_{n,4}.  
\end{align*}
Now, Proposition \ref{Bn} and limits \eqref{eq:eq5.4},
\eqref{eq:eq5.7} and \eqref{eq:eq5.9} imply
\begin{align}
n^{-9/2}Q_{n,3,1}^{(1)}=&\,\big(n^{-9/4}S_{n,3}+n^{-9/4}S_{n,5}\big)\big(n^{-9/4}S_{n,3}
\big)-\big(n^{-9/4}S_{n,1}+n^{-9/4}S_{n,3}\big)\big(n^{-9/4}S_{n,4}\big)
\nonumber \\
&+\big(n^{-2}S_{n,5}-n^{-2}S_{n,1}\big)\big(n^{-5/2}T_n\big)\stoch
\sigma^2_{\alpha,\beta}\big(\kappa^{(2)}_{\alpha,\beta}-\kappa^{(1)}_{\alpha,\beta}\big),
\label{eq:eq7.2} \\
n^{-17/4}Q_{n,2}^{(2)}=&\,\big(n^{-17/4}S_{n,1}\big)\big(n^{-17/8}S_{n,2}\big)-
\big(n^{-17/8}S_{n,5}\big)^2+\big(n^{-2}S_{n,2}-n^{-2}S_{n,5}\big)
\big(n^{-9/4}S_{n,3}\big) \nonumber \\
&+\big(n^{-2}S_{n,1}-n^{-2}S_{n,5}\big)
\big(n^{-9/4}S_{n,4}\big) \stoch 0 \nonumber
\end{align}
as \ $n\to\infty$, \ and using the same ideas one can find the limits of the
remaining entries of \ $Q_n^{(1)}$ \ and coordinates of \
$Q_n^{(2)}$. \ In this way
\begin{equation*}
n^{-9/2}Q_n^{(1)} \stoch 
\sigma^2_{\alpha,\beta}{\mathcal H}^{\top}{\overline{\mathcal K}}_{\alpha,\beta}
\qquad \text{and} \qquad n^{-17/4}Q_n^{(2)} \stoch \big (0, \ 0, \ 0\big)^{\top}
\end{equation*}
as \ $n\to\infty$, \ that together  with \eqref{eq:eq7.1}, Slutsky's lemma and
Propositions \ref{An} and \ref{Cnlim}  implies the first statement of
Proposition \ref{LimAn}.

Further, let \ $(\alpha, \beta, \gamma) \in {\mathcal E}_{1+}$. \ As in this
case \ ${\overline \Sigma}_{\alpha,\beta}=(1-\gamma^2)(0, \ 1, \ -1
)^{\top} (0, \ 1, \ -1 )$, \ short calculation shows
\begin{align}
n^{-7}{\overline B}_nA_n=&\,\Big(n^{-7}{\overline B}_nA-
\big(2(1\!-\!\gamma^2)\big)^{-2}{\overline \Sigma}_{\alpha,\beta}\Big)
\big(n^{-1}A_n\big)+ \Big(\big(2(1\!-\!\gamma^2)\big)^{-2}{\overline
  \Sigma}_{\alpha,\beta}\Big) \big(n^{-1}A_n\big) \nonumber \\
=&\,\big(n^{-11/2}Q_n^{(1)}\big)\big(n^{-3/2}A_n\big)+\Big(n^{-6}Q_n^{(2)}
-\big(4(1\!-\!\gamma^2)\big)^{-1}(0, \ 1, \ -1)^{\top}
\Big)\big(n^{-1}C_n^{(2)}\big) \label{eq:eq7.3}\\
&+\big(4(1\!-\!\gamma^2)\big)^{-1}(0, \ 1, \
-1)^{\top} \big(n^{-1}C_n^{(2)}\big),  \nonumber
\end{align}
where now \ $Q_n^{(1)}$ \ is a three-by-three matrix with entries
\begin{align*}
Q_{n,1,1}^{(1)}&:=S_{n,1}T_n-S_{n,3}^2, \qquad \qquad \qquad 
Q_{n,2,2}^{(1)}=Q_{n,3,3}^{(1)}:=0,\\
Q_{n,1,2}^{(1)}&=Q_{n,2,1}^{(1)}:=S_{n,3}S_{n,4}-S_{n,5}T_n,\\
Q_{n,1,3}^{(1)}&=Q_{n,3,1}^{(1)}:=\big(S_{n,3}+S_{n,5}\big)S_{n,3}
-\big(S_{n,1}+S_{n,3}\big)S_{n,4} +\big(S_{n,5}-S_{n,1}\big)T_n ,\\
Q_{n,2,3}^{(1)}&:=S_{n,4}S_{n,5}-\big(S_{n,2}+S_{n,4}\big)S_{n,3}+S_{n,5}T_n,
\\
Q_{n,3,2}^{(1)}&:=S_{n,1}S_{n,2}+\big(S_{n,2}-S_{n,5}\big)S_{n,3}
+\big(S_{n,1}-S_{n,5}\big)S_{n,4}+\big(S_{n,1}-S_{n,5}\big)T_n
 \\
&\phantom{:=}-\big(S_{n,3}+S_{n,5}\big)S_{n,3}
+\big(S_{n,1}+S_{n,3}\big)S_{n,4}-S_{n,5}^2,  
\end{align*} 
and \ $Q_n^{(2)}=\big(0 \ Q_{n,2}^{(2)} \ -Q_{n,3}^{(2)} \big)^{\top}$ \ with
\begin{align*}
Q_{n,2}^{(2)}:=&S_{n,2}T_n-S_{n,4}^2, \\
Q_{n,3}^{(2)}:=&S_{n,1}S_{n,2}+\big(S_{n,2}-S_{n,5}\big)S_{n,3}
+\big(S_{n,1}-S_{n,5}\big)S_{n,4}+\big(S_{n,1}-S_{n,5}\big)T_n
+\big(S_{n,2}-S_{n,5}\big)T_n  
 \\
&-\big(S_{n,4}+S_{n,5}\big)S_{n,4}-\big(S_{n,3}+S_{n,5}\big)S_{n,3}
+\big(S_{n,1}+S_{n,3}\big)S_{n,4}+\big(S_{n,2}+S_{n,4}\big)S_{n,3}
-S_{n,5}^2.     
\end{align*}
Using Proposition \ref{Bn} and limits \eqref{eq:eq5.10} and
\eqref{eq:eq5.11}, similarly to \eqref{eq:eq7.2} one can show that
\begin{equation*}
n^{-11/2}Q_{n,i,j}^{(1)} \stoch 0, \qquad  i,j=1,2,3, 
\end{equation*}
and
\begin{equation*}
n^{-6}Q_{n,2}^{(2)} \stoch \big(4(1-\gamma^2)\big)^{-1}, \qquad
n^{-6}Q_{n,3}^{(2)} \stoch \big(4(1-\gamma^2)\big)^{-1} 
\end{equation*}
as \ $n\to\infty$, \ that together  with \eqref{eq:eq7.3}, Slutsky's lemma and
Propositions \ref{An} and \ref{Cnlim} implies the second statement of
Proposition \ref{LimAn}.

Finally, if \ $(\alpha, \beta, \gamma) \in {\mathcal E}_{2+}$ \ we have
${\overline \Sigma}_{\alpha,\beta}=(1-\gamma^2)(1, \ 0, \ -1
)^{\top} (1, \ 0, \ -1 )$. \ Hence, similarly to the previous case one
can prove
that the limiting distribution of \ $n^{-7}{\overline B}_nA_n$ \
equals that of \ $\big(4(1-\gamma^2)\big)^{-1}(1, \ 0, \
-1)^{\top} \big(n^{-1}C_n^{(1)}\big)$ \ which completes the proof. \proofend

\section{Proof of Theorem \ref{main}}
   \label{sec:sec8}
 
Cases \ $(\alpha, \beta, \gamma) \in {\mathcal F}_{++}$ \ and \
$(\alpha, \beta, \gamma) \in {\mathcal E}_{1+}\cup {\mathcal E}_{2+}$
\ are direct consequences of Propositions 
\ref{Cnlim} and \ref{LimAn}.

Consider now the case \ $(\alpha, \beta, \gamma) \in {\mathcal V}_+$. \ Using 
notations \eqref{eq:eq3.5} introduced in Section \ref{sec:sec3}  let  
\begin{equation*}
S_n:=\begin{bmatrix}
  S_{n,1}&0\\0&S_{n,2}\\-S_{n,1}&-S_{n,2} \end{bmatrix}. 
\end{equation*}

As by Proposition \ref{Cnlim}
\begin{equation}
  \label{eq:eq8.1}
n^{-3/2}C_n \distr {\mathcal N}\big(0,{\mathcal I_2}/2\big)
\qquad \text{as \ $n\to\infty$,} 
\end{equation}
to prove the last statement of Theorem \ref{main} it suffices to show 
\begin{equation}
  \label{eq:eq8.2}
n^{3/2}\left (B_n^{-1}A_n - \frac {S_nC_n}{S_{n,1}S_{n,2}}\right )
    \stoch \big (0, \ 0, \ 0 \big)^{\top} \qquad \text{as \ $n\to\infty$.}  
\end{equation}

Further, denote by \ $R_{n,i,j}, \ i,j=1,2,3,$ \ the entries of the matrix \
${\overline B}_n-T_nS_nH$.  \ Short calculations show that
\begin{align*}
R_{n,1,1}&:=-S_{n,3}^2, \qquad R_{n,2,2}:=-S_{n,4}^2, \qquad
R_{n,1,2}=R_{n,2,1}:=S_{n,3}S_{n,4}-T_nS_{n,5}, \\
R_{n,1,3}&=R_{n,3,1}:=S_{n,3}\big (S_{n,3}+S_{n,5}\big)
-S_{n,4}\big(S_{n,1}+S_{n,3}\big)+T_nS_{n,5},\\
R_{n,2,3}&=R_{n,3,2}:=S_{n,4}\big (S_{n,4}+S_{n,5}\big)
-S_{n,3}\big(S_{n,2}+S_{n,4}\big)+T_nS_{n,5}, \\
R_{n,3,3}&:=\big(S_{n,2}+S_{n,4}-S_{n,3}-S_{n,5}\big)
\big(S_{n,1}+S_{n,3}-S_{n,4}-S_{n,5}\big) \\
&\phantom{=}+\big(S_{n,3}+S_{n,4}+S_{n,5}\big)
\big(S_{n,1}+S_{n,2}-2S_{n,5}\big)-2T_n S_{n,5}.
\end{align*}
Now, Slutsky's lemma together with \eqref{eq:eq3.7}, \eqref{eq:eq3.10} and
\eqref{eq:eq3.11} implies
\begin{align*}
n^{-13/2}&R_{n,3,3}=n^{-13/4}\big(S_{n,2}+S_{n,4}-S_{n,3}-S_{n,5}\big)
n^{-13/2}\big(S_{n,1}+S_{n,3}-S_{n,4}-S_{n,5}\big) \\
&+n^{-13/2}\big(S_{n,3}+S_{n,4}+S_{n,5}\big)
n^{-13/2}\big(S_{n,1}+S_{n,2}-2S_{n,5}\big)-2n^{-4}T_n n^{-5/2}S_{n,5}
\stoch 0 
\end{align*}
as \ $n\to\infty$, \ and obviously the same result can be proved for the
remaining 8 entries of \  ${\overline B}_n-T_nS_nH$. \ Combining this result
with Proposition \ref{An} we obtain 
\begin{equation}
   \label{eq:eq8.3}
n^{-17/2}\big({\overline B}_nA_n-T_nS_nC_n\big) \stoch \big (0, \ 0, \ 0
\big)^{\top} \qquad \text{as \ $n\to\infty$.}  
\end{equation}
Using again Slutsky's lemma together with \eqref{eq:eq3.7},
\eqref{eq:eq3.10} and \eqref{eq:eq3.11} from \eqref{eq:eq5.2} we have
\begin{equation*}
n^{-10} \big(T_nS_{n,1}S_{n,2}-\det (B_n)\big)\stoch 0 \qquad \text{as \
  $n\to\infty$,} 
\end{equation*}
that together with \eqref{eq:eq3.10} and \eqref{eq:eq8.1} gives us
\begin{align}
   \label{eq:eq8.4}
\frac 1{n^{17/2}}\bigg(T_nS_nC_n -\frac {S_nC_n}{S_{n,1}S_{n,2}}\det
(B_n)\bigg)=&\,\frac
{(n^{-3}S_n)(n^{-3/2}C_n)}{(n^{-3}S_{n,1})((n^{-3}S_{n,2})} \\
&\times\frac 1{n^{10}}
\big(T_nS_{n,1}S_{n,2}-\det (B_n)\big) \!\stoch\! \big (0, \ 0, \ 0
\big)^{\top} \nonumber
\end{align}
as  \ $n\to\infty$. \  In this way \eqref{eq:eq8.2} follows from
Proposition \ref{DetB} and limits \eqref{eq:eq8.3} 
and \eqref{eq:eq8.4}. 
\proofend

\bigskip
\noindent
{\bf Acknowledgments.} \ Research has been supported by 
the Hungarian  Scientific Research Fund under Grant No. OTKA
T079128/2009 and partially supported by T\'AMOP
4.2.1./B-09/1/KONV-2010-0007/IK/IT project. 
The project is implemented through the New Hungary Development Plan
co-financed by the European Social Fund, and the European Regional
Development Fund.

\end{document}